# The Seiberg-Witten equations and the Weinstein conjecture

Clifford Henry Taubes[†]

Department of Mathematics
Harvard University
Cambridge MA  02133

chtaubes@math.harvard.edu

Let M denote a compact, oriented 3-dimensional manifold and let a denote a contact 1-form on M; thus a ∧ da is nowhere zero.  This article proves that the vector field that generates the kernel of da has a closed integral curve.

[†]Supported in part by the National Science Foundation

# 1. Introduction

Let M denote a compact, orientable 3-manifold and let a denote a smooth 1-form on M such that $a \wedge da$ is nowhere zero. Such a 1-form is called a contact form. The associated Reeb vector field is the section, v, of TM that generates the kernel of da and pairs with a to give 1. The generalized Weinstein conjecture in dimension three asserts that v has at least one closed integral curve (see [W]). The purpose of this article is to prove this conjecture and somewhat more. To state the result, remark that the kernel of the 1-form a defines an oriented 2-plane subbundle $K^{-1} \subset TM$. Since an oriented 3-manifold is spin, the first Chern class of this two-plane bundle is divisible by 2.

**Theorem 1:** *Fix an element* $e \in H^2(M; \mathbb{Z})$ *that differs from half the first Chern class of* K *by a torsion element. There is a non-empty set of closed integral curves of the Reeb vector field, and a positive integer weight assigned to each orbit in this set such that the resulting formal weighted sum of loops represents the Poincaré dual of* e *in* $H_1(M; \mathbb{Z})$.

Note that various special cases of the Weinstein conjecture have already been established. For example Hofer [Hof1] proved the Weinstein conjecture in the case where $M = S^3$ or where $\pi_2(M) \neq 0$, or where the associated contact plane field, ker(a), is *over twisted*. The most recent results known to the author are those of Etnyre and Gay [G], Colin-Honda [CH] and Abbas-Cielebak-Hofer-[ABH]. See [Hof2], [Hof3] and [Hof4], [Hon] for references to older papers about the Weinstein conjecture in dimension 3.

The proof of Theorem 1 invokes a version of the Seiberg-Witten Floer homology described by Peter Kronheimer and Tom Mrowka [KM]. In so doing, it borrows a strategy from [T1] and [T2] that is used to identify the Seiberg-Witten and Gromov invariants of a compact, 4-dimensional symplectic manifold. This said, note that a sequel to this article is planned to connect the story told here with the 4-dimensional story that is told in [T1] (see [T3]). This planned sequel will identify a version of the Seiberg-Witten Floer homology for a given compact, oriented 3-manifold with a variant of the Eliashberg-Givental-Hofer contact homology [EGH], a variant along the lines of Hutching's embedded contact homology (see [HS]). By the way, the equivalence of the Seiberg-Witten invariants and Gromov invariants was used by Chen [Ch] to prove some special cases of Theorem 1. However, the approach taken here is along very different lines than that taken by Chen.

As remarked above, the proof of Theorem 1 uses ideas from [T1] and [T2] to identify the Seiberg-Witten and Gromov invariants of a compact, symplectic 4-manifold. However, there is one crucial new ingredient to the story told here with no analog in the 4-dimensional story, and this involves the notion of spectral flow. In particular, a proof is given in what follows of an apparently novel estimate for the spectral flow of a family of Dirac operators on a 3-manifold. This spectral flow result, Proposition 5.5, has generalizations that may be of independent interest [T4].

Before turning to the details, there is an acknowledgement due: A immense debt is owed to Tom Mrowka and Peter Kronheimer for generously sharing their encyclopedic knowledge of Seiberg-Witten Floer homology and the like. As should be evident, this article owes much to their work. Moreover, the approach taken here was sparked by some comments of Tom Mrowka. A great debt is also owed to Michael Hutchings for his many sage comments, suggestions and support.

**a) An outline of the proof**

Sections 2-5 supply various parts of the proof; Sections 6 and 7 tie up loose ends from Sections 2 and 3; and Section 8 puts the parts from Sections 2-5 together to complete the story. What follows outlines how the parts from Sections 2-5 are used to prove Theorem 1.

The Seiberg-Witten equations on M are a system of equations for a connection on a complex line bundle and a section of a related $\mathbb{C}^2$ bundle of spinors over M. The spinor solves the Dirac equation with covariant derivative defined by the connection and a conveniently chosen Riemannian metric. Meanwhile, the curvature of the connection is must equal a 2-form that is a quadratic function of the spinor. The strategy taken from [T2] is as follows: Deform the Seiberg-Witten equations on the 3-manifold M by adding a constant multiple of -ida to the curvature equation. The multiplying factor is denoted by r. Consider a sequence of values of r that limit to ∞ and a corresponding sequence of solutions to the resulting equations. Under optimal circumstances, the spinor component of a solution to a large r version of the equations is nearly zero only on a set that closely approximates a closed integral curve of the Reeb vector field. As r → ∞ along the sequence, a subsequence of such sets limits to the desired closed integral curve. A precise definition of 'optimal circumstances' and the corresponding existence theorem for a closed integral curve is stated in Theorem 2.1. The rest of Section 2 provides a quick introduction to the Seiberg-Witten equations.

Theorem 2.1 is, of course, useless without a proof that all large r versions of the equations have solutions. This is where the Seiberg-Witten Floer homology enters. Kronheimer and Mrowka [KM] describe $\mathbb{Z}$-graded versions of this theory with non-zero homology in an infinite set of degrees, a set that is unbounded from below. The cycles for this homology theory are the solutions to various allowed deformations of the Seiberg-Witten equations. In particular, the deformations just described are allowed. Since the homology is non-zero, there are solutions to the deformed equations for all values of r that represent any given fixed, but sufficiently negative degree in the Seiberg-Witten Floer homology. Note that the particular classes from $H_1(M; \mathbb{Z})$ that appear in Theorem 1 arise, in part, from the use here of a $\mathbb{Z}$-graded version of Seiberg-Witten Floer homology. A sequel to this article will explain how the $\mathbb{Z}/p\mathbb{Z}$ graded Seiberg-Witten Floer homologies in [KM] can be used to find other homology classes that are generated by integral curves of v. Salient features of the $\mathbb{Z}$-graded Seiberg-Witten Floer homology from [KM] are presented in Section 3.

Even granted solutions of the deformed equations for all values of r, nothing of consequence can be said if these solutions do not meet the 'optimal circumstances' requirements that are demanded by Theorem 2.1. The problematic requirement involves a certain functional on the space of connections. This function associates to a connection the integral over M of the wedge product of the curvature 2-form with ia, where a is the contact 1-form. This function is denoted by E. Given the unbounded sequence of r values and the corresponding sequence of solutions, consider the sequence of numbers whose n'th element is the value of E on the connection for the n'th solution. This sequence of numbers must be bounded to obtain a closed integral curve limit. Thus, an argument is needed that gives such a bound. Here, things get subtle, for Michael Hutchings has convincingly argued that there are sequences that can represent a Seiberg-Witten Floer homology class of fixed degree for which the corresponding sequence of E values diverges.

The argument given here that guarantees sequences with uniformly bounded E values requires the introduction of another function on the space of connections, this the Chern-Simons functional. This function is denoted here by $\mathfrak{cs}$. Up to a factor of -1, the Chern-Simons functional realizes the goal of defining a number from a connection by integrating over M the wedge of the connection 1-form with its curvature 2-form. A precise definition is given in Section 3.

To see how $\mathfrak{cs}$ enters the story, fix attention on a non-zero Seiberg-Witten Floer homology class. Consider an unbounded sequence of r values and a corresponding sequence of solutions to the deformed Seiberg-Witten equations where each solution is a generator that appears in a cycle representative of the given homology class. Let $\mathcal{C}$ denote the set of all such pairs of sequences for the fixed homology class. For each sequence of pairs from $\mathcal{C}$, define another sequence of numbers as follows: The n'th number in this new sequence is the value of $\mathfrak{cs}$ on the connection from the n'th pair in the given sequence of pairs. Now suppose that the attending sequence of E values diverges for any choice of a sequence of pairs from $\mathcal{C}$. As explained in Section 4, this can happen only if there exists a sequence of pairs from $\mathcal{C}$ whose associated sequence of $\mathfrak{cs}$ values diverges as $\mathcal{O}(r^2)$ as $r \to \infty$.

The argument for this uses a third functional, this a perturbation of $\frac{1}{2}$ (rE - $\mathfrak{cs}$). This perturbed function is denoted by $\mathfrak{a}$. Section 4 describes a 'min-max' procedure that assigns a value of $\mathfrak{a}$ to any large r and any Seiberg-Witten homology class. With the homology class fixed, the resulting function of r is continuous and piecewise differentiable. Properties of its derivative where it is differentiable follow from properties of solutions to the deformed Seiberg-Witten equations. In particular, these properties imply the assertion about the divergence as $\mathcal{O}(r^2)$ of an associated sequence of $\mathfrak{cs}$ values.

A digression is need to explain the shift of focus to the sequence of 𝔠𝔰 values. To start the digression, note that the degree of a given Seiberg-Witten cycle as defined by a deformed version of the equations is determined by the spectral flow for a path of self adjoint operators. This path starts at the Dirac operator defined by a fiducial connection and ends at the operator that gives the formal linearization of the Seiberg-Witten equations at any solution that appears as a generator in the given cycle. If the cycle has degree k, then this spectral flow is k. This understood, here is how 𝔠𝔰 comes in: Section 5 proves that this spectral flow differs from $\frac{1}{8\pi^2}$ 𝔠𝔰 by $o(r^2)$. Thus, if 𝔠𝔰 is $\mathcal{O}(r^2)$, then the degree of the cycle is very large.

With the preceding understood, here is how the proof of Theorem 1 ends: Fix a non-trivial Seiberg-Witten Floer homology class. Use this class to define the set $\mathcal{C}$ of sequence pairs. Suppose that the attending sequence of E values diverges for all sequence pairs that come from $\mathcal{C}$. If this is the case, then there exists a sequence pair from $\mathcal{C}$ whose associated sequence of 𝔠𝔰 values diverges as $\mathcal{O}(r^2)$. As a consequence, the degree of the representative cycle for the given homology class must be $\mathcal{O}(r^2) \pm o(r^2)$ and thus the degree is increasing with r. But this is nonsense because the degree is fixed since the homology class is fixed. To avoid this nonsense, there exists a sequence pair from $\mathcal{C}$ where the attending sequence of E values is bounded. Theorem 2.1 uses this sequence to find the desired closed Reeb orbit.

## 2. The Seiberg-Witten equations

Let M here denote a compact, oriented Riemannian 3-manifold. Fix a $\text{Spin}_{\mathbb{C}}(3)$ structure on M. The latter constitutes an equivalence class of lifts of the orthonormal frame bundle, Fr → M, to a principle, U(2) bundle, F → M. The set of such lifings can be placed in a 1-1 correspondence with $H^2(M; \mathbb{Z})$. With a lift chosen, let $\mathbb{S} = F \times_{U(2)} \mathbb{C}^2$.

The bundle $\mathbb{S}$ inherits from $\mathbb{C}^2$ a canonical hermitian inner product. Choose a hermitian connection on $\det(\mathbb{S}) = F \times_{U(2)} \mathbb{C}$ and the latter with the Levi-Civita connection on Fr give $\mathbb{S}$ a connection that respects the inner product. The associated covariant derivative is denoted here by $\nabla$; it sends $C^\infty(M; \mathbb{S})$ to $C^\infty(M; \mathbb{S} \otimes T^*M)$. There is also a canonical anti-hermitian action of $T^*M$ on $\mathbb{S}$, this being Clifford multiplication. The map from $T^*M$ to $\text{End}(\mathbb{S})$ is denoted in what follows by cl.

Granted the preceding, the Seiberg-Witten equations on M are equations for a pair consisting of a connection on $\det(\mathbb{S})$ and a section, $\psi$, of $\mathbb{S}$. The simplest version of these equations read

$$*F = \psi^\dagger \tau \psi \quad and \quad \hat{c}(\nabla \psi) = 0,$$

(2.1)

where the notation is as follows: First, ∗F denotes the Hodge dual of the curvature 2-form of the chosen connection, and $\psi^\dagger\tau\psi$ denotes the section of $iT^*$ that is the metric dual to the homomorphism $\psi^\dagger cl(\cdot)\psi$: $T^*M \to i\mathbb{R}$. Meanwhile, $\hat{c}$: $\mathbb{S} \otimes T^*M \to \mathbb{S}$ denotes the endomorphism that is induced by cl.

**a) Variants of the Seiberg-Witten equations**

Certain variants of (2.1) play central roles in the discussions that follow. To say more, suppose that a is a smooth, nowhere vanishing vector field on M. In what follows, a is going to be a contact form, but there is no need yet to restrict a ∧ da. As a is nowhere zero, it induces the splitting $\mathbb{S} = E \oplus E'$ into eigenbundles for cl(a). Convention taken here is that cl(a) acts as i|a| on the first factor and as -i|a| on the second. There is a *canonical* $\text{Spin}_\mathbb{C}(3)$ structure determined by a, that where the bundle E is the trivial bundle. Use $\mathbb{S}_I$ to denote the canonical $\text{Spin}_\mathbb{C}$ structure's version of $\mathbb{S}$. The splitting for $\mathbb{S}_I$ is written as $I_\mathbb{C} \oplus K^{-1}$ where $I_\mathbb{C} \to M$ denotes the trivial $\mathbb{C}$-bundle. The bundle K is called the *canonical line bundle*. The bundle $K^{-1}$ is isomorphic to the 2-plane subbundle in TM whose vectors are annihilated by the 1-form M. Note that the specification of a canonical $\text{Spin}_\mathbb{C}$-structure allows one to write the bundle $\mathbb{S}$ for any other $\text{Spin}_\mathbb{C}$-structure as

$$\mathbb{S} = E \oplus K^{-1}E$$

(2.2)

where $E \to M$ is a complex line bundle. Thus, $\det(\mathbb{S}) = K^{-1}E^2$ in general. By the way, assigning E's first Chern class to the given $\text{Spin}_\mathbb{C}$ structure provides a 1-1 correspondence between the set of $\text{Spin}_\mathbb{C}$ structures and $H^2(M; \mathbb{Z})$. Note that the first Chern class of $\det(\mathbb{S})$ is a torsion class if and only if $E^2$ differs from K by a torsion class.

Let $1_\mathbb{C}$ denote a unit normed, trivializing section of $I_\mathbb{C}$. There is a unique connection on $\det(\mathbb{S}_I) = K^{-1}$ with the property that the section $\psi = (1_\mathbb{C}, 0)$ of $\mathbb{S}_I$ is annihilated by the associated Dirac operator. This connection is called the *canonical connection*. When necessary, this connection is denoted by γ. Note that with $\mathbb{S}$ as in (2.1), any given connection on $\det(\mathbb{S})$ can be written as γ + 2A where A is a connection on E. The Dirac operator on $C^\infty(\mathbb{S})$ that is defined by a given connection A on E is denoted below by $D_A$.

Now let a denote a contact form on an orientable 3-manifold M. Orient M so that a ∧ da is a positive 3-form. Fix a Riemannian metric on M so that a has unit length and so that da = 2∗a. Use a to define the canonical $\text{Spin}_\mathbb{C}$ structure, the canonical bundle K, and the canonical connection, γ, on K. Let E denote a given complex line bundle over M.

The model for the variants of the Seiberg-Witten equations that are of concern in what follows is a system of equations for a pair (A, ψ) of connection on E and section ψ of $\mathbb{S} = E \oplus K^{-1}E$. These equations require the specification of a constant, $r \in [0, \infty)$. With r chosen, the equations read:

- $B_A = r(\psi^\dagger \tau \psi - ia) + i\varpi_K$,
- $D_A \psi = 0$.

(2.3)

Here, $B_A$ denotes the metric Hodge dual of the curvature 2-form of the connection A, and $\varpi_K$ is the harmonic 1-form on M with the property that the Hodge dual of $-\frac{1}{\pi} *\varpi_K$ represents the image in $H^2(M; \mathbb{Q})$ of the first Chern class of K.

    The other variants are perturbations of (2.3). To describe the latter, introduce Conn(E) to denote the space of connections on E. The generic sort of perturbed equation is defined with the help of a smooth, gauge invariant function $\mathfrak{g}$: Conn(E) $\times\, C^\infty(\mathbb{S}) \to \mathbb{R}$. The adjective 'gauge invariant' means that $\mathfrak{g}(A - u^{-1}du, u\psi) = \mathfrak{g}(A, \psi)$ for any choice of pair $(A, \psi) \in$ Conn(E) $\times\, C^\infty(\mathbb{S})$ and any smooth map u: $M \to S^1$. The function $\mathfrak{g}$ is chosen so that its gradient defines a smooth section of $C^\infty(M; iT^*M \oplus \mathbb{S})$. This is to say that there is a smooth map, $(\mathfrak{T}, \mathfrak{S})$: Conn(E) $\times\, C^\infty(\mathbb{S}) \to C^\infty(M; iT^*M \oplus \mathbb{S})$, which is such that that $\frac{d}{dt} \mathfrak{g}(A+tb, \psi + t\eta)|_{t=0} = \int_M (b \wedge *\mathfrak{T} - \frac{1}{2}(\eta^\dagger \mathfrak{S} + \mathfrak{S}^\dagger \eta))$ for any pair $(A, \psi)$ in Conn(E) $\times\, C^\infty(M; \mathbb{S})$ and any $(b, \eta) \in C^\infty(M; iT^*M \oplus \mathbb{S})$. Here, $*$ denotes the metric's Hodge dual. The allowed functions $\mathfrak{g}$ are of the sort that are introduced in Chapter 11 of [KM]. In particular, they are *tame* in the sense used by [KM]. The most general perturbation also requires the choice of a harmonic 1-form. This is denoted below by $\varpi$.

    With $\mathfrak{g}$ and $\varpi$ given, the corresponding perturbed equations are:

- $B_A = r(\psi^\dagger \tau \psi - ia) + \mathfrak{T}(A, \psi) + i\varpi$,
- $D_A \psi = \mathfrak{S}(A, \psi)$.

(2.4)

The terms $\mathfrak{T}$, $\mathfrak{S}$ and $\varpi$ are deemed the *perturbation terms*. Except in this subsection, the next subsection and Section 6, the form $\varpi$ is taken to be the form $\varpi_K$ that appears in (2.3).

    Note that the equations in (2.4) are gauge invariant. This means the following: Suppose that $(A, \psi) \in$ Conn(E) $\times\, C^\infty(M; \mathbb{S})$ is a solution to (2.4) and u is a smooth map from M to $S^1$. Then $(A - u^{-1}du, u\psi)$ is also a solution to (2.4). A solution $(A, \psi)$ to (2.4) with $\psi$ not identically zero is called *irreducible*. The stabilizer in $C^\infty(M; S^1)$ of an irreducible solution is the constant map to $1 \in S^1$. That of a solution $(A, \psi = 0)$ consists of the circle of constant maps to $S^1$.

    Note for reference below that there is a special solution to (2.4) when $\varpi = 0$ and $\mathfrak{g}$ is chosen appropriately. To describe the solution, take $E = I_\mathbb{C}$ and so the Spin$_\mathbb{C}$ structure is canonical. Let $A_I$ denote the connection on $I_\mathbb{C}$ for which the section $1_\mathbb{C}$ is covariantly

constant. Define the section $\psi_I$ of $\mathbb{S}$ by writing it as $(1_\mathbb{C}, 0)$ with respect to the splitting $\mathbb{S}_I = I_\mathbb{C} \oplus K^{-1}$. Then $(A_I, \psi_I)$ is a solution to (2.4) if $\mathfrak{g}$ is chosen so that $\mathfrak{T}$ and $\mathfrak{S}$ vanishes on any $(A, \psi)$ that is gauge equivalent to $(A_I, \psi_I)$.

A solution to (2.4) is deemed to be *reducible* if $\psi$ is identically zero. Reducible solutions to (2.4) arise when $*(iB_A + \varpi)$ is zero in $H^2(M; \mathbb{R})$. For example, in the case of (2.5) below where det($\mathbb{S}$) has torsion first Chern class, the reducible solutions have $\psi = 0$ and $A = A_E - \frac{1}{2} i r a + \mu$ where $A_E$ is a connection on E whose curvature 2-form is $*i\varpi_K$.

What follows explains how solutions to certain versions of (2.4) can lead to closed integral curves of the vector field v.

**Theorem 2.1:** *Fix a complex line bundle* E *so as to define a Spin$_\mathbb{C}$-structure on* M *with spinor bundle* $\mathbb{S}$ *given by (2.2). Let* $\{r_n\}_{n=1,2,\ldots}$ *denote an increasing, unbounded sequence of postive real numbers and for each n, let* $\mu_n$ *denote a co-exact 1-form on M and let* $\varpi_n$ *denote a harmonic 1-form. For each n, let* $(A_n, \psi_n) \in \text{Conn}(E) \times C^\infty(M; \mathbb{S})$ *denote a solution to the* $r = r_n$ *version of (2.4) as defined using the perturbation defined by* $\varpi_n$ *and the pair* $(\mathfrak{T}_n = *d\mu_n, \mathfrak{S}_n = 0)$. *Suppose that there is an n-independent bound for the* $C^3$ *norm of* $\mu_n$ *and* $L^2$ *norm of* $\varpi_n$. *In addition, suppose that there exists an index* n *independent upper bound for* $i\int_M a \wedge *B_{A_n}$ *and a strictly positive,* n *independent lower bound for* $\sup_M (1 - |\psi_n|)$. *Then there exists a non-empty set of closed, integral curves of the Reeb vector field. Moreover, there exists a positive, integer weight assigned to each of these integral curves such that the corresponding formal, integer weighted sum of loops in* M *gives the class in* $H_1(M; \mathbb{Z})$ *that is Poincaré dual to the first Chern class of the bundle* E.

This theorem can be proved using results from [T2]. However, most of the heavy analysis in [T2] is not required here given that this theorem concerns dimension 3, not dimension 4. This being the case, an independent proof of Theorem 2.1 is given in Section 6.

Most of this article uses a version of (2.4) with a very simple perturbation term:

- $B_A = r (\psi^\dagger \tau \psi - i a) + i(*d\mu + \varpi_K)$,
- $D_A \psi = 0$.

(2.5)

Here, $\mu$ is a co-closed 1-form that is $L^2$-orthogonal to all harmonic 1-forms. In what follows, a 1-form $\mu$ is said to be *co-exact* when $d*\mu = 0$ and $\mu$ is orthogonal to all harmonic 1-forms.

**b) Apriori bounds**

The proof of Theorem 2.1 and much of the facts about solutions to (2.4) exploit just two fundamental apriori bounds for solutions to the large r versions of (2.4). To state these bounds, use the splitting in (2.2) to write a given section $\psi$ of $\mathbb{S}$ as $\psi = (\alpha, \beta)$ where $\alpha$ is a section of E and where $\beta$ is a section of $K^{-1}E$.

The next lemma supplies the fundamental estimates

**Lemma 2.2**: *Given $c > 0$, there is a constant $\kappa \geq 1$ with the following significance: Let $\mu$ denote a co-exact 1-form and let $\varpi$ denote a harmonic 1-form. Assume that both the $C^3$ norm of $\mu$ and the $L^2$ norm of $\varpi$ bounded by c. Fix $r \geq 1$ and suppose that $(A, \psi = (\alpha, \beta))$ is a solution to the version of (2.4) given by r and the perturbation data $(\mathfrak{T} = *d\mu, \mathfrak{S} = 0)$ and $\varpi$. Then*
- $|\alpha| \leq 1 + \kappa\, r^{-1}$ ,
- $|\beta|^2 \leq \kappa\, r^{-1}(|1 - |\alpha|^2| + r^{-1})$ .

This lemma is proved in Section 6a.

The bounds in Lemma 2.2 have various implications that are used in later arguments. The first concerns the derivatives of $\alpha$ and $\beta$. To state the latter, suppose that A is a given connection on E. Introduce $\nabla$ to denote the associated covariant derivative. The covariant derivative on $K^{-1}E$ that is defined by A and the canonical connection, $\gamma$, is denoted in what follows by $\nabla'$.

**Lemma 2.3:** *For each integer $q \geq 1$ and constant $c > 0$, there is a constant $\kappa \geq 1$ with the following significance: Let $\mu$ denote a co-exact 1-form and let $\varpi$ denote a harmonic 1-form. Assume that both the $C^{3+q}$ norm of $\mu$ and the $L^2$ norm of $\varpi$ are bounded by c. Fix $r \geq 1$ and suppose that $(A, \psi = (\alpha, \beta))$ is a solution to the version of (2.4) defined by r and the perturbation data $(\mathfrak{T} = *d\mu, \mathfrak{S} = 0)$ and $\varpi$. Then*
- $|\nabla^q \alpha| \leq \kappa\,(r^{q/2} + 1)$.
- $|\nabla'^q \beta| \leq \kappa\,(r^{(q-1)/2} + 1)$.

This lemma is proved in Section 6b.

There is one more apriori estimate that plays a prominent role in what follows, this an estimate for the connection A itself. To this end, introduce the functional $\mathrm{E}$ on Conn(E) that sends any given connection A to

$$\mathrm{E}(A) = i \int_M a \wedge *B_A$$

(2.6)

It is an immediate consequence of Lemma 2.2 and (2.5) that if $(A, \psi)$ is a solution to (2.5), then

$$-c \leq E \leq r \, \text{vol}(M) + c ,$$
(2.7)

where an upper bound for c depends only on the $C^3$ norm of $\mu$ and the $L^2$ norm of $\varpi$.

**Lemma 2.4**: *Fix a connection $A_E$ on E and $c > 0$. There exists a constant $\kappa \geq 1$ with the following significance: Let $\mu$ denote a co-exact 1-form and let $\varpi$ denote a harmonic 1-form. Assume that both the $C^3$ norm of $\mu$ and the $L^2$ norm of $\varpi$ bounded by c. Fix $r \geq 1$ and suppose that $(A, \psi = (\alpha, \beta))$ is a solution to the version of (2.4) defined by r and the perturbation data $(\mathfrak{T} = *d\mu, \mathfrak{S} = 0)$ and $\varpi$. Then there exists a smooth map $u: M \to S^1$ such that $\hat{a} = A - u^{-1}du - A_E$ obeys $|\hat{a}| \leq \kappa (r^{2/3} E^{1/3} + 1)$.*

This lemma is proved in Section 6c.

### c) Instantons

Assume now that the bundle $\det(\mathbb{S})$ has torsion first Chern class. This means that the respective images in $H^2(M; \mathbb{Q})$ of the first Chern class of E and the form $-\frac{1}{\pi} *\varpi_K$ are equal. In what follows here and in all other sections but Section 6, the 2-form $\varpi$ that appears in (2.4) is the 2-form $\varpi_K$. This choice for $\varpi$ is assumed implicitly where not stated explicitly in what follows.

Two other functionals are used in the search for sequences that satisfy Theorem 2.1's conditions. The first is also a functional on Conn(E), this the *Chern Simons* functional. The definition of the latter requires a preliminary choice of a fiducial connection, $A_E$, on E. In this regard, choose $A_E$ so that its curvature 2-form is $*i\varpi_K$. The choice of $A_E$ identifies Conn(E) with $C^\infty(M; iT^*M)$. Use this identification to write a given connection A as $A_E + \hat{a}$. Then

$$\mathfrak{cs}(A) = -\int_M \hat{a} \wedge d\hat{a} .$$
(2.8)

Note that $\mathfrak{cs}$ is gauge invariant in the sense that $\mathfrak{cs}(A - iu^{-1}du) = \mathfrak{cs}(A)$ when u is a smooth map from M to $S^1$.

The second of the required functionals is denoted by $\mathfrak{a}$. To define $\mathfrak{a}$, introduce the function $\mathfrak{g}: \text{Conn}(E) \times C^\infty(M; \mathbb{S}) \to \mathbb{R}$ that defines the perturbation terms $\mathfrak{T}$ and $\mathfrak{S}$ that appear in (2.4). With $\mathfrak{g}$ given, the function $\mathfrak{a}$ is defined on $\text{Conn}(E) \times C^\infty(M; \mathbb{S})$ and it sends any given pair $\mathfrak{c} = (A, \psi)$ to

$$\mathfrak{a}(\mathfrak{c}) = \tfrac{1}{2}(\mathfrak{cs}(A) + 2\mathfrak{g} - r E(A)) + r\int_M \psi^\dagger D_A \psi \;.$$

(2.9)

Note that if $\mathfrak{g}$ is independent of the section of $\mathbb{S}$, thus $\mathfrak{g} = \mathfrak{g}(A)$, then

$$\mathfrak{a} = \tfrac{1}{2}(\mathfrak{cs} + 2\mathfrak{g} - rE)$$

(2.10)

if $\mathfrak{c} = (A, \psi)$ is a solution to (2.4). In particular, (2.10) holds in the case where $\mathfrak{g}$ leads to the equations in (2.5). In any event, a pair $(A, \psi)$ is a critical point of $\mathfrak{a}$ if and only if the pair satisfies (2.4).

The 'gradient flow' lines of the functional $\mathfrak{a}$ in $\mathrm{Conn}(E) \times C^\infty(M; \mathbb{S})$ also play a role in this story. Such a gradient flow line is called an *instanton* when it has an $s \to \infty$ limit and an $s \to -\infty$ limit, and both limits are solutions to (2.4). By definition, a gradient flow line of $\mathfrak{a}$ is a smooth map, $s \to (A(s, \cdot), \psi(s, \cdot))$, from $\mathbb{R}$ into $\mathrm{Conn}(E) \times C^\infty(M; \mathbb{S})$ that obeys the equation

- $\frac{\partial}{\partial s} A = -B_A + r(\psi^\dagger \tau^k \psi - ia) + \mathfrak{T}(A, \psi) + i\varpi_K,$
- $\frac{\partial}{\partial s} \psi = -D_A \psi + \mathfrak{S}(A, \psi).$

(2.11)

The latter equations can be written as $\partial_s(A,\psi) = -\nabla \mathfrak{a}|_{(A,\psi)}$ where the gradient of $\mathfrak{a}$ is defined using the $L^2$ inner product on $C^\infty(M; iT^*M \oplus \mathbb{S})$. This is to say that $\nabla \mathfrak{a}|_\mathfrak{c}$ is the section of $iT^*M \oplus \mathbb{S}$ with the property that $\frac{d}{dt} \mathfrak{a}(\mathfrak{c} + \mathfrak{b}) = \langle \nabla \mathfrak{a}, \mathfrak{b} \rangle_{L^2}$ for all sections $\mathfrak{b}$ of $iT^*M \oplus \mathbb{S}$. Here, $\langle\,,\,\rangle_{L^2}$ denotes the $L^2$ inner product.

Of interest in what follows are the instanton solutions to (2.11). For those new to (2.11), note that if $(A, \psi)$ is a solution, then so is $(A - u^{-1}du, u\psi)$ if $u$ is any smooth map from $M$ to $S^1$. Solutions that differ in this way are deemed to be gauge equivalent.

**d) A Banach space of perturbations**

Kronheimer and Mrowka introduce the notion of a *large*, separable Banach space of tame perturbations for use in (2.4). This notion is made precise in their Definition 11.6.3. Such a Banach space is denoted here by $\mathcal{P}$. In particular, functions in $\mathcal{P}$ are, in a suitable sense, dense in the space of functions on $\mathrm{Conn}(E) \times C^\infty(M; \mathbb{S})$. What follows summarizes some other features of $\mathcal{P}$ that are used here. To set the stage for this summary, note that a smooth 1-form, $\mu$, on $M$ defines the gauge invariant function $\mathfrak{e}_\mu$: $\mathrm{Conn}(E) \to \mathbb{R}$ by the rule

$$\mathfrak{e}_\mu(A) = i\int_M \mu \wedge *B_A \;.$$

(2.12)

View $\mathfrak{e}_\mu$ as a function on $\text{Conn}(E) \times C^\infty(M; \mathbb{S})$ with no dependence on the $C^\infty(M; \mathbb{S})$ factor. To complete the stage setting, let $\Omega_0$ denote the vector space of finite linear combinations of the eigenvectors of operator $*d$ on $C^\infty(M; iT^*M)$ whose eigenvalue is non-zero. All forms in $\Omega_0$ are co-exact, and $\Omega_0$ is dense in the space of co-exact 1-forms on M. In the proposition that follows, $\|\cdot\|_2$ denotes the $L^2$ norm of the indicated section of $C^\infty(M; iT^*M \oplus \mathbb{S})$ and $\|\cdot\|_{\mathcal{P}}$ denotes the norm on the indicated element in $\mathcal{P}$.

**Proposition 2.5**: *There is a large, separable Banach space of tame perturbations with the following properties*:
- *If $\mu \in \Omega_0$, then $\mathfrak{e}_\mu \in \mathcal{P}$.*
- *Let $r \geq 0$ and let $\{\mathfrak{g}_n\}_{n=1,2,\ldots}$ be a convergent sequence from $\mathcal{P}$, and let $\mathfrak{g} \in \mathcal{P}$ denote the limit. For each $n \in \{1, 2,\ldots\}$, let $\mathfrak{c}_n$ denote a solution to the version of (2.4) that is defined by $\mathfrak{g}_n$. Then there is a subsequence of $\{\mathfrak{c}_n\}$ and a corresponding sequence of gauge transformations such that the result converges in $\text{Conn}(E) \times C^\infty(M; \mathbb{S})$ to a solution to the $(r, \mathfrak{g})$ version of (2.4).*
- *The functions in $\mathcal{P}$ are smooth with respect to the Sobolev $L^2_1$ topology on $\text{Conn}(E) \times C^\infty(M; \mathbb{S})$. In particular there exists $\kappa_* > 0$ such that if $\mathfrak{g} \in \mathcal{P}$ and $\mathfrak{c} = (A, \psi) \in \text{Conn}(E) \times C^\infty(\mathbb{S})$ and $\mathfrak{b} \in C^\infty(M; iT^*M \oplus \mathbb{S})$, then*

$$\left|\tfrac{d}{dt} \mathfrak{g}(\mathfrak{c} + t\mathfrak{b})_{t=0}\right| \leq \kappa_* \|\mathfrak{g}\|_{\mathcal{P}} (1 + \|\psi\|_2) \|\mathfrak{b}\|_2.$$

*Proof of Proposition 2.5*: Kronheimer and Mrowka describe how to construct a large, separable $\mathcal{P}$ in Section 11.5 of their book [KM], and their constructions readily accommodate the requirement in the first bullet. The second bullet follows from Propositions 10.7.2 and 11.6.4 in [KM]. The third is from Item (iv) of Theorem 11.6.1 in [KM]

The pertubations that are used in what follows are assumed to come from the Banach space $\mathcal{P}$. Let $\Omega$ denote the closure of $\Omega_0$ using the norm on $\mathcal{P}$. It is important to note that all forms in $\Omega$ are smooth. In particular, a Cauchy sequence in $\Omega$ has a convergent subsequence in $C^\infty(M; iT^*M)$. Keep in mind in what follows that a perturbation term in (2.4) defined by a $\mu \in \Omega$ version of $\mathfrak{e}_\mu$ and $\varpi = \varpi_K$ leads to the equations that are depicted in (2.5).

### 3. The Seiberg-Witten Floer homology for a contact form

The purpose of this section is to explain what is meant here by Seiberg-Witten Floer homology for a contact 1-form. It is defined here for a $\text{Spin}_\mathbb{C}$ structures that arises when the bundle det(𝕊) has torsion first Chern class. Assume now that E in (2.2) is chosen so that this is the case. The resulting homology is denoted in what follows by cSWF homology.

Most of the subsequent exposition summarizes material from Kronheimer and Mrowka's elegant exposition. This being the case, the reader is referred to [KM] for the assertions that are little more than restatements of aspects of their general theory for Seiberg-Witten Floer homology.

**a) The cycles for the cSWF homology**

The cSWF homology is defined by a differential on a $\mathbb{Z}$-graded chain complex whose cycles are formal linear combinations of equivalence classes of irreducible solutions to (2.4). In this regard, remember that term $\varpi = \varpi_K$ in all of what is said here. The equivalence relation that defines the generators of the chain complex pairs $(A, \psi)$ and $(A´, \psi´)$ when $A´ = A - u^{-1}du$ and $\psi´ = u\psi$ with u a smooth map from M to $S^1$.

The set of these equivalence classes generate the cycles for the complex if the solutions to (2.4) (and the instanton solutions to (2.11)) obey certain extra conditions. A solution that obeys these conditions is said in what follows to be *non-degenerate*. As it turns out, these conditions are present if the perturbation term that appears in these equations is chosen in a suitably generic fashion. The condition on the solutions to (2.4) are discussed momentarily

The digression that follows is needed before more can be said. To start the digression, fix $(A, \psi) \in \text{Conn}(E) \times C^\infty(M; \mathbb{S})$. Use the latter to define a certain linear operator that maps $C^\infty(M; iT^*M \oplus \mathbb{S} \oplus i\mathrm{I}_\mathbb{R})$ to $C^\infty(M; iT^*M \oplus \mathbb{S} \oplus i\mathrm{I}_\mathbb{R})$, where $\mathrm{I}_\mathbb{R} = M \times \mathbb{R}$. The operator $\mathfrak{L}$ sends a triple $(b, \eta, \phi)$ to the pair whose respective $iT^*M$, $\mathbb{S}$ and $i\mathrm{I}_\mathbb{R}$ components are:

- $*db - d\phi - 2^{-1/2} r^{1/2} (\psi^\dagger \tau \eta + \eta^\dagger \tau \psi) - \mathfrak{t}_{(A,\psi)}(b, \eta)$,
- $D_A \eta + 2^{1/2} r^{1/2}(\text{cl}(b)\psi + \phi\psi) - \mathfrak{s}_{(A,\psi)}(b, \eta)$,
- $*d*b - 2^{-1/2} r^{1/2} (\eta^\dagger \psi - \psi^\dagger \eta)$,

(3.1)

where the pair $(\mathfrak{t}_{(A,\psi)}, \mathfrak{s}_{(A,\psi)})$ denotes the operator on $C^\infty(M; iT^*M \oplus \mathbb{S})$ that sends a given section $(b, \eta)$ to $(\frac{d}{dt} \mathfrak{T}(A+tb, \psi+t\eta), \frac{d}{dt} \mathfrak{S}(A+tb, \psi+t\eta))|_{t=0}$. Denote this operator by $\mathfrak{L}$. In the case of (2.5), the terms $\mathfrak{t}$ and $\mathfrak{s}$ are absent in (3.1).

In general, the operator $\mathfrak{L}$ defines an unbounded, self adjoint operator on the $L^2$ closure of $C^\infty(M; iT^*M \oplus \mathbb{S} \oplus i\mathrm{I}_\mathbb{R})$. The domain of this operator is the $L^2_1$ closure of $C^\infty(M; iT^*M \oplus \mathbb{S} \oplus i\mathrm{I}_\mathbb{R})$. The spectrum of the self-adjoint extension of $\mathfrak{L}$ is discrete, with no accumulation points. Moreover, the spectrum is unbounded in both directions on $\mathbb{R}$. See, Chapter 12.3.2 in [KM].

A solution (A, ψ) to (2.4) with ψ not identically zero is said to be *non-degenerate* when the kernel of $\mathfrak{L}$ is trivial. In the case where ψ = 0, a solution is deemed to be non-degenerate when the kernel of $\mathfrak{L}$ is spanned by (b = 0, η = 0, φ = i). One of the purposes of introducing the perturbation data term in (2.4) is to ensure that all solutions to these equations are non-degenerate. The proof of the next lemma is in Section 7a.

**Lemma 3.1:** *Given* r ≥ 0 *there is a residual set of* μ ∈ Ω *such that all of the irreducible solutions to the corresponding version of (2.5) are non-degenerate. There is an open dense set of*
$\mathfrak{g}$ ∈ $\mathcal{P}$ *such that all solutions to corresponding version of (2.4) are non-degenerate.*

The gauge equivalence class of a non-degenerate solution is isolated. It is also the case that a non-degenerate solution persists when the equations are deformed. These notions are made precise in the next lemma. The lemma reintroduces the functional $\mathfrak{e}_\mu$ on Conn(E) that appears in (2.12).

**Lemma 3.2**: *For a given* r ≥ 0 *and* μ ∈ Ω *and* $\mathfrak{q}$ ∈ $\mathcal{P}$, *suppose that* $\mathfrak{c}$ = (A, ψ) *is a non-degenerate solution to the* r *and* $\mathfrak{g}$ = $\mathfrak{e}_\mu$ + $\mathfrak{q}$ *version of (2.4). Then the following is true:*
- *There exist* ε > 0 *such that if* (A´, ψ´) *is a solution to (2.4) that is not gauge equivalent to* (A, ψ), *then the* $L^2_1$ *norm of* (A - A´, ψ - ψ´) *is greater than* ε.
- *There is a smooth map,* $\mathfrak{c}(\cdot)$, *from the ball of radius* ε *centered at the origin in* $\mathcal{P}$ *to* Conn(E) × $C^\infty$(M; $\mathbb{S}$) *such that* $\mathfrak{c}(0)$ = $\mathfrak{c}$ *and such that* $\mathfrak{c}(\mathfrak{p})$ *solves the version of (2.4) that is defined by* r *and the perturbation defined by* $\mathfrak{g}$ = $\mathfrak{e}_\mu$ + $\mathfrak{q}$ + $\mathfrak{p}$.

*Proof of Lemma 3.2*: This is a corollary of Lemmas 12.5.2 and 12.6.1 in [KM].

Each non-degenerate, irreducible solution to (2.4) has a degree that is determined by the spectral flow for the operator $\mathfrak{L}$. For those new to the notion of spectral flow, a more detailed discussion is given in Section 5a. Suffice it to say now that the spectral flow for a continuous family s → $\mathcal{L}_s$ of self-adjoint, Fredholm operators parametrized by s ∈ [0, 1] is canonically defined if the kernels of $\mathcal{L}_0$ and $\mathcal{L}_1$ are trivial. In this case, the spectral flow is the number of eigenvalues (counting multiplicity) that cross 0 ∈ $\mathbb{R}$ from below as t increases minus the number that cross 0 from above as t increases. As is explained below, in the case where K has torsion first Chern class and E = $I_\mathbb{C}$, each non-degenerate, irreducible solution to (2.4) has a canonical degree. In other cases, the degree defined below requires some auxilliary choices. In any case, there is a canonical relative degree that can be assigned to any ordered pair of irreducible, non-degenerate solutions to (2.4).

The relative degree between an ordered pair, $(\mathfrak{c}, \mathfrak{c}')$, of non-degenerate, irreducible solutions to (2.4) is deemed to be minus the spectral flow for a suitable $s \in [0, 1]$ parametrized family of self-adjoint Fredholm operators with the $s = 0$ operator equal to $\mathfrak{c}$'s version of $\mathfrak{L}$ and with the $s = 1$ operator equal to the version of $\mathfrak{L}$ defined by $\mathfrak{c}'$. As the first Chern class of $\det(\mathbb{S})$ is torsion, this notion of degree is gauge invariant. In particular, the ordered pair $\mathfrak{c} = (A, \psi)$ and $\mathfrak{c}' = (A - u^{-1}du, u\psi)$ have the same relative degree for any given $u \in C^{\infty}(M; S^1)$. Thus, the notion of a relative degree descends to the cycles that define the cSWF homology.

The definition of the absolute degree in the case $E = I_{\mathbb{C}}$ requires the following lemma.

**Lemma 3.3**: *Suppose that $\mathfrak{T}$, $\mathfrak{S}$, and $\varpi$ in (2.4) are all zero and that $E = I_{\mathbb{C}}$. There exists $r_* \geq 0$ with the following significance: The solution $(A_I, (1_{\mathbb{C}}, 0))$ to the resulting version of (2.4) is non-degenerate for all $r > r_*$.*

This lemma follows from results about $\mathfrak{L}$ that are discussed in Section 5. Its proof is deferred to Section 5e.

Granted Lemma 3.3, the spectral flow from any large r version of the operator $\mathfrak{L}$ as defined by $(A_I, (1_{\mathbb{C}}, 0))$ with $\mathfrak{t} = \mathfrak{s} = 0$ in (3.1) endows each non-degenerate solution to (2.4) with an absolute degree. These degree assignments descend to the cycles that define the cSWF complex and give this complex its canonical $\mathbb{Z}$ grading.

When E is not trivial (but $\det(\mathbb{S})$ has torsion first Chern class), the absolute degree is defined as follows: Let $A_E$ denote the connection on E that was chosen just prior to (2.8) so as to identify Conn(E) with $C^{\infty}(M; iT^*M)$. Thus, the curvature of 2-form of $A_E$ is $*i\varpi_K$. Choose a section $\psi_E$ of $\mathbb{S}$ and some perturbation data for which the resulting $r = 1$ and $\mathfrak{t} = \mathfrak{s} = 0$ version of the operator $\mathfrak{L}$ in (3.1) has trivial kernel. It is not necessary that the pair $(A_E, \psi_E)$ obey the Seiberg-Witten equations. The spectral flow between this $(A_E, \psi_E)$ version of $\mathfrak{L}$ and that defined by any given non-degenerate solution $\mathfrak{c}$ to (2.4) for any given r and perturbation data from $\mathcal{P}$ is then well defined, and minus this number is deemed to be the degree of $\mathfrak{c}$. Because $\det(\mathbb{S})$ has torsion first Chern class, these degree assignments descend to the cSWF cycles and so give the complex its $\mathbb{Z}$-grading.

In the case when the equations are given by some $\mu \in \Omega$ version of (2.5), the version of the operator $\mathfrak{L}$ for any reducible $(A = A_E - \frac{1}{2}ira + \mu,\ 0)$ has a kernel. Even so, upper and lower bounds are available for the spectral flow to a reducible solution solution of this type. These bounds are stated in the next proposition; they play a central role in subsequent arguments.

**Proposition 3.4:** *Given* $c > 0$, *there is a constant*, $\kappa > 0$ *with the following significance. Fix* $r > 0$ *and a smooth 1-form*, $\mu$, *on* M *with* $C^3$ *norm less than* $c$. *Let* $A_E$ *denote a connection on* E *whose curvature 2-form is* $*i\varpi_K$. *There exists an open neighborhood in* $\text{Conn}(E) \times C^\infty(M; \mathbb{S})$ *of the pair* $(A_E - \frac{1}{2} i r a + \mu, 0)$ *such that the degree of any non-degenerate, irreducible solution in this neighborhood to the* r *and* $\mu$ *version of (2.5) differs from* $-\frac{1}{4\pi^2} r^2 \int_M a \wedge da$ *by at most* $\kappa(r^{31/16} + 1)$.

This proposition is a corollary of Proposition 5.5.
      As the set of solutions to (2.4) is compact modulo gauge equivalence, this last proposition and Lemma 3.2 have the following important consequence:

**Proposition 3.5**: *Given an integer* k, *there exists* $r_k \geq 0$ *with the following significance. Suppose that* $\mu \in \Omega$ *has* $C^3$-*norm less than 1 and that* $r \geq r_k$.
- *All solutions to the* r *and* $\mu$ *version of (2.5) with degree k or greater are irreducible.*

*In addition, there exists* $\varepsilon > 0$ *such that if* $\mathfrak{q} \in \mathcal{P}$ *has norm less than* $\varepsilon$, *then*
- *All solutions to the* r *and* $\mathfrak{g} = \mathfrak{e}_\mu + \mathfrak{q}$ *version of (2.4) with degree k or greater are also irreducible.*
- *There is a neighborhood of the set of reducible solutions to the* r *and* $\mathfrak{g} = \mathfrak{e}_\mu + \mathfrak{q}$ *version of (2.4) such that the spectral flow from any of degree k or greater solution of the* r *and* $\mathfrak{g} = \mathfrak{e}_\mu + \mathfrak{q}$ *version of (2.4) to any non-degenerate configuration in this neighborhood is greater than* $\frac{1}{16\pi^2} r^2 \int_M a \wedge da$.
- *If all solutions to the* r *and* $\mathfrak{g} = \mathfrak{e}_\mu + \mathfrak{q}$ *version of (2.4) with degree k or greater are non-degenerate, then there are only finitely many such solutions modulo gauge equivalence.*
- *If all solutions to the* r *and* $\mu$ *version of (2.4) with degree k or greater are non-degenerate,*
   a) *There is a 1-1 correspondence between the set of solutions to the* r *and* $\mu$ *version of (2.5) with degree* k *or greater and the set of solutions to the* r *and* $\mathfrak{g} = \mathfrak{e}_\mu + \mathfrak{q}$ *version of (2.4) with degree k or greater.*
   b) *In particular, if* $\mathfrak{c}$ *is a solution to the* r *and* $\mu$ *version of (2.5) with degree k or greater, then there exists a smooth map,* $\mathfrak{c}(\cdot)$, *from the radius* $\varepsilon$ *ball in* $\mathcal{P}$ *into* $\text{Conn}(E) \times C^\infty(M; \mathbb{S})$ *such that* $\mathfrak{c}(\mathfrak{q})$ *solves the* r *and* $\mathfrak{e}_\mu + \mathfrak{q}$ *version of (2.4) and such that* $\mathfrak{c}(0) = \mathfrak{c}$.

*Proof of Proposition 3.5*: Granted Proposition 3.3, this follows from Lemma 3.2 and Proposition 2.5.

This last proposition has the following consequence. Fix a line bundle E so that the resulting version of det($\mathbb{S}$) has torsion first Chern class. Suppose that r is large, that $\mu \in \Omega$ has $C^3$ norm less than 1, and that all solutions of degree k or greater to the r and $\mu$ version of (2.5) are non-degenerate. If $\mathfrak{q} \in \mathcal{P}$ has sufficiently small norm, then the complex for the cSWF homology in degrees greater than or equal to k as defined by the solutions to the r and $\mathfrak{g} = \mathfrak{e}_\mu + \mathfrak{q}$ version of (2.4) is finitely generated.

**b) The differential in the cSWF homology**

The differential in the cSWF homology is defined by counting the instantons with appropriate algebraic weights. To make such a count, it is necessary that all instanton solutions to (2.11) with these asymptotic constraints satisfy certain constraints. To describe the relevant constraints, suppose that $s \to \mathfrak{c}(s) = (A(s, \cdot), \psi(s, \cdot))$ is an instanton solution to (2.11) as defined by some $r \geq 0$ and a given $\mathfrak{g} \in \mathcal{P}$. Such an instanton is said to be non-degenerate when there are no not everywhere zero maps, $s \to \mathfrak{b}(s)$, from $\mathbb{R}$ to $C^\infty(M; iT^*M \oplus \mathbb{S} \oplus i\mathbb{I}_\mathbb{R})$ that obey the equations

- $-\frac{\partial}{\partial s} \mathfrak{b} + \mathfrak{L}|_{\mathfrak{c}(s)} \mathfrak{b} = 0$,
- $\lim_{|s| \to \infty} \|\mathfrak{b}(s)\|_2 \to 0$.

(3.2)

Here, the notation has $\mathfrak{L}|_{\mathfrak{c}(s)}$ denoting the version of (3.1) with $A = A(s, \cdot)$ and $\psi = \psi(s, \cdot)$; and it uses $\|\cdot\|_2$ as before to denote the $L^2$ norm on M.

Suppose that $\mathfrak{c}$ and $\mathfrak{c}'$ are non-degenerate, irreducible solutions to some $r \geq 0$ and $\mathfrak{g} \in \mathcal{P}$ version of (2.4). Let $\mathcal{M} = \mathcal{M}(\mathfrak{c}, \mathfrak{c}')$ to denote the set of instantons with $s \to -\infty$ limit equal to $\mathfrak{c}$ and $s \to +\infty$ limit equal to $u\mathfrak{c}'$ for $u \in C^\infty(M; S^1)$. As Chapter 14.4 of [KM] explains, the set $\mathcal{M}$ has the local structure of the zero locus of a smooth map from $\mathbb{R}^{n+\iota}$ to $\mathbb{R}^n$ where $\iota = \text{degree}(\mathfrak{c}) - \text{degree}(\mathfrak{c}')$. Because the equations in (2.11) are invariant with respect to the constant translations of s, the space $\mathcal{M}$ has an $\mathbb{R}$-action that has a fixed point if and only if $\mathcal{M}$ is the equivalence class of the constant map $s \to \mathfrak{c}$ with $\mathfrak{c}$ a solution to (2.4).

If all instantons between $\mathfrak{c}$ and $\mathfrak{c}'$ are non-degenerate, then $\mathcal{M}$ is a smooth manifold of dimension $\iota$ with a smooth $\mathbb{R}$ action. Furthermore $\mathcal{M}/\mathbb{R}$ is compact and a smooth manifold of dimension $\iota-1$ unless $\iota = 0$. In this case, the following are true: First, there are no instantons from $\mathfrak{c}$ to $\mathfrak{c}'$ if $\text{degree}(\mathfrak{c}') > \text{degree}(\mathfrak{c})$. Second, if the degrees are equal, then $\mathfrak{c} = \mathfrak{c}'$ and $\mathcal{M}$ consists of the equivalence class of the constant instanton, this the map $s \to \mathfrak{c}$. Finally, in the case where $\text{degree}(\mathfrak{c}) = \text{degree}(\mathfrak{c}') + 1$, each instanton with $s \to -\infty$ limit $\mathfrak{c}$ and $s \to \infty$ limit $u\mathfrak{c}'$ for $u \in C^\infty(M; S^1)$ is the translate via the $\mathbb{R}$ action of some representative of the finite set $\mathcal{M}/\mathbb{R}$. These assertions restate a part of Proposition 14.5.7 in [KM].

**Lemma 3.6:** *Given* $k \in \mathbb{Z}$, *there exists* $r_k > 0$ *with the following significance: If* $r \geq r_k$ *then there exists an open, dense set of* $\mu \in \Omega$ *with* $C^3$ *norm less than 1 for which the following is true*:
- *Each solution to the* $r$ *and* $\mu$ *version of* (2.5) *with degree greater than or equal to* $k$ *is irreducible and non-degenerate.*
- *Given* $\mu$ *for which the preceding conclusions hold, there exists* $\varepsilon > 0$ *and a dense, open subset of the radius* $\varepsilon$ *ball in* $\mathcal{P}$ *such that if* $\mathfrak{q}$ *is in this set, and if* $\mathfrak{c}$ *and* $\mathfrak{c}'$ *are solutions to the* $r$ *and* $\mathfrak{g} = \mathfrak{e}_\mu + \mathfrak{q}$ *version of* (2.4) *with the degrees* $\mathfrak{c}$ *and* $\mathfrak{c}'$ *greater than or equal to* $k$ *and* $\iota(\mathfrak{c}, \mathfrak{c}') \leq 1$, *then each instanton in* $\mathcal{M}(\mathfrak{c}, \mathfrak{c}')$ *is non-degenerate.*

*Proof of Lemma 3.6*: The first bullet follows from Lemmas 3.1, 3.2 and Proposition 3.5. The second with the 'open dense' replaced by 'residual' restates part of Theorem 15.1.1 in [KM]. The fact that the set in question is open for the case $\iota(\mathfrak{c}, \mathfrak{c}') \leq 1$ follows from Theorem 16.1.3 in [KM] given that there are only a finite many solutions to the $r$ and $\mu$ version of (2.4).

A pair $(\mu, \mathfrak{q}) \in \Omega \times \mathcal{P}$ as described by Lemma 3.6 for a given integer $k$ and $r > r_k$ is called *(r, k)-admissable*.

Assume now that $(\mu, \mathfrak{q})$ is $(r, k)$ admissable for given integer $k$ and $r \geq r_k$. Let $\mathfrak{c}$ and $\mathfrak{c}'$ be solutions to (2.5) with degree($\mathfrak{c}$) $> k$ and degree($\mathfrak{c}'$) = degree($\mathfrak{c}$) - 1. In this case, each point in $\mathcal{M}(\mathfrak{c}, \mathfrak{c}')/\mathbb{R}$ has a well defined associated sign, either +1 or -1 (see, Chapter 22.1 in [KM]). Let $\sigma(\mathfrak{c}, \mathfrak{c}')$ denote the sum of these plus and minus ones with the understanding that $\sigma = 0$ when $\mathcal{M}(\mathfrak{c}, \mathfrak{c}') = \emptyset$.

Assuming that $(\mu, \mathfrak{q})$ is $(r, k)$ admissable, what follows describes the differential on the cSWF complex in degrees greater than $k$ when the generators are the gauge equivalence classes of solutions to the $r$ and $\mathfrak{g} = \mathfrak{e}_\mu + \mathfrak{q}$ version of (2.4). The differential sends a given generator, $\mathfrak{c}$, to

$$\delta \mathfrak{c} = \sum_{\mathfrak{c}'} \sigma(\mathfrak{c}, \mathfrak{c}') \, \mathfrak{c}'$$

(3.3)

where the sum is over all irreducible cycles with degree one less than that of $\mathfrak{c}$. Thus, the differential decreases degree by 1.

**Proposition 3.7:** *Given* $k \in \mathbb{Z}$ *there exists* $r_k \geq 0$ *with the following significance: Fix* $r > r_k$ *and an* (r, k)-*admissable pair* $(\mu, \mathfrak{q})$ *to define the generators and differential on the cSWF complex in degrees greater than* $k$. *Then* $\delta^2 = 0$ *on all cycles of degree greater* *than* k. *Moreover, given two* (r, k) *admissible pairs, there exists an isomorphism between the corresponding homology groups in degrees greater than* k. *In addition, the homology so defined in degrees greater than k for different values of* $r > r_k$ *are isomorphic.*

***Proof of Proposition 3.7***: The result that $\delta^2 = 0$ is Proposition 22.1.4 in [KM] given Proposition 3.5. The invariance of the homology with respect to a change in the perturbation term is Corollary 23.1.6 in [KM] granted Proposition 3.5.

Some explicit isomorphisms between the cSWF homology groups for different values of r and q are described in the next three subsections.
 The next proposition is central to all that follows.

**Proposition 3.8**: *Given* $k \in \mathbb{Z}$, *there exists* $k´ < k$ *such that the cSWF homology in degree* $k´$ *is non-trivial.*

***Proof of Proposition 3.8***: Kronheimer and Mrowka introduce in Chapter 3 of their book three $\mathbb{Z}$-graded Seiberg-Witten Floer homology groups which are denoted here by $\overline{H}_*$, $\widehat{H}_*$, and $\widecheck{H}_*$. Chapter 35.1 in [KM] says quite a bit about the groups $\overline{H}_*$ $\widehat{H}_*$ and $\widecheck{H}_*$. In particular, Corollary 35.1.4 from [KM] finds that $\widecheck{H}_*$ is non-zero in an infinite set of degrees, a set that is bounded from below, but unbounded from above. Meanwhile $\widehat{H}_*$ is non-zero in an infinite set of degrees, a set that is bounded from above and unbounded from below. Since there are only a finite set of irreducible solutions, the reducible solutions to (2.4) supply all but a finite set of classes to both $\widehat{H}_*$ and $\widecheck{H}_*$. (As is indicated by Theorem 35.1.1 in [KM], the group $\overline{H}_*$ is mostly determined by the reducible solutions.)
 These Seiberg-Witten Floer homology groups are defined using both the irreducible solutions to (2.4) and the reducible solutions, with the corresponding instantons solutions to (2.11). In this regard, the complex is defined using the gauge equivalence classes of such solutions given that r and the perturbation $\mathfrak{g}$ are such that all solutions to (2.4) are non-degenerate, and all instantons with $\iota \leq 1$ are non-degenerate. Note that there is residual set of such $\mathfrak{g}$ for any given choice of r. In particular, for fixed k and $r > r_k$, there exist pairs $(\mu, q)$ that can be used to compute both the cSWF homology in degrees less than k and also the three Seiberg-Floer homology groups. Moreover, the generators of the cSWF complex in degrees less than k are also generators in the larger Seiberg-Witten Floer complexes for $\widehat{H}_*$ and $\widecheck{H}_*$, and the differential for the cSWF homology gives a part of the respective differentials for these other two homology theories. What is missing from the cSWF differential are the instantons that limit as either $s \to \infty$ or $s \to -\infty$ to a reducible solution to (2.4).

With the preceding understood, note that if $r > r_k$ and the perturbation $\mathfrak{g} = \mathfrak{e}_\mu + \mathfrak{q}$ has very small norm, then it follows from Proposition 3.4 that the reducible solutions to (2.4) have very negative degree. This means that the contribution from the reducibles to $\widehat{H}_*$ starts at a correspondingly negative degree. As a consequence, there must be a correspondingly large set of negative degrees where the homology in $\widehat{H}_*$ comes from the irreducible solutions to (2.4). In particular, the cSFW homology can not be zero in these degrees. As r increases, there are more and more such degrees since the set of degrees where $\widehat{H}_* \neq 0$ is unbounded from below. Proposition 3.8 follows directly from this.

### c) Identifying homology defined for distinct r values: Part I

This and the subsequent two subsections describe some isomorphisms between the different r versions of cSWF homology. To set the stage here, fix $k \in \mathbb{Z}$, suppose that $r > r_k$, and suppose that $\mu \in \Omega$ is such that all solutions to the r and $\mu$ version of (2.5) with degree k or greater are non-degenerate. As there are but a finite set of gauge equivalence classes of such solutions, this condition holds for all r´ in some neighborhood of r. Moreover, the solutions vary smoothly as r´ varies in this neighborhood. This follows as a special case of Proposition 3.5. Here is a precise statement:

**Lemma 3.9**: *Fix $k \in \mathbb{Z}$, and $r > r_k$. Suppose that $\mu \in \Omega$ has $C^3$ norm less than 1, and is such that all solutions to the r and $\mu$ version of (2.5) with degree k or greater are non-degenerate.*
- *There exists a maximal open set $(r_0, r_1)$ with $r_k \leq r_0 < r < r_1 \leq \infty$ such that all solutions to the r and $\mu$ version of (2.5) with degree k or greater are non-degenerate.*
- *For each $r´ \in (r_0, r_1)$ there is a 1-1 correspondence between the respective sets of solutions wih degree k or greater to the r´ and $\mu$ version of (2.5) and to the r and $\mu$ version.*
- *In particular, if $\mathfrak{c}$ is a solution to the r and $\mu$ version of (2.5) with degree k or greater, then there is a smooth map $\mathfrak{c}(\cdot): (r_0, r_1) \to \text{Conn}(E) \times C^\infty(M; \mathbb{S})$ with $\mathfrak{c}(r) = \mathfrak{c}$, and such that $\mathfrak{c}(r´)$ obeys the r´ and $\mu$ version of (2.5) for each $r´ \in (r_0, r_1)$.*

Let k and $\mu$ be as described in this last lemma. It may not be possible to define the cSWF homology complex in degrees greater than k using a given $r´ \in (r_1, r_0)$ and $\mu$ to define the perturbation term for use in (2.4) and (2.11). The point is that the pair $(\mu, 0)$ need not be (r´, k)-admissable. In fact, the pair $(\mu, 0)$ need not be (r´, k) admissable for any r´. However, there are $\mathfrak{q} \in \mathcal{P}$ with positive norm, but as small as desired, such that $(\mu, \mathfrak{q})$ is (r´, k) admissable. Such a pair $(\mu, \mathfrak{q})$ can be used to define the cSWF complex in degrees greater than or equal to k at r = r´.

Granted this, note that if $\mathfrak{q}$'s norm is less than the constant ε in Proposition 3.5, then $\mathfrak{q}$ is in the domain of the maps $\mathfrak{c}(\cdot)$ from Proposition 3.5, and the latter identify the generators of the cSWF complex in degrees k or greater as defined using the r´ and (µ, q) versions of (2.4) and (2.11) with the gauge equivalence classes of the degree k or greater solutions to the r´ and µ version of (2.5).

The preceding motivates the definition that follows.

**Definition 3.10**: *Suppose that k and µ are as in Lemma 3.9. Let r´ ∈ ($r_0$, $r_1$) and suppose that q ∈ $\mathcal{P}$ has norm less than the constant ε in Proposition 3.5. In addition, assume that the pair (µ, q) is (r´, k) admissable. The identification just described is used in what follows to label the generators for the cSWF homology complex in degrees k and greater as defined by the solutions to the r´ and $\mathfrak{g} = \mathfrak{e}_\mu + q$ version of (2.4) by the elements of the set of gauge equivalence classes of solutions to the r and µ version of (2.5) with degree k or greater. This labeling is deemed the <u>canonical labeling</u>.*

Where appropriate, and unless directed otherwise, the canonical labeling should be assumed in the ensuing discussion.

It is not likely that any one µ ∈ Ω will be such that for all r > $r_k$, the r and µ version of (2.5) has only non-degenerate solutions in degrees k or greater. The next proposition says, among other things, that the failure here can be assumed to occur for a discrete set of r ∈ ($r_k$, ∞).

**Proposition 3.11**: *Fix k ∈ ℤ, and there is a residual subset in Ω with $C^3$ norm less than 1 and with the following properties: Let µ denote a form from this subset. There is a locally finite set {$\rho_j$} ⊂ ($r_k$, ∞) with $\rho_1 < \rho_2 < \cdots$ such that if r > $r_k$ and r ∉ {$\rho_j$}, then*
  1) *Each solution with degree k or greater to the r and µ version of (2.5) is non-degenerate.*
  2) *Define $\mathfrak{a}$ using the r and $\mathfrak{g} = \mathfrak{e}_\mu$ version of (2.9). If $\mathfrak{c}$ and $\mathfrak{c}´$ are solutions with degree k or greater to the r and µ version of (2.5) that are not gauge equivalent, then $\mathfrak{a}(\mathfrak{c}) \neq \mathfrak{a}(\mathfrak{c}´)$.*

This proposition is proved in Section 7b.

Proposition 3.11 motivates the following terminology: Let k and µ be as in this proposition. Suppose that q ∈ $\mathcal{P}$ has small norm. Say that (µ, q) is *strongly (r, k)-admissable* for a given r ∉ {$\rho_j$} when (µ, q) is (r, k)-admissable, when q is in the ball of radius ε as described in Proposition 3.5, thus in the domain of the the various maps $\mathfrak{c}(\cdot)$ from Proposition 3.5; and when the following is true: Let $\mathfrak{c}$ and $\mathfrak{c}´$ denote solutions to the r and µ version of (2.5) with degree k or greater and such that $\mathfrak{a}(\mathfrak{c}) > \mathfrak{a}(\mathfrak{c}´)$ when $\mathfrak{a}$ is defined via (2.9) using r and $\mathfrak{g} = \mathfrak{e}_\mu$. Then $\mathfrak{a}(\mathfrak{c}) > \mathfrak{a}(\mathfrak{c}´)$ when $\mathfrak{a}$ is defined via (2.9) by $\mathfrak{g} = \mathfrak{e}_\mu + q$.

Let k and µ be as in Proposition 3.11. As r varies in some ($\rho_j$, $\rho_{j+1}$), the generators of the cSWF complex as defined for any strongly (r, k) admissable (µ, q) in degrees greater than or equal to k are labeled in an r and q-independent manner as follows: Fix a degree n ≥ k. The generators of the r and q version of the cSWF complex are labeled as $\{c_\upsilon\}_{\upsilon=1,2,\ldots}$ so that $\mathfrak{a}(c_\upsilon) > \mathfrak{a}(c_{\upsilon+1})$ where $\mathfrak{a}$ can be either the r and $\mathfrak{g} = \mathfrak{e}_\mu + q$ version of (2.9) and $c_\upsilon$ is represented by a solution to the r and $\mathfrak{g} = \mathfrak{e}_\mu + q$ version of (2.4), or, equivalently, $\mathfrak{a}$ is the r and $\mathfrak{g} = \mathfrak{e}_\mu$ version of (2.9) and $c_\upsilon$ is represented by a solution to the r and $\mu$ version of (2.5). If not stated to the contrary, this labeling of a basis for the cSWF complex in degrees k and greater is implicit in what follows. The basis labeled in this way is called the *canonical basis*.

Even with a canonical labeling of the generators, there may not be an r-independent choice for the representatives of a given homology class as r varies in a given interval $(\rho_j, \rho_{j+1})$. This is because the differential still requires the choice of an appropriate, small normed element $q \in \mathcal{P}$ for any given value of r so that $(\mu, q)$ is strongly (r, k)-admissable. The next proposition describes how the representative of a given class can change as r varies.

**Proposition 3.12**: *Let* k *and* $\mu$ *be as in Proposition 3.11. Fix* $\rho_j \in (r_k, \infty)$ *from the set described in Proposition 3.11. There exists a possibly empty, but contiguous set* $\mathbb{J}(i) \subset \mathbb{Z}$, *and a corresponding sequence* $\{t_m\}_{m \in \mathbb{J}(i)} \in (\rho_i, \rho_{i+1})$ *with the following properties:*
- *The sequence is increasing, and it has no accumulation points in the open interval.*
- *For any given* $m \in \mathbb{J}(i)$, *there exists* $q_m \in \mathcal{P}$ *of small norm such that* $(\mu, q_m)$ *is strongly* (r, k) *admissable for all* $r \in [t_m, t_{m+1}]$. *As a consequence, the differential of the cSWF complex in degrees greater than k can be assumed to be independent of r as r varies in* $[t_m, t_{m+1}]$. *This differential is denoted by* $\delta_m$.
- *Let* $m \in \mathbb{J}(i)$. *In each degree greater than equal to k, there is an upper triangular, integer valued matrix,* $\mathbb{A}$, *with* 1 *on the diagonal such that* $\delta_m = \mathbb{A}^{-1} \delta_{m-1} \mathbb{A}$. *Here, both* $\delta_m$ *and* $\delta_{m-1}$ *are written with respect to the canonical basis.*

To orient the reader who is familiar with Morse/Cerf theory in finite dimensions, the change induced by the matrix $\mathbb{A}$ is the analog of a handle slide. This proposition is proved in Section 7c.

Fix $\rho_i \in \{\rho_j\}$, and let $\{t_n\}_{n \in \mathbb{J}(i)}$ be as described by Proposition 3.12. Let $t_m \in \{t_n\}$. According to Proposition 3.12, the cSWF complex in degrees greater than k have an r-independent definition as long as r varies in $[t_m, t_{m+1})$. This definition is used, often implicitly, in what follows when reference is made to the 'cSWF complex' or to a particular 'cSWF homology class' as defined for $r \in [t_m, t_{m+1})$. The isomorphisms that are supplied by Proposition 3.12 for any given $t_m \in \{t_n\}$ are used now to extend these notions so as to be able to talk about a particular cSWF homology class for values of r in $(\rho_i, \rho_{i+1})$. The following makes this precise:

**Definition 3.13**: *Let* k *and* μ *be as in Proposition 3.12. Suppose that* $\rho_i \in \{\rho_j\}$. *A class* θ *in degree greater than k for the cSWF complex as defined for the interval* $(\rho_i, \rho_{i+1})$ *is, first of all, represented in any given* $[t_m, t_{m+1})$ *for* $m \in \mathbb{J}(i)$ *by a* $\delta_m$ *closed cycle. However, if this closed cycle for* $[t_{m-1}, t_m)$ *has the form* $\sum_v z_v c_v$, *then it is represented in* $[t_m, t_{m+1})$ *by the cycle* $\sum_{v,v'} z_v (\mathbb{A}^{-1})_{v,v'} c_{v'}$ *where* $\mathbb{A}$ *is the upper triangular matrix that is supplied for* $t_m$ *by Proposition 3.12*

The identifications given in this definition are used, sometimes implicitly, to talk about a cSWF homology class defined for the interval $(\rho_i, \rho_{i+1})$.

The next subsection provides what is necessary to describe an isomorphism between the cSWF homologies in degrees greater than k as defined on intervals $(\rho_{i-1}, \rho_i)$ and on $(\rho_i, \rho_{i+1})$. An actual isomorphism is described in Subsection 3e.

### d) Identifying homology for distinct values of r: Part II

Fix an integer k, and a form μ as described in Proposition 3.11. Let $\{\rho_j\} \subset (r_k, \infty)$ denote the set that is described in this proposition, and fix $\rho_i \in \{\rho_j\}$. The purpose of this subsection and the next is to relate the cSWF homology in degrees greater than k as r crosses $\rho_i$. Here is a preview of what is in store: The strategy is to consider a path of perturbations where the changes to the Floer differential occur at discrete times along the path and such that each change is one of a handful of standard operations. Each operation has its finite dimensional Cerf theory analogy, and the latter are as follows:

- *The disappearance or appearance of a pair of flow lines between a pair of critical points that contribute with opposing signs to the differential.*
- *A handle slide.*
- *Two critical points on the same level set.*
- *The cancellation of a single pair of critical points*

(3.4)

To set the stage, fix $r_- \in (\rho_{i-1}, \rho_i)$ and $r_+ \in (\rho_i, \rho_{i+1})$. Let $m \in \mathbb{J}(i-1)$ be such that $r_- \in [t_m, t_{m+1})$, and set $q_- = q_m$. Let $m' \in \mathbb{J}(i)$ be such that $r_+ \in [t_{m'}, t_{m'+1})$ and set $q_+ = q_{m'}$. Given $\varepsilon > 0$, and both $r_-$ and $r_+$ sufficiently close to $\rho_i$, it can be assumed that both $q_-$ and $q_+$ are in the radius $\varepsilon$ ball about 0 in $\mathcal{P}$. The next task is to choose a path $r \to q(r)$ in this ball with certain desired properties. The path is parameterized by $r \in [r_-, r_+]$, it obeys $q(r_-) = q_-$ and $q(r_+) = q_+$. If $|\rho_i - r_-|$ and $|\rho_i - r_+|$ are sufficiently small, the path can be chosen to have the five properties listed next. Section 7d describes how to find a path with these properties.

*Property 1*: Let $\mathfrak{g}(r) = \mathfrak{e}_\mu + \mathfrak{q}(r)$, and let $\mathfrak{a}_{\mathfrak{g}(r)}$ denote the action functional as depicted in (2.9) using the function $\mathfrak{g}(r)$. For any $r \in (r_-, r_+)$, the value of $\mathfrak{a}_{\mathfrak{g}(r)}$ on any solution to the r and $\mathfrak{g}(r)$ version of (2.4) is within $\varepsilon^2$ of the value of the $\mathfrak{g} = \mathfrak{e}_\mu$ version of $\mathfrak{a}$ on some solution to the $r = \rho_i$ and $\mu$ version of (2.5). Moreover, there is a finite, increasing subset $\{y_n\} \subset (r_-, r_+)$ such that all solutions to the r and $\mathfrak{g}(r)$ version of (2.4) are non-degenerate when $r \notin \{y_n\}$ and such that the values of $\mathfrak{a}_{\mathfrak{g}(r)}$ distinguish the gauge equivalence classes of solutions to the r and $\mathfrak{g}(r)$ version of (2.4).

*Property 2*: Let $I \subset [r_-, r_+] - \{y_n\}$ denote a component. There is a consecutively labeled, increasing set, $\{w_n\}_{n \in \mathbb{K}(I)}$, in the interior of I that is finite or countable, but with no accumulation points in I. For each $m \in \mathbb{K}(I)$, there exists a perturbation $\mathfrak{p}_m \in \mathcal{P}$ of very small norm such that $(\mu, \mathfrak{q}(r) + \mathfrak{p}_m)$ is (k, r)-admissable at each $r \in [w_m, w_{m+1}]$. Also, $\mathfrak{p}_m$ is such that the gauge equivalence classes of solutions to the r and $\mathfrak{g}(r, m) = \mathfrak{e}_\mu + \mathfrak{q}(r) + \mathfrak{p}_m$ version of (2.4) with degree k or greater are in 1-1 correspondence with those of the r and $\mathfrak{g}(r)$ version of (2.4) with degree k or greater for all $r \in [t_m, t_{m+1}]$. This equivalence is given by the maps in Lemma 3.2 and it is such that the value of $\mathfrak{a}_{\mathfrak{g}(r,m)}$ on a given equivalence class of r and $\mathfrak{g}(r, m)$ solutions to (2.4) is very much closer to the value of $\mathfrak{a}_{\mathfrak{g}(r)}$ on its partner equivalence class of r and $\mathfrak{g}(r)$ solutions to (2.4) then it is to the value of $\mathfrak{a}_{\mathfrak{g}(r)}$ on any other r and $\mathfrak{g}(r)$ equivalence class. In particular, the ordering of the r and $\mathfrak{g}(r)$ solutions given by the values $\mathfrak{a}_{\mathfrak{g}(r)}$ is the same as that defined by $\mathfrak{a}_{\mathfrak{g}(r,m)}$ via the equivalence.

Fix $I \subset [r_-, r_+] - \{y_n\}$ and $m \in \mathbb{K}(I)$. The cSWF homology in degrees greater than k can be defined for $r \in [w_m, w_{m+1}]$ by using the r and $\mathfrak{g}(r, m)$ versions of (2.4) and (2.11). Note in this regard that the vector space of cycles in a given degree can be identified using Property 2 with a fixed vector space, this defined by the solutions to the r and $\mathfrak{g}(r)$ version of (2.4) and the latter labeled by their ordering using $\mathfrak{a}_{\mathfrak{g}(r)}$. Here, the convention is to label the basis of cycles with the larger numbered ones having smaller values of $\mathfrak{a}_{\mathfrak{g}(r)}$. This fixed, r-independent basis is called the I-canonical basis.

*Property 3*: This next property is summarized by

**Lemma 3.14**: *Fix an interval* $I \subset [r_-, r_+] - \{y_n\}$ *and* $w_m \in \mathbb{K}(I)$. *As* r *varies in* $[w_m, w_{m+1}]$, *the differentials as written for the I-canonical basis of the cSWF complex in degrees greater than k, and as defined by the* r *and* $\mathfrak{g}(r,m)$ *version of (2.11) are independent of* r. *Moreover, there exists an upper triangular, degree preserving matrix,* $\mathbb{A} = \mathbb{A}(m)$ *with 1's on the diagonal such that the differential,* $\delta_{m-1}$ *defined on* $[w_{m-1}, w_m]$ *and the differential* $\delta_m$ *defined on* $[w_m, w_{m+1}]$ *are related, when written using the I-canonical basis, by the rule* $\delta_m = \mathbb{A}^{-1} \delta_m \mathbb{A}$.

This lemma is proved in Section 7d. The behavior that is described here corresponds to the first two items that appear in (3.4).

The next property addresses behavior of the solutions to the r and $\mathfrak{g}(r)$ version of (2.4) at any given $y \in \{y_n\}$. In what follows, a solution to the $r = y$ and $\mathfrak{g}(y)$ version of (2.4) is said to be *degenerate* when the relevant version of the operator in (3.1) has a kernel. Such a degenerate solution is said to have degree k or greater if it is the limit as r $\to$ y of a degree k or greater solution to the r and $\mathfrak{g}(r)$ version of (1.4).

What follows also uses $I_-$ to denote the component of $(r_-, r_+) - \{y_n\}$ whose closure adds y as its upper endpoint; and it uses $I_+$ to denote the the component whose closure adds y as its lower endpoint.

*Property 4*: One and only one of the following two assertions is relevant:

- *All solutions to the $r = y$ and $\mathfrak{g}(y)$ version of (2.4) with degree k or greater are non-degenerate, and there is precisely one pair of distinct, gauge equivalence classes of solutions to the $r = y$ and $\mathfrak{g}(y)$ version of (2.4) that are not distinguished by $\mathfrak{a}_{\mathfrak{g}(y)}$. In addition, there exist $y_- \in I_-$ and $y_+ \in I_+$ such that if $y_0 \in [y_-, y_+]$ and if $\mathfrak{c}$ is a solution to the $r = y_0$ and $\mathfrak{g}(y_0)$ version of (2.4), then there is a smooth map, $\mathfrak{c}(\cdot)$: $[y_-, y_+] \to$ Conn(E) $\times C^\infty(\mathbb{S})$ such that $\mathfrak{c}(y_0) = \mathfrak{c}$ and $\mathfrak{c}(r)$ solves the r and $\mathfrak{g}(r)$ version of (2.4) for each $r \in [y_-, y_+]$.*
- *The function $\mathfrak{a}_{\mathfrak{g}(y)}$ distinguishes the gauge equivalence classes of solution to the $r = y$ and $\mathfrak{g}(y)$ version of (2.4) with degree k or greater. Also, there is exactly one gauge degree k or greater gauge equivalence class of solution to the $r = y$ and $\mathfrak{g}(y)$ version of (2.4) that is degenerate. In addition,*
  1) *The operator $\mathfrak{L}$ for any solution in the one anomolous gauge equivalence class has kernel dimension 1.*
  2) *The number of gauge equivalence classes of degree k or greater solutions to the r and $\mathfrak{g}(r)$ version of (2.4) changes by two as r crosses y. In addition, the number of gauge equivalence classes of solutions to the $r = y$ and $\mathfrak{g}(y)$ version of (2.4) with degree k or greater differs by 1 from the number on either side of y.*
  3) *Let $I \in \{I_-, I_+\}$ denote the component with the greater number of equivalence classes. Then there are respective representatives, $\mathfrak{c}(r)$ and $\mathfrak{c}'(r)$, of distinct equivalence classes of solutions to the r and $\mathfrak{g}(r)$ version of (2.4) that vary smoothly with $r \in I$ and converge in Conn(E) $\times C^\infty(\mathbb{S})$ as $r \to y$ to the one anomolous $r = y$ equivalence class. Also, these classes are such that $\iota(\mathfrak{c}, \mathfrak{c}') = 1$.*
  4) *Let $\mathfrak{n}$ denote a solution to the $r = y$ and $\mathfrak{g}(y)$ version of (2.4) that is not gauge equivalent to the one anomolous gauge equivalence class.. Then there is a smooth map $\mathfrak{n}(\cdot)$: $I_- \cup \{y\} \cup I_+ \to$ Conn(E) $\times C^\infty(\mathbb{S})$ such that $\mathfrak{n}(y) = \mathfrak{n}$, and such that $\mathfrak{n}(r)$ is a solution to the r and $\mathfrak{g}(r)$ version of (2.4) for all $r \in I_- \cup \{y\} \cup I_+$.*

(3.5)

*Property 5*: What follows describes how the generators of the cSWF homology in degrees greater than k change as r crosses a given y ∈ {$y_n$}. To this end, define the respective $I_-$ and $I_+$ versions of the cSWF complex and homology in degrees greater than k using the points $y_-$ and $y_+$. This is to say that $y_-$ is in some $I_-$ version of [$w_m$, $w_{m+1}$], and use the corresponding r = $y_-$ and $\mathfrak{g}$(r, m) to define the cSWF homology in degrees greater than k using these points. Use the analogous construction for $y_+$. The story here is told in three parts.

Part 1: Assume here that the first bullet in (3.5) is relevant for y. Use the maps $\mathfrak{c}(\cdot)$ to extend the $I_+$-canonical basis as defined at $y_+$ to give a new basis for the cSWF complex at $y_-$. Let $\mathfrak{c}$ and $\mathfrak{c}'$ denote the two generators that are not distinguished by $\mathfrak{a}_{\mathfrak{g}(y)}$. If $\mathfrak{c}$ and $\mathfrak{c}'$ have different degrees, then this new basis at $y_-$ is the same as the $I_-$-canonical basis. If $\mathfrak{c}$ and $\mathfrak{c}'$ have the same degree, make the convention that $\mathfrak{c}(y_+) = \mathfrak{c}_n$ and $\mathfrak{c}'(y_+) = \mathfrak{c}_{n+1}$ where $\mathfrak{c}_n$ and $\mathfrak{c}_{n+1}$ are $I_+$-canonical basis elements at $y_+$. With respect to the $I_-$-canonical basis at $y_-$, either $\mathfrak{c}(y_-) = \mathfrak{c}_n$ and $\mathfrak{c}'(y_-) = \mathfrak{c}_{n+1}$, or else $\mathfrak{c}(y_-) = \mathfrak{c}_{n+1}$ and $\mathfrak{c}'(y_-) = \mathfrak{c}_n$. If the labelings do not change, then the respective $I_-$ and $I_+$ canonical basis for the cSWF complexes as defined at $y_-$ and $y_+$ agree. If these canonical basis agree, either for this reason, or because $\mathfrak{c}$ and $\mathfrak{c}'$ have distinct degrees, then the differential, $\delta_-$, at $y_-$ is related to the differential, $\delta_+$, defined at $y_+$ as follows: $\delta_+ = \mathbb{A}^{-1}\delta_-\mathbb{A}$, where $\mathbb{A}$ is a degree preserving, upper triangular matrix with 1's on the diagonal.

Suppose now that $\mathfrak{c}$ and $\mathfrak{c}'$ have the same degree and the labelings change as r crosses y. Let d denote the degree of $\mathfrak{c}$ and $\mathfrak{c}'$. In this case, the differentials are again related by $\delta_+ = \mathbb{A}^{-1}\delta_-\mathbb{A}$, where $\mathbb{A}$ is a degree preserving matrix of the following sort: In degrees not equal to d, the matrix $\mathbb{A}$ is upper triangular with 1's on the diagonal. In degree d,

- $\mathbb{A}_{n,n} = \mathbb{A}_{n+1,n+1} = 0$ *and* $\mathbb{A}_{n,n+1} = \mathbb{A}_{n+1,n} = 1$.
- $\mathbb{A}_{v,v} = 1$ *if* $v \neq n$ *or* n+1.
- $\mathbb{A}_{v,v'} = 0$ *if* $v > v'$ *and* $(v, v') \neq$ (n+1, n).

(3.6)

Part 2: Assume here that the second bullet in (3.5) describes the situation and that $I = I_-$. Let $\mathfrak{c}$ and $\mathfrak{c}'$ denote respective representatives of the two equivalence classes that are defined at $y_-$ and do not extend across y; and let d+1 and d denote their respective degrees. The maps that are supplied by the Item 4) of the second bullet in (3.5) are used in what follows to identify the remaining generators for the $I_-$-canonical basis at $y_-$ with the generators of the $I_+$-canonical basis at $y_+$. This identifies the full $I_-$-canonical basis at $y_-$ with the full $I_+$-canonical basis at $y_+$ in degrees different from d and d+1, and does so as the identity map. In degree d+1, the canonical basis at $y_+$ is obtained from that at $y_-$ by deleting the generator $\mathfrak{c}$; and in degree d, the change is deletion of the generator $\mathfrak{c}'$. Note that this identification preserves the ordering given by the value of $\mathfrak{a}_{\mathfrak{g}(r)}$. Let $\mathbb{V}_+$ denote the vector space of cycles as defined for $y_+$. What with the identifications just made, the vector space of cycles for $y_-$ is then $\mathbb{Z}\mathfrak{c} \oplus \mathbb{Z}\mathfrak{c}' \oplus \mathbb{V}_+$. Let $\delta_+$ denote the cSWF differential on $\mathbb{V}_+$ and let $\delta_-$ denote that on $\mathbb{Z}\mathfrak{c} \oplus \mathbb{Z}\mathfrak{c}' \oplus \mathbb{V}_+$.

**Lemma 3.15**: *There is a degree preserving homomorphism,* $\mathbb{T}: \mathbb{Z}\mathfrak{c} \oplus \mathbb{Z}\mathfrak{c}' \oplus \mathbb{V}_+ \to \mathbb{V}_+$ *with the following properties*:
- $\mathbb{T}\delta_- = \delta_+\mathbb{T}$.
- $\mathbb{T}$ *induces an isomorphism on homology*
- $\mathbb{T}$ *maps* $\mathbb{V}_+$ *to itself as an upper triangular matrix with 1's on the diagonal.*
- *The value of* $\mathfrak{a}_{\mathfrak{g}(y)}$ *on any generator that appears in* $\mathbb{T}\mathfrak{c}$ *is less than* $\mathfrak{a}_{\mathfrak{g}(y)}(\mathfrak{c})$.
- *The value of* $\mathfrak{a}_{\mathfrak{g}(y)}$ *on any generator that appears in* $\mathbb{T}\mathfrak{c}'$ *is less than* $\mathfrak{a}_{\mathfrak{g}(y)}(\mathfrak{c}')$.

This lemma is proved in Section 7d. What follows states a key implication.

**Lemma 3.16**: *Let* $\mathfrak{u} \in \mathbb{V}_+$ *denote the class such that* $\mathbb{T}\mathfrak{u} = \mathbb{T}\mathfrak{c}$ *and let* $\mathfrak{v} \in \mathbb{V}_+$ *denote the class such that* $\mathbb{T}\mathfrak{v} = \mathbb{T}\mathfrak{c}'$. *Then there exists* $A \in \{\pm 1\}$ *such that* $\delta_-(\mathfrak{c} - \mathfrak{u}) = A(\mathfrak{c}' - \mathfrak{v})$. *As a consequence, there exists* $\mathfrak{n} \in \mathbb{V}_+$ *of degree d with* $\mathfrak{a}_{\mathfrak{g}(y)}(\cdot) < \mathfrak{a}_{\mathfrak{g}(y)}(\mathfrak{c})$ *on the generators that appear in* $\mathfrak{n}$, *and such that* $\delta_-\mathfrak{c} = A\mathfrak{c}' + \mathfrak{n}$.

*Proof of Lemma 3.16*: Let $\mathfrak{v} \in \mathbb{V}_+$ be the class with degree d such that $\mathbb{T}\mathfrak{v} = \mathbb{T}\mathfrak{c}'$. Note that $\mathfrak{a}_{\mathfrak{g}(y)}(\cdot) < \mathfrak{a}_{\mathfrak{g}(y)}(\mathfrak{c}')$ on all generators that appear in $\mathfrak{v}$. The first three bullets in Lemma 3.15 imply that $\delta_-(\mathfrak{c}' - \mathfrak{v}) = 0$. Since $\mathbb{T}(\mathfrak{c}' - \mathfrak{v}) = 0$, the class $\mathfrak{c}' - \mathfrak{v}$ must be exact so as not to run afoul of the second bullet in the lemma. Thus, $\mathfrak{c}' - \mathfrak{v} = \delta_-(\mathfrak{w} + K\mathfrak{c})$ with $\mathfrak{w} \in \mathbb{V}_+$ and $K \in \mathbb{Z}$. Now let $\mathfrak{u} \in \mathbb{V}_+$ denote the class in degree d + 1 with $\mathbb{T}\mathfrak{u} = \mathbb{T}\mathfrak{c}$. Note that $\mathfrak{a}(\cdot) \leq \mathfrak{a}(\mathfrak{c})$ on all generators that appear in $\mathfrak{u}$. The first bullet of Lemma 3.15 demands that $\delta_+\mathbb{T}(\mathfrak{w} + K\mathfrak{c}) = 0$ where $\mathfrak{w} \in \mathbb{V}_+$ and $K \in \mathbb{Z}$. The second bullet of the lemma then requires that $\mathbb{T}(\mathfrak{w} + K\mathfrak{c}) = \delta_+\mathbb{T}(\mathfrak{o})$ with $\mathfrak{o}$ a degree d + 2 class. Another appeal to the first bullet finds that $\delta_-\mathfrak{o} = \mathfrak{w} + K\mathfrak{c} + A(\mathfrak{u} - \mathfrak{c})$ for some $A \in \mathbb{Z}$. Thus, $A\delta_-(\mathfrak{c} - \mathfrak{u}) = \mathfrak{c}' - \mathfrak{v}$. Since $\mathfrak{c}'$ is a generator, $|A| = 1$; and thus this last equation can be rewritten as $\delta_-\mathfrak{c} = A\mathfrak{c}' + \mathfrak{n}$ where $\mathfrak{n} = \delta_-\mathfrak{u} - A\mathfrak{v}$. The fact that $\delta_-$ decreases $\mathfrak{a}$ implies that $\mathfrak{a}_{\mathfrak{g}(y)}(\cdot) < \mathfrak{a}_{\mathfrak{g}(y)}(\mathfrak{c})$ on all generators that appear in $\mathfrak{n}$.

  Part 3: Suppose that the second bullet in (3.5) describes the situation, but assume now that $I = I_+$. Let $\mathfrak{c}$ and $\mathfrak{c}'$ denote respective representatives of the two equivalence classes that do not extend across y and let d+1 and d denote their respective degrees. Use the maps supplied by the Item 4) of the second bullet in (3.5) to identify the remaining generators for the $I_+$-canonical basis at $y_+$ with the generators of the $I_-$-canonical basis at $y_-$. As before, this identification preserves the ordering given by the value of $\mathfrak{a}_{\mathfrak{g}(y)}$. Let $\mathbb{V}_-$ denote the vector space of cycles as defined at $y_-$. With the preceding identification understood, the vector space of cycles at $y_+$ is $\mathbb{Z}\mathfrak{c} \oplus \mathbb{Z}\mathfrak{c}' \oplus \mathbb{V}_-$. Let $\delta_-$ denote the cSWF differential on $\mathbb{V}_-$ and let $\delta_+$ denote cSWF differential on $\mathbb{Z}\mathfrak{c} \oplus \mathbb{Z}\mathfrak{c}' \oplus \mathbb{V}_-$.

The lemma that follows describes what can be said in this case.

**Lemma 3.17**: *There is a degree preserving homomorphism* $\mathbb{T}: \mathbb{V}_- \to \mathbb{Z}\mathfrak{c} \oplus \mathbb{Z}\mathfrak{c}' \oplus \mathbb{V}_-$ *with the following properties:*
- $\mathbb{T}\delta_- = \delta_+\mathbb{T}$.
- $\mathbb{T}$ *induces an isomorphism on homology*
- $\mathbb{T}$ *is upper triangular with ones on the diagonal in degrees different from d+1 and d,*
- *If* $\mathfrak{u}$ *has degree d+1, then* $\mathbb{T}\mathfrak{u} = \mathbb{A}\mathfrak{u} + \kappa_\mathfrak{u}\mathfrak{c}$ *where* $\mathbb{A}: \mathbb{V}_+ \to \mathbb{V}_+$ *is an upper triangular matrix with 1's on the diagonal. Here,* $\kappa_\mathfrak{u} = 0$ *for a generator* $\mathfrak{u}$ *if* $\mathfrak{a}_{\mathfrak{g}(y)}(\mathfrak{u}) < \mathfrak{a}_{\mathfrak{g}(y)}(\mathfrak{c})$.
- *If* $\mathfrak{v}$ *has degree d, then* $\mathbb{T}\mathfrak{v} = \mathbb{A}\mathfrak{v} + \kappa_\mathfrak{v}\mathfrak{c}'$ *where* $\mathbb{A}: \mathbb{V}_+ \to \mathbb{V}_+$ *is an upper triangular matrix with 1's on the diagonal. Here,* $\kappa_\mathfrak{v} = 0$ *for a generator* $\mathfrak{v}$ *if* $\mathfrak{a}_{\mathfrak{g}(y)}(\mathfrak{v}) < \mathfrak{a}_{\mathfrak{g}(y)}(\mathfrak{c})$.

This lemma is also proved in Section 7d.

**e) Identifying homology for distinct values of r: Part III**

Fix an integer k, and a form $\mu$ as described in Proposition 3.11. Let $\{\rho_j\} \subset (r_k, \infty)$ denote the set that is described in this proposition, and fix $\rho_i \in \{\rho_j\}$. The purpose of this subsection is to complete the story started in the previous subsection by describing some explicit isomorphisms that relate the respective cSWF homology groups in degrees k or less for $(\rho_{i-1}, \rho_i)$ and $(\rho_i, \rho_{i+1})$.

For this purpose, pick $\varepsilon > 0$ but very small, and then pick $r_- \in (\rho_{i-1}, \rho_i)$ and $r_+ \in (\rho_i, \rho_{i+1})$, both very close to $\rho_i$ as described at the start of Section 3d. Fix $r \to q(r)$ for $r \in [r_-, r_+]$ to be the path in $\mathcal{P}$ as described in Section 3d. Let $\{y_n\}$ be as described in Property 1 of Section 3d. Lemma 3.14 identifies the various cSWF complexes and their differentials as r varies in any given component of $[r_-, r_+] - \{y_n\}$. Use the constructions from Property 5 and Lemmas 3.15 and 3.17 to identify the cSWF homology in consecutive intervals of $[r_-, r_+] - \{y_n\}$. Compose these homomorphism to obtain a homomorphism, $\mathbb{U}$, between the cSWF homology defined for $r_-$ with that defined for $r_+$. Note that $\mathbb{U}$ depends on the choice of the path $q(\cdot)$, and on the data from Properties 2 and 5 as well. A choice for $\mathbb{U}$ is described in the next section.

## 4. Max/min and estimates for $\mathfrak{a}$, $E$, and $\mathfrak{cs}$

This section studies the r-dependence of the values taken by the functionals $\mathfrak{a}$, $E$ and $\mathfrak{cs}$ on cycles that represent cSWF homology classes.

### a) Continuity with respect to r of the functional $\mathfrak{a}$

As is indicated by Proposition 3.12 and Lemmas 3.14, 3.15 and 3.17, a given cSWF homology class may not admit r dependent choices for its cycle representative that vary with r in a continuous fashion. Even so, it is possible to assign to such a class a continuous function that is defined up to any given large value of r, and whose value at all but a discrete set of r is the value of $\mathfrak{a}$ on some generator in a representing cycle. Here and in the remainder of this section, the functional $\mathfrak{a}$ is defined using $\mathfrak{g} = \mathfrak{e}_\mu$ in (2.9). The definition that follows describes how this is done.

**Definition 4.1:** *Fix an integer k, and a form $\mu$ as described in Proposition 3.11. Let $\{\rho_j\} \subset (r_k, \infty)$ be as described in this same proposition, and fix $\rho_i \in \{\rho_j\}$. Let $\{t_n\}_{n \in J(i)}$ be as described in Proposition 3.12. Fix $t_m \in \{t_n\}_{n \in J(i)}$ and introduce the perturbation $\mathfrak{q}_m$ from Proposition 3.12. Given $r \in [t_m, t_{m+1})$, use r and the perturbation $\mathfrak{g} = \mathfrak{e}_\mu + \mathfrak{q}_m$ to define the cSWF complex in degrees k and greater. Use the canonical labeling from Definition 3.10 to identify the generators with the solutions to the r and $\mu$ version of (2.5). Let $\theta$ denote a non-zero class with degree greater than k in the resulting cSWF homology. Suppose that $\mathfrak{n} = \sum_v z_v \mathfrak{c}_v$ is a cycle defined for the given value of r and $\mu$ that represents the class $\theta$. Define $\hat{\mathfrak{a}}(\mathfrak{n},r)$ to be the maximum value of $\mathfrak{a}$ on the set of gauge equivalence classes of solutions to the r and $\mu$ version of (2.5) that appear in the sum for $\mathfrak{n}$. Then define $\mathfrak{a}_\theta(r)$ to be the minimum over all such $\mathfrak{n}$ of the values of $\hat{\mathfrak{a}}(\mathfrak{n},r)$.*

The r-dependence of $\mathfrak{a}_\theta$ has a crucial role in this story. The next proposition addresses this issue.

**Proposition 4.2:** *Fix k and $\mu$ as in Definition 4.1. Let $\rho_i \in \{\rho_j\}_{j \geq 1}$. Use Definition 3.13 to identify the cSWF complexes as defined by Proposition 3.12's perturbations $\{\mathfrak{q}_m\}_{m \in J(i)}$ on the components of $(\rho_i, \rho_{i+1}) - \{t_n\}_{n \in J(i)}$. With this identification understood, let $\theta$ denote a cSWF homology class in degree greater than k. Then the function $r \to \mathfrak{a}_\theta(r)$ as described above for $r \in (\rho_i, \rho_{i+1}) - \{t_n\}_{n \in J(i)}$ defines a piece-wise differentiable, continuous function on $(\rho_i, \rho_{i+1})$. Moreover, for each index i, there is an isomorphism between the respective cSWF homologies in degree greater k as defined for $(\rho_{i-1}, \rho_i)$ and $(\rho_i, \rho_{i+1})$ such that with these isomorphisms understood, the following is true: Fix a cSWF homology class $\theta$ with degree greater than k. Then the function $r \to \mathfrak{a}_\theta(r)$ defines a piecewise differentiable, continuous function on $(r_k, \infty)$.*

***Proof of Proposition 4.2***: Consider first the behavior in an interval $(t_m, t_{m+1}) \subset (\rho_i, \rho_{i+1})$. The first point is that $\mathfrak{a}_\theta(\cdot)$ varies smoothly since the generators of the cSWF complex and the differential do not change with r in such an interval. The second point is that $\mathfrak{a}_\theta(\cdot)$ has a unique limit from above and also from below as r approaches any tiven $t_m$. This follows from Proposition 2.5. Consider next the behavior as r crosses a given $t_m$. Let $\mathfrak{n} = \sum_v Z_v \mathfrak{c}_v$ represent $\theta$ for $r \in (t_{m-1}, t_m)$. Let $\mathbb{A}$ denote the matrix supplied by Definition 3.13. The matrix $\mathbb{A}^{-1}$ acts so as to add to any given $\mathfrak{c}_v$ only multiples of generators on which $\mathfrak{a}$ has value less than $\mathfrak{a}(\mathfrak{c}_v)$. As a consquence, $\hat{\mathfrak{a}}(\mathbb{A}^{-1}\mathfrak{n},r) = \hat{\mathfrak{a}}(\mathfrak{n},r)$; and so $\mathfrak{a}_\theta$ is continuous at $t_m$. This proves that $\mathfrak{a}_\theta$ extends as a continuous and piecewise differentiable on the whole of any $(\rho_i, \rho_{i+1})$.

Consider now the behavior of $\mathfrak{a}_\theta$ when r crosses $\rho_i \in \{\rho_j\}$. In this regard, note first that $\mathfrak{a}_\theta(\cdot)$ has a unique limit from above and a unique limit from below as r approaches $\rho_i$. Again, this follows from Proposition 2.5. Granted that such is the case, the next task is to identify the respective cSWF homologies for r just less than $\rho_i$ and for r just greater than $\rho_i$. For this purpose, return to the milieu of Section 3d and its notation. Fix $\varepsilon > 0$ but very small, and pick $r_- \in (\rho_{i-1}, \rho_i)$ and $r_+ \in (\rho_i, \rho_{i+1})$, both very close to $\rho_i$ as described at the start of Section 3d. Fix $r \to \mathfrak{q}(r)$ for $r \in [r_-, r_+]$ to be a path in $\mathcal{P}$ as described in Section 3d. Fix a component $I \subset [r_-, r_+]-\{y_m\}$, and let $\{w_n\}_{m \in \mathbb{K}(I)} \subset I$ be as described in Property 2 of Section 3d. For $r \in [w_m, w_{m+1})$, use $\mathfrak{g}(r) = \mathfrak{e}_\mu + \mathfrak{q}(r) + \mathfrak{p}_m$ to define the cSWF complex in degrees k and greater. As indicated in Lemma 3.14, the respective cSWF complexes and differentials as defined for different values of r in any given fixed interval of $I-\{w_n\}$ do not vary with r when the canonical basis is used. Meanwhile, Parts 1, 2 and 3 from Property 5 of Section 3d with (3.6) and Lemmas 3.15 and 3.17 provide homomorphisms between the cSWF complexes that are defined in contiguous components of $[r_-, r_+]-\{y_m\}$, and that descend as isomorphisms to the respective homology groups. The composition of consecutive isomorphisms gives an isomorphism between the cSWF homology defined for any component of $[r_-, r_+]-\{y_m\}$, These isomorphisms are used implicitly in what follows.

With the preceding understood, for each component $I \in [r_-, r_+]-\{y_m\}$ and then each $r \in I-\{w_n\}$, define the function $r \to \mathfrak{a}_{\mathfrak{g}(r),\theta}$ using the prescription in Definition 4.1 with $\mathfrak{a}_{\mathfrak{g}(r)}$ replacing $\mathfrak{a}$. As is explained next, Proposition 4.2 follows from

**Lemma 4.3**: *If $\varepsilon > 0$ is small, then the function $r \to \mathfrak{a}_{\mathfrak{g}(r),\theta}$ defines a continuous, piecewise differentiable function on $[r_-, r_+]$. Moreover, its total change between $r_-$ and $r_+$ is less than $\varepsilon$.*

To see why this lemma implies the proposition, choose a decreasing sequence $\{\varepsilon_p\}_{p=1,2,...}$ with $\lim_{p \to \infty} \varepsilon_p = 0$, an increasing sequence $\{r_{p-}\} \subset (\rho_{i-1}, \rho_i)$ with limit $\rho_i$, and a decreasing sequence $\{r_{p+}\} \subset (\rho_i, \rho_{i+1})$ with limiit $\rho_i$ for use in Section 3d. For each p, let $\mathbb{U}_p$ denote the resulting isomorphism from the $r_{p-}$ version of the cSWF homology in a given degree greater than k to its analog at $r_{p+}$ as described in Section 3e. Compose the latter with the isomorphisms given by Proposition 3.12 for the intervals $[r_{1-}, r_{p-}]$ and $[r_{p+}, r_{1+}]$ to obtain an isomorphism, $\mathbb{W}_p$, from the cSWF homology as defined at $r_{1-}$ to that defined at $r_{1+}$. As

$$|\mathfrak{a}_\theta(r_{p-}) - \mathfrak{a}_{\mathbb{U}_p\theta}(r_{p+})| \leq \varepsilon_p,$$

(4.1)

the proposition's claims about continuity as r crosses $\rho_i$ follow provided that the sequence $\{\mathbb{W}_p\}$ can be suitably modified so as to have a limit as $p \to \infty$. To see that such is the case, define a partial ordering on the cSWF classes as defined for $r \in (\rho_i, \rho_{i+1})$ as follows: Say that $\theta \geq \theta'$ when the $r \to \rho_i$ limit of the function $\mathfrak{a}_\theta(\cdot)$ is no less than that of $\mathfrak{a}_{\theta'}(\cdot)$. It follows from (4.1) that if p and p´ are both sufficiently large, then $\mathbb{U}_p = \mathbb{A}\mathbb{U}_{p'}$ where $\mathbb{A}$ is an isomorphism that preserves this ordering. Note that isomorphisms with this property form a subgroup of the group of degree preserving isomorphisms of the cSWF homology. This understood, the sequence $\{\mathbb{U}_p\}$ can be modified by composing with such isomorphisms so as to be constant for p sufficiently large. After this modification, the resulting, now modified version of the sequence $\{\mathbb{W}_p\}$ converges, and it has the properties claimed by Proposition 4.2.

**Proof of Lemma 4.3**: Suppose that I is a component of $[r_-, r_+]-\{y_m\}$. As noted by Lemma 3.14, when using the canonical basis, the following is true: As r varies in any interval of $I-\{w_n\}$, neither the generators nor the differentials change for the version of the cSWF complex in degrees k or greater as defined by r and the $\mathfrak{g}(r) = \mathfrak{e}_\mu + \mathfrak{q}(r)$ versions of (2.4) and (2.11). This being the case, $\mathfrak{a}_{\mathfrak{g}(\cdot),\theta}$ defines a smooth function on any interval of $I-\{w_n\}$. Proposition 2.5 guarantees that this function has a unique limit as $r \to w \in \{w_n\}$ from above, and also a unique limit as $r \to w$ from below.

Consider now what happens as r crosses a given w from $\{w_n\}$. Lemma 3.14 guarantees that $\mathfrak{a}_{\mathfrak{g}(\cdot),\theta}$ is continuous across w because an upper triangular isomorphism doesn't change the value of this function.

Consider next the behavior as r crosses a given $y \in \{y_m\}$. Suppose that the first bullet of (3.5) is relevant at y. Since the isomorphism that is described in Part 1 of Property 5 of Section 3d acts either as an upper triangular matrix, or a matrix that switches two generators on which $\mathfrak{a}_{\mathfrak{g}(\cdot)}$ agree but is otherwise upper triangular, it follows that $\mathfrak{a}_{\mathfrak{g}(\cdot),\theta}$ is continuous across y.

Finally, consider what happens in the case when the second bullet in (3.5) is relevant for the given $y \in \{y_m\}$. Suppose first that the situation is that described in Part 2 of Section 3d. Return to the notation used in Part 2. If a class $\theta$ has degree different from either d or d+1, then it follows from the fact that Lemma 3.15's matrix $\mathbb{T}$ is upper triangular that $\mathfrak{a}_{g(\cdot),\theta}$ is continuous across y.

Suppose that θ has degree d, and is represented by for r just less than y by the cycle $\phi = B\mathfrak{c}' + \mathfrak{w}$ where $\mathfrak{w} \in \mathbb{V}_+$ and $B \in \mathbb{Z}$. If $B = 0$, then $\mathbb{T}$ acts on $\mathfrak{v}$ as an upper triangular matrix and so the limiting value $\hat{\mathfrak{a}}_{g(r)}(\cdot,r)$ on $\phi$ as $r \to y$ from below is the same as the limiting value as $r \to y$ from above of this function on $\mathbb{T}\phi$. In the case where $B \neq 0$, there are two cases to consider. Either the $r < y$ version of $\hat{\mathfrak{a}}_{g(r)}(\phi,r)$ is greater than all generators in $\mathfrak{v}$ or not. If not, then it follows from the third and fifth bullets in Lemma 3.15 that the limiting value $\hat{\mathfrak{a}}_{g(r)}(\cdot,r)$ on $\phi$ as $r \to y$ from below is the same as the limiting value as $r \to y$ from above of this function on $\mathbb{T}\phi$. If so, then it follows from Lemma 3.16 that $\phi' = B\mathfrak{v} + \mathfrak{w}$ represents θ also, and that the $r < y$ version of $\hat{\mathfrak{a}}_{g(r)}(\cdot,r)$ on $\phi'$ is less than on $\phi$. This then implies that $\mathfrak{a}_{g(\cdot),\theta}$ is continuous across y using the following two facts: First, $\phi' \in \mathbb{V}_+$ and $\mathbb{T}$ acts on $\mathbb{V}_+$ as an upper triangular matrix. Second, $\mathbb{T}\phi = \mathbb{T}\phi'$.

Finally, suppose that θ has degree d + 1 and is represented by $\phi = B\mathfrak{c} + \mathfrak{w}$ where $\mathfrak{w} \in \mathbb{V}_+$. If $B = 0$, then $\phi \in \mathbb{V}_+$ and the continuity of $\mathfrak{a}_{g(\cdot),\theta}$ across y again follows from the fact that $\mathbb{T}$ is upper triangular on $\mathbb{V}_+$. On the other hand, suppose that $B \neq 0$. It then follows from Lemma 3.16 that $\delta_- \mathfrak{w} = -B(A\mathfrak{c}' + \mathfrak{n})$, and since $\delta_-$ decreases the value of $\mathfrak{a}_{g(r)}$, it must be the case that the value of the $r < y$ version of $\hat{\mathfrak{a}}_{g(r)}(\cdot,r)$ on $\phi$ is that of $\mathfrak{a}_{g(r)}$ on a generator in $\mathfrak{w}$, and hence in $\mathbb{V}_+$. This then implies the continuity of $\mathfrak{a}_{g(\cdot),\theta}$ as r crosses y using the third and fourth bullets in Lemma 3.15.

Now consider the case where the second bullet in (3.7) is relevant for $y \in \{y_m\}$ and the situation is as described in Part 3 in Property 5 of Section 3d. As before, there is no change in $\mathfrak{a}_{g(r),\theta}$ as r crosses y if the degree of θ differs from either d or d+1. This is because $\mathbb{T}$ acts in these degrees as an upper triangular matrix. There is also no change when θ has degree d or d+1 because in these degrees, $\mathbb{T}$ acts on a given generator $\mathfrak{w} \in \mathbb{V}_+$ so as to give $\mathfrak{w} + \phi$, where $\phi$ is a sum of generators of $\mathbb{Z}\mathfrak{c} \oplus \mathbb{Z}\mathfrak{c}' \oplus \mathbb{V}_+$ on which $\mathfrak{a}_{g(r)}$ is less than $\mathfrak{a}_{g(r)}(\mathfrak{w})$.

The fact that the total change in $\mathfrak{a}_{g(\cdot),\theta}$ is bounded by ε follows from the preceding given Property 1 from Section 3d.

The isomorphisms provided by Proposition 4.2 are used now implicitly in the discussion that follows when reference is made to a particular cSWF homology class. This understood, let θ denote a cSWF homology class. The function $\mathfrak{a}_\theta$ is used next to define three more functions.

**Definition 4.4**: *Fix an integer k, and let* μ *be as described by Proposition 3.11. Let* $\{\rho_j\}$ $\subset (r_k, \infty)$ *be as described in this same proposition. Let* θ *denote a non-zero cSWF* *homology class of degree greater than* k. *Fix* $\rho_i \in \{\rho_j\}$ *and* $r \in (\rho_i, \rho_{i+1}) - \{t_n\}_{n \in J(i)}$. *Suppose that* $\mathfrak{n} = \sum_v z_v \mathfrak{c}_v$ *is a cycle defined for the given value of* r *and* $\mu$ *that represents the class* $\theta$ *and is such that* $\hat{\mathfrak{a}}(\mathfrak{n},r) = \mathfrak{a}_\theta(r)$. *Let* $\hat{E}(r, \mathfrak{n})$ *denote the infimum of the values of* E *on the configurations* $\mathfrak{c} \in \{\mathfrak{c}_v\}$ *that appear in the sum for* $\mathfrak{n}$ *and have* $\mathfrak{a}(\mathfrak{c}) = \mathfrak{a}_\theta(r)$. *Then, define*

- $\hat{E}(r)$ *to be the infimum of the set* $\{\hat{E}(\mathfrak{n}, r)\}$ *over all such* $\mathfrak{n}$,
- $\mathfrak{v}(r) = 2\mathfrak{a}_\theta(r) + r\hat{E}(r)$,
- $\mathfrak{f}(r) = -2r^{-1}\mathfrak{a}_\theta(r) = \hat{E}(r) - r^{-1}\mathfrak{v}(r)$.

Note that $\mathfrak{v}(r)$ is the value of $\mathfrak{cs}(\cdot) + 2\mathfrak{e}_\mu(\cdot)$ on some degree k generator $\mathfrak{c}$ with $\mathfrak{a}(\mathfrak{c}) = \mathfrak{a}_\theta(r)$ and $E(\mathfrak{c}) = \hat{E}(r)$. Here, $\mathfrak{e}_\mu$ is defined as in (2.12). Although $\mathfrak{f}$ is a continuous function of r, neither $\mathfrak{v}$ nor $\hat{E}$ need be continuous.

The next definition is motivated directly by the appearance of the energy function E in the statement of Theorem 2.1

**Definition 4.5:** *Fix an integer* k *and let* $\mu$ *be as described by Proposition 3.11. Let* $\theta$ *denote a non-trivial, cSWF homology class in degree greater than* k. *The class* $\theta$ *is said to be a <u>divergence class</u> when the conditions that are stated next hold. Given* $E > 0$, *there exists* $\rho_E \geq 0$ *with the following significance*: *Use the form* $\mu$ *to define the function* $\hat{E}$. *Then* $\hat{E}(r) > E$ *when* $r > \rho_E$.

What follows is the key observation about divergence classes.

**Proposition 4.6**: *Fix an integer* k *and let* $\mu$ *be as described by Proposition 3.11. Suppose that* $\theta$ *is a divergence class of degree greater than* k *in the cSWF homology. The class* $\theta$ *determines positive constants* c *and* $r_\theta$ *with the following significance*: *Fix* $r´ > r_\theta$ *and there exists* $r > r´$ *such that* $\mathfrak{v}(r) > \frac{1}{10} r \hat{E}(r)$ *and* $\hat{E}(r) > cr$.

The next subsection contains the proof of Proposition 4.6. Here is a key corollary:

**Corollary 4.7**: *Fix an integer* k *and let* $\mu$ *be as described by Proposition 3.11. Suppose that* $\theta$ *is a divergence class of degree greater than* k *in the cSWF homology. The class* $\theta$ *determines a constant,* $c > 0$, *with the following significance*: *Fix* $r´ > r_\theta$ *and there exists* $r > r´$ *and a solution,* $(A, \psi)$, *to the version of (2.5) determined by* r *and* $\mu$ *that has the same degree as* $\theta$ *and is such that* $\mathfrak{cs}(A) > \frac{1}{16} r E(A)$ *and* $E(A) > cr$.

*Proof of Corollary 4.7*: Let r be as in Proposition 4.5. By definition, there exists a solution $\mathfrak{c} = (A, \psi)$ to (2.5) with $\theta$'s degree, and with $\mathfrak{a}(\mathfrak{c}) = \mathfrak{a}_\theta(r)$ and $E(A) = \hat{E}(r)$. Thus, $\mathfrak{v}(r) = \mathfrak{cs}(A) + 2\mathfrak{e}_\mu(A)$. Meanwhile, $|\mathfrak{e}_\mu(A)| \leq \kappa\, E(A)$ where $\kappa$ is independent of $r$ and $\mu$. This follows from Lemma 2.2, for this lemma implies that

$$\int_M |B_A| \leq \kappa\,(E(A) + 1)$$

(4.2)

where $\kappa$ is independent of $r$ and $\mu$. As a consequence, $\mathfrak{cs}(A) \geq \frac{1}{10}(r - \kappa)E(A) - \kappa$. Because $\theta$ is a divergence class, this is larger than $\frac{1}{16} r E(A)$ if $r_\theta$ is not too small.

**b) The proof of Proposition 4.6**

The proof starts with a digression to derive some preliminary facts. With k fixed and $\mu$ as in Proposition 3.11, let $\{\rho_j\} \subset (r_k, \infty)$ denote the set that is specified by Proposition 3.11. Fix j and let $\{t_m\}_{m \in \mathbb{J}} \subset (\rho_j, \rho_{j+1})$ be as described in Proposition 3.12. Let I denote a given component of the complement in $(\rho_j, \rho_{j+1})$ of $\{t_m\}_{m \in \mathbb{J}}$. Let $r \to \mathfrak{c}(r)$ denote a path of solutions to the r and $\mu$ version of (2.5) defined for $r \in I$ as described in Lemma 3.9. It then follows that

$$\tfrac{d}{dr} E = r^{-1} \tfrac{d}{dr}(\mathfrak{cs} + 2\mathfrak{e}_\mu),$$

(4.3)

To see why (4.3) holds, view the tangent vector to this path at a given value of r as a section of $iT^*M \oplus \mathbb{S}$. Write this vector as $(b, \zeta)$. Then

$$\tfrac{d}{dr}\mathfrak{cs} = -2\int_M b \wedge *B_A,$$

(4.4)

and this is equal to

$$-2r\int_M \psi^\dagger cl(b)\psi + 2ri\int_M b \wedge *a - 2\int_M b \wedge d\mu$$

(4.5)

Note next that the left most integral in (4.5) vanishes since $\psi$ is $L^2$-orthogonal to the image of $D_A$ while $cl(b)\psi = -D_A\zeta$. Meanwhile, the middle term on the right hand side of (4.4) is $r\tfrac{d}{dr} E$ as can be seen with the help of an integration by parts and an application of the identity $da = 2*a$. Finally, the right most term in $-2\tfrac{d}{dr}\mathfrak{e}_\mu$.

Equations (4.3) and (2.10) imply that

$$\tfrac{d}{dr}\mathfrak{a} = -\tfrac{1}{2} E.$$

(4.6)

at all points in I. What follows is a consequence of (4.6).

**Lemma 4.8**: *Fix an integer k and let μ be as described by Proposition 3.11. Suppose that θ is a divergence class of degree greater than k in the cSWF homology. There exists $r_θ > r_k$ with following significance: The corresponding function $\mathfrak{a}_θ$ is less than -r and monotonically decreasing where $r > r_θ$.*

**Proof of Lemma 4.8**: Fix $E \geq 1$ and let $\rho_E$ be as described in Definition 4.5. It follows from (2.7) and Lemmas 2.2 and 2.4 that there exists a constant $c_E$ that is independent of μ and is such that $\mathfrak{a}(\mathfrak{c}) < c_E$ for all solutions $\mathfrak{c}$ to the $r = \rho_E$ and μ version of (2.5). Because $\mathfrak{a}_θ$ is continuous, it follows from this last fact and (4.6) that

$$\mathfrak{a}_θ(r) \leq - \int_{\rho_E}^{r} s\hat{E}(s)ds + c_E \leq -rE + (c_E + \rho_E E).$$

(4.7)

This last equation proves the first assertion. The second follows directly from (4.6) given that $\mathfrak{a}_θ$ is piecewise differentiable and $\hat{E}(r)$ is positive for $r \geq \rho_E$.

Lemma 4.8 implies that the function $\mathfrak{f}$ is positive and $r\mathfrak{f}(r)$ is increasing for $r \geq r_θ$. Moreover, it follows from (4.3) and (2.10) that if $[x´, x] \subset [r_k+1, \infty)$, then

$$\mathfrak{f}(x) = \int_{x´}^{x} s^{-2}\mathfrak{v}(s)ds + \mathfrak{f}(x´)$$

(4.8)

This equation ends the preliminary digression.

To get on with the proof of the proposition, note first that if $\mathfrak{v}(r) \geq \frac{1}{6} r\hat{E}(r)$, then $\hat{E}(r) \geq cr$ where c is a constant that is independent of r and μ if $r \geq 1$. To see why suppose that $(A, \psi)$ is a solution to (2.5) with $r \geq 1$ and $E(A) \geq 1$. It then follows from Lemma 2.4 and (4.2) that

$$\mathfrak{cs}(A) \leq \kappa\, r^{2/3} E(A)^{4/3}$$

(4.9)

where κ is a constant that is independent of r, μ, and the pair $(A, \psi)$. Now, suppose that $\mathfrak{cs}(A) + 2\mathfrak{e}_\mu(A) > \frac{1}{6} rE(A)$. As noted in the proof of Corollary 4.7, the value of κ can be taken such that $|\mathfrak{e}_\mu(A)| \leq \kappa\, E(A)$. As a consequence of these last two inequalities, (4.9) implies that

$$\tfrac{1}{6} r\, E(A) \leq \kappa´\, r^{2/3} E(A)^{4/3}\,,$$

(4.10)

where κ´ is also independent of r, μ and $(A, \psi)$. This last inequality can hold only when $E(A) \geq (6\kappa´)^{-3} r$.

Granted the point made in the previous paragraph,, suppose that the assertion of the proposition is false. If this is the case, then there exists $r'$ such that $\mathfrak{v}(r) \leq \varepsilon\, r\hat{E}(r)$ for all $r \in [r', \infty)$ with $\varepsilon \leq \frac{1}{10}$. In this case, the integrand for the integral in (4.8) is no larger than $2\varepsilon\, s^{-2}(s\mathfrak{f}(s))$. Indeed, this inequality holds where $\mathfrak{v}(s) < 0$ since $\mathfrak{f}(s)$ is positive. Meanwhile, where $\mathfrak{v}(s)$ is non-negative, then $\mathfrak{v} \leq \varepsilon s \hat{E} \leq 2\varepsilon(s\hat{E} - \mathfrak{v})$ which is just $2\varepsilon s\mathfrak{f}(s)$. Now, let $x_m = 2^m r'$. Since $s\mathfrak{f}(s)$ is an increasing function of $s$, (4.8) implies that

$$\mathfrak{f}(x_{m+1}) \leq 2\varepsilon\, \mathfrak{f}(x_{m+1})\, x_{m+1} \int_{x_m}^{x_{m+1}} s^{-2}\, ds + \mathfrak{f}(x_n)\, .$$

(4.11)

This then implies that $\mathfrak{f}(x_{m+1}) \leq (1-2\varepsilon)^{-1}\mathfrak{f}(x_m)$. Iterating this finds

$$\mathfrak{f}(x_m) \leq (1-2\varepsilon)^{-m}\, \mathfrak{f}(r')\, .$$

(4.12)

Save this for the moment.

Note next that (4.2) and (4.9) imply that there exists $\kappa \geq 1$ that is independent of $r$ and $\mu$ such that

$$|\mathfrak{v}(r)| \leq \kappa\, r^{2/3}\hat{E}(r)^{4/3} \leq 16\kappa\, r^{2/3}\mathfrak{f}(r)^{4/3}\, .$$

(4.13)

This with (4.8) implies that

$$\mathfrak{f}(x_{m+1}) \leq 16\kappa \int_{x_m}^{x_{m+1}} s^{-2}(s^{2/3}\mathfrak{f}(s)^{4/3})ds + \mathfrak{f}(x_m)\, .$$

(4.14)

Using the fact that $s\mathfrak{f}(s)$ is an increasing function of $s$, this last inequality leads to

$$\mathfrak{f}(x_{m+1}) \leq \kappa'\, (\mathfrak{f}(x_{m+1})/x_{m+1})^{1/3}\, \mathfrak{f}(x_{m+1}) + \mathfrak{f}(x_m)\, ,$$

(4.15)

where $\kappa'$ is independent of $r$ and $\mu$. What with (4.12), the preceeding requires that

$$\mathfrak{f}(x_{m+1}) \leq \kappa'\, (\mathfrak{f}(r')/r')^{1/3}\, (2(1-2\varepsilon))^{-m/3}\, \mathfrak{f}(x_{m+1}) + \mathfrak{f}(x_m)\, .$$

(4.16)

To procede from here, note that $\mu \leq \frac{1}{10}$ and so $2(1-2\varepsilon) \geq \frac{6}{5}$. Thus, (4.16) can be written as

$$\mathfrak{f}(x_{m+1}) \leq \mathfrak{z}\lambda^m\, \mathfrak{f}(x_{m+1}) + \mathfrak{f}(x_m)$$

(4.17)

where $\lambda = (\frac{5}{6})^{1/3} < 1$ and where $\mathfrak{z} = \kappa'\, (\mathfrak{f}(r')/r')^{1/3}$.

To finish the argument, note that there exists m(r´) such that $_3\lambda^m < 1$ when $m \geq$ m(r´). As a consequence, (4.17) finds that

$$\mathfrak{f}(x_{m+1}) = \prod_{m(r´) \leq j \leq m} (1 - {}_3\lambda^j)^{-1} \mathfrak{f}(x_{m(r´)}) .$$
(4.18)

This then implies that there exists c(r´) such that

$$\mathfrak{f}(x_m) \leq c(r´) \quad \text{for all} \quad m.$$
(4.19)

However, there can be no such uniform bound if $\theta$ is a divergence class and if $\mathfrak{v}(r) \leq \frac{1}{10} r \hat{E}(r)$ for all $r > r´$. Indeed, the latter condition finds $\mathfrak{f}(x_m) = \hat{E}(x_m) - x_m^{-1}\mathfrak{v}(x_m) \geq \frac{1}{2} \hat{E}(x_m)$; and the divergence class condition requires that $\hat{E}(x_m)$ be very large when $x_m$ is very large.

## 5. Spectral flow estimates

What follows is the principle result of this section:

**Proposition 5.1**: *Given $c > 0$ there exists a constant $\kappa$ with the following significance: Suppose that $\mu$ has $C^3$ norm less than c. Suppose that $r \geq 0$ and that $\mathfrak{c} = (A, \psi)$ is a non-degenerate solution to the $r$ and $\mu$ version of (2.5). Then the degree of $\mathfrak{c}$ differs by less than $\kappa r^{31/16}$ from $-\frac{1}{4\pi^2} \mathfrak{cs}(A)$.*

*Proof of Proposition 5.1*: As explained in Section 3, the degree of this solution is defined in terms of the spectral flow between two versions of (3.1)'s operator $\mathcal{L}$. The first version is written with the given value of r, the given form $\mu$, and the given solution $(A, \psi)$; and the second is written using some other value of r, some other one form, $\mu´$, and some fiducial pair in $\text{Conn}(E) \times C^\infty(\mathbb{S})$. As the degree is defined by the spectral flow, the proposition is proved by giving an estimate for the spectral flow. This understood, Proposition 5.1 follows directly from Lemmas 5.3 and Proposition 5.5. The proofs of the latter occupy most of the remaining subsections of Section 5.

**a) The definition of spectral flow**

To define what is meant in this article by spectral flow, suppose that $\mathbb{H}$ is a separable Hilbert space, that $\mathcal{L}$ is an unbounded, self-adjoint operator on $\mathbb{H}$ such that the operator $\mathcal{L}^2 + 1$ has compact inverse. Let $s \to q_s$ denote a real analytic map from [0, 1] into the space of self-adjoint, bounded operators on $\mathbb{H}$. Let $\mathcal{L}_s = \mathcal{L} + q_s$. Then each $\mathcal{L}_s$ is self-adjoint. In addition, each $\mathcal{L}_s$ has purely discrete spectrum, all eigenvalues are real, each has finite multiplicity, and there are no accumulation points in $\mathbb{R}$.

The spectral flow for the family $\{\mathcal{L}_s\}_{s \in [0,1]}$ is defined with the help of a certain stratified, real-analytic set in $\mathbb{R} \times [0, 1]$. This set is denoted by $\mathcal{E}$, and its stratification is depicted by

$$\mathcal{E} = \mathcal{E}_1 \supset \mathcal{E}_2 \supset \cdots .$$

(5.1)

where $\mathcal{E}_k$ consists of the set of pairs $(\lambda, s)$ such that $\lambda$ is an eigenvalue of $\mathcal{L}_s$ with multiplicity k or greater. Each $\mathcal{E}_k$ is a closed set.

It is a now standard result, see for example Chapter 7 of [Ka], that $\mathcal{E}_{k*} = \mathcal{E}_k - \mathcal{E}_{k+1}$ is an open, real analytic submanifold in $\mathbb{R} \times [0, 1]$. The collection $\{\mathcal{E}_{k*}\}$ are called the *smooth strata* of $\mathcal{E}$. The aforementioned results from [Ka] imply that when the 1-dimensional smooth strata are oriented by the pull-back from $\mathbb{R} \times [0, 1]$ of the 1-form ds, then the zero dimensional strata can be consistently oriented so that the formal, weighted sum $\mathcal{E}_* = \mathcal{E}_{1*} + 2\mathcal{E}_{2*} + \cdots$ defines a locally closed cycle in $\mathbb{R} \times [0, 1]$. This means the following: Let f denote a smooth function on $\mathbb{R} \times (0, 1)$ with compact support. Then

$$\sum_{k=1,2,\ldots} k \int_{\mathcal{E}_{k*}} df = 0$$

(5.2)

Sard's theorem finds a dense, open set $\mathbb{U} \subset \mathbb{R}$ with the property that the respective maps from a point, $*$, to $\mathbb{R} \times [0, 1]$ that send $*$ to $(\lambda, 0)$ and to $(\lambda, 1)$ are transversal to the smooth strata of $\mathcal{E}$ for all $\lambda \in \mathbb{U}$. In this language, the spectral flow for the family $\{\mathcal{L}_s\}_{s \in [0,1]}$ is defined as follows: Fix $\lambda_0 \in \mathbb{U}$ with $\lambda_0 > 0$. By Sard's theorem, there exist smooth, oriented paths $\sigma \subset \mathbb{R} \times [0, 1]$ that start at $(\lambda_0, 0)$, end at $(\lambda_0, 1)$, and are transversal to the smooth strata of $\mathcal{E}$. Such a path has a well defined intersection number with $\mathcal{E}$, this being

$$f_{\lambda_0} = \sum_{k=1,2,\ldots} \sum_{p \in \sigma \cap \mathcal{E}_{k*}} (-1)^{o(p)} k ,$$

(5.3)

where $o(p) \in \{0, 1\}$. In the case where $\sigma$ is the graph of a smooth function from $[0, 1]$ to $\mathbb{R}$, the sign $(-1)^{o(p)}$ is obtained as follows: The pull-back to a smooth, 1-dimensional stratum of $\mathcal{E}$ of the 1-form $d\lambda$ from $\mathbb{R} \times [0, 1]$ at a point $(\lambda, s)$ can be written as $\lambda' ds$ with

$$\lambda' = \langle \varsigma, (\tfrac{d}{ds} \mathfrak{q}_s) \varsigma \rangle_{\mathbb{H}} .$$

(5.4)

Here, the notation uses $\varsigma$ to denote a unit length eigenvalue of $\mathcal{L}_s$ whose eigenvalue is $\lambda$, and $\langle\,,\,\rangle_{\mathbb{H}}$ denotes the inner product on $\mathbb{H}$. The sign of $\lambda'$ at an intersection point with the image of a graph is the factor $(-1)^{o(\cdot)}$ that appears in (5.3).

The fact that $q_s$ varies with s in a real analytic fashion implies that $f_{\lambda_0}$ is independent of $\lambda_0$ if $\lambda_0$ is sufficiently close to 0. This is so when 0 is an eigenvalue of one or both of $\mathcal{L}_0$ and $\mathcal{L}_1$. This noted, the spectral flow for the family is defined to be

$$f = \lim_{\lambda_0 \to 0} f_{\lambda_0}.$$

(5.5)

This definition readily generalizes to the case where the family of operators has a continuous and piece-wise real analytic parametrization. This is to say that the parametrization has the form $s \to \mathcal{L} + q_s$ for values of s in a finite union of closed intervals, $[0, s_1] \cup [s_1, s_2] \cup \cdots \cup [s_{N-1}, s_N]$ where $s \to q_s$ is real analytic and of the form described above on each of these closed intervals.

**b) Estimating spectral flow**

What follows describes the strategy from [T4] that is employed here to estimate the spectral flow for a family $\{\mathcal{L}_s = \mathcal{L} + q_s\}_{s \in [0,1]}$. Take x either $\infty$ or positive and finite, and fix an orientation preserving diffeomorphism $\Phi: \mathbb{R} \to (-x, x)$ that sends 0 to 0. In the first application that follows, $\Phi$ is the identity map from $\mathbb{R}$ to $\mathbb{R}$. The second application takes $x < \infty$ and thus $\Phi$ more interesting. In any event, fix now $T \in (0, x)$ and let S denote the circle that is obtained from the interval $[-T, T]$ by identifying the endpoints. This circle has a fiducial point, $T_*$, that given by $\{\pm T\}$, and an orientation given by the orientation of $(-T, T)$.

Now let $\sigma = \Phi^{-1}(T)$. For each $s \in [0, 1]$, let $\mathfrak{n}_s$ denote the maximal number of linearly independent eigenvectors of $\mathcal{L}_s$ whose eigenvalue lies in $[-\sigma, \sigma]$. Use $\mathfrak{n}$ in what follows to denote the maximum from the set $\{\mathfrak{n}_s\}_{s \in [0,1]}$. An estimate for the spectral flow for the family $\{\mathcal{L}_s\}_{s \in [0,1]}$ is obtained by considering the trajectories of $\mathfrak{n}$ particles on S whose paths vary continuously and piecewise differentiably as functions of $s \in [0, 1]$.

To elaborate, introduce $\mathcal{E}^\sigma$ to denote the set $\{(\lambda, s) \in \mathcal{E}: |\lambda| < \sigma\}$, and for each k, use $\mathcal{E}^\sigma_{k*}$ to denote $\mathcal{E}_{k*} \cap \mathcal{E}^\sigma$. Each point $(\lambda, s) \in \mathcal{E}^\sigma_{k*}$ corresponds to k particles on S all at the point $\Phi(\lambda)$. If $\mathcal{E}^\sigma_{k*}$ is 1-dimensional, then these k particles all move together near s, and the common tangent vector to their trajectories is $\lambda'(\frac{d}{d\lambda}\Phi)|_\lambda$ with $\lambda'$ as in (5.4). The set of all such trajectories that limit to a given zero-dimensional stratum, $\mathcal{E}^\sigma_{k'*}$ as s limits to some $s_*$ can be joined at this stratum to obtain a set of $k'$ continuous, piecewise smooth, trajectories that are defined for s near $s_*$. This follows from (5.2). There is no canonical way to do this joining, but any method will suffice.

At any given value of s, what was just described accounts for at most $\mathfrak{n}_s$ of the particles. The remaining particles are at the point $T_* \in S$. Particles move off or onto the point $T_*$ at values of s for which either of the points $(-\sigma, s)$ or $(\sigma, s)$ are in the closure of $\mathcal{E}^\sigma$. The particles that move on or off $T_*$ and the direction in S that they move are determined by which smooth strata of $\mathcal{E}^\sigma$ have $(-\sigma, s)$ or $(\sigma, s)$ in their closure. The rules for this are essentially the same as those given in the preceding paragraph.

Granted the preceding, let $s \to z(s) \in S$ denote the trajectory of a given particle. The total change let $\Delta z = z(1) - z(0)$, this is the net change in z as s increases from 0 to 1. The intersection number with the point $0 \in S$ of this trajectory is, at most, the least integer that is greater than $\frac{1}{2T} \Delta z$, thus at most $\frac{1}{2T} \Delta z + 1$. Meanwhile, this intersection number is at least the greatest integer less than $\frac{1}{2T} \Delta z$, thus at least $\frac{1}{2T} \Delta z - 1$. As a consequence, the spectral flow for the family $\{\mathcal{L}_s\}_{s \in [0,1]}$ differs by at most $\mathfrak{n}$ from the integral between 0 and 1 of the function

$$\wp(s) = \tfrac{1}{2T} \sum_\varsigma \langle \varsigma, (\tfrac{d}{ds} \mathfrak{q}_s) \varsigma \rangle_\mathbb{H} \, (\tfrac{d}{d\lambda} \Phi)|_{\lambda_\varsigma} \, ,$$

(5.6)

where the sum is over a basis of orthonormal eigenvectors of $\mathcal{L}_s$ whose eigenvalue has absolute value no greater than $\sigma$, and where the notation uses $\lambda_\varsigma$ to denote the eigenvalue of the indicated eigenvector.

The estimates derived below for the spectral flow are obtained by deriving suitable estimates for the function $s \to \wp(s)$ and upper bounds for the number $\mathfrak{n}$.

**c) An upper bound on $\mathfrak{n}$.**

Considered here is a generic sort of operator, $\mathcal{L}$, on $C^\infty(M; iT^*M \oplus \mathbb{S} \oplus i\mathbb{I}_\mathbb{R})$ of the form described next. Fix a connection, $A \in \text{Conn}(E)$, and a hermitian endomorphism, $\Theta$, of $iT^*M \oplus \mathbb{S} \oplus i\mathbb{I}_\mathbb{R}$. Write the covariant derivative on $iT^*M \oplus \mathbb{S} \oplus i\mathbb{I}_\mathbb{R}$ as $\nabla$. This derivative is defined using the Riemannian connection and the connection A. Assume that $r \geq 1$ and that

$$|B_A| + r^{1/2}|\Theta| + |\nabla \Theta| \leq cr \, .$$

(5.7)

The operator $\mathcal{L}$ has the form $\mathcal{L} = \mathbb{L} + \Theta$, where $\mathbb{L}$ sends a given triple, $(b, \eta, \phi)$ to the triple whose respective components in $iT^*M$, $\mathbb{S}$ and $i\mathbb{I}_\mathbb{R}$ are

- $*db - d\phi$ .
- $D_A \eta$ .
- $*d*b$ .

(5.8)

With $\mathcal{L}$ understood, consider:

**Proposition 5.2:** *Fix a constant c and there exists a constant $\kappa$ with the following significance: Define $\mathcal{L}$ as above such that (5.7) holds. Given $\sigma \geq 0$, let $\mathfrak{n}(\sigma)$ denote the* *number of linearly independent eigenvectors of $\mathcal{L}$ whose eigenvalue has absolute value no greater than $\sigma$. For any $R \geq 0$, the number $\mathfrak{n}(Rr^{1/2})$ is bounded by $\kappa\, r^{-3/2}(R^3 + 1)$.*

**Proof of Proposition 5.2**: The proof that follows uses the heat equation for the operator $\mathcal{L}^2$. The idea follows a strategy introduced by [CL]. To start the story, let $j = (b, \eta, \phi) \in C^\infty(M; iT^*M \oplus \mathbb{S} \oplus i\mathbb{I}_\mathbb{R})$, and note that $\mathcal{L}^2 j$ has the form

$$\mathcal{L}^2 j = \nabla^\dagger \nabla j + \mathcal{R}_1 \nabla j + \mathcal{R}_0 j,$$

(5.9)

where $\mathcal{R}_1$ and $\mathcal{R}_0$ are endomorphisms with $r^{1/2}|\mathcal{R}_1| + |\mathcal{R}_0| \leq c\,r$. Introduce the heat kernel for $\mathcal{L}^2$; for each $t \geq 0$, this is the bounded operator on $L^2(iT^*M \oplus \mathbb{S} \oplus i\mathbb{I}_\mathbb{R})$ that is given by

$$E_t = \sum_\varsigma e^{-\lambda_\varsigma^2 t}\, \varsigma \otimes \varsigma^\dagger$$

(5.10)

where the sum is over an orthonormal basis of eigenvectors of $\mathcal{L}^2$. It is well known that $E_t$ is trace class for $t > 0$. Let $\mathrm{Tr}(E_t)$ denote the trace of $E_t$ in $L^2(iT^*M \oplus \mathbb{S} \oplus i\mathbb{I}_\mathbb{R})$. Then

$$\mathrm{Tr}(E_t) = \sum_\varsigma e^{-\lambda_\varsigma^2 t} \geq \mathfrak{n}(Rr^{1/2})\, e^{-R^2 r t}$$

(5.11)

This equation provides an upper bound to $\mathfrak{n}(Rr^{1/2})$.

Standard parametrix techniques (see, for example, [Se], [BGM], [BGV]) prove that $E_t$ has an integral kernel that is smooth for $t > 0$. The value of this kernel at a given $(x, y) \in M \times M$ is denoted in what follows by $E_t(x, y)$. In this regard, $E_t(x, y)$ is a homomorphism from $(iT^*M \oplus \mathbb{S} \oplus i\mathbb{I}_\mathbb{R})|_y$ to $(iT^*M \oplus \mathbb{S} \oplus i\mathbb{I}_\mathbb{R})|_x$. With $y$ fixed and $(t, x)$ allowed to vary, this homomorphism obeys the equation

$$\tfrac{\partial}{\partial t} E_t = -\mathcal{L}^2 E_t\,.$$

(5.12)

with initial condition $E_0 = \mathbb{I}\,\delta_y$. Here, $\delta_y$ denotes the delta function measure at $y$ and $\mathbb{I}$ denotes the identity automorphism of $(iT^*M \oplus \mathbb{S} \oplus i\mathbb{I}_\mathbb{R})|_y$. Taking the inner product of both sides of this with $E_t(\cdot, y)$ finds that the function, $f$, of $t$ and $x$ given by $f(\cdot) = |E_t(\cdot, y)|$ obeys (in the weak sense) the inequality

$$\tfrac{\partial}{\partial t} f \leq -d^\dagger d f + c\, r\, f\,.$$

(5.13)

As a consequence, the function $h = f\, e^{-crt}$ obeys the inequality

$$\tfrac{\partial}{\partial t} h \leq -d^\dagger d h\,.$$

(5.14)

Note also that $h \to \sqrt{6}\, \delta_y$ as $t \to 0$. A standard application of the comparison principle for the heat equation (see [BGM], [Pa], [M]) can now be applied to see that

$$h_t(x) \le c_* (t^{-3/2} + 1)\, e^{c_* t}\, \exp(-\text{dist}(x,y)^2/4t)$$

(5.15)

for $t \le 1$. Here, and below, $c_*$ denotes a constant that depends only on the Riemannian metric. Its value will change from appearance to appearance. Granted (5.15), it follows from what has been said that

$$E_t(x, x) \le c_* (t^{-3/2} + 1)\, e^{crt},$$

(5.16)

Thus,

$$\text{Tr}(E_t) \le c_* (t^{-3/2} + 1)\, e^{c_* t}\, e^{crt}.$$

(5.17)

Taking $t = (R^{-2} + 1) r^{-1}$ in (5.11) and (5.17) gives the claim in Proposition 5.2.

**d) Spectral flow when rescaling $\psi$**

The spectral flow between the two versions of $\mathcal{L}$ as defined by different pairs in $\text{Conn}(E) \times C^\infty(\mathbb{S})$, different values of $r$, and different small normed elements in $\mathcal{P}$ is estimated in what follows using a continuous, but only piecewise real analytic family of operators. This subsection considers this family on an interval where the factor that multiplies $\psi$ in the $r$ and $\mu$ version of (3.1).

To this end, fix $(A, \psi) \in \text{Conn}(E) \oplus C^\infty(M; \mathbb{S})$ and $r \ge 0$. Consider the family of operators on $C^\infty(M; iT^*M \oplus \mathbb{S} \oplus i\mathbb{I}_\mathbb{R})$ that is parametrized by $s \in [0, 1]$ and whose member at a given value of $s$ sends $(b, \eta, \phi)$ to

- $*db - d\phi - s 2^{-1/2} r^{1/2} (\psi^\dagger \tau \eta + \eta^\dagger \tau \psi)$,
- $D_A \eta + s 2^{1/2} r^{1/2} (\text{cl}(b)\psi + \phi\psi)$,
- $*d*b + s 2^{-1/2} r^{1/2} (\eta^\dagger \psi - \psi^\dagger \eta)$.

(5.18)

Let $s \to \mathcal{L}_s$ denote this family. The following lemma summarizes most of what is needed about the spectral flow for $\{\mathcal{L}_s\}_{s \in [0,1]}$. The spectral flow for special choices of $(A, \psi)$ considered in the next subsection.

**Lemma 5.3:** *Given* $c > 0$, *there exists a constant* $\kappa$ *with the following significance: Let* $(A, \psi) \in \text{Conn}(E) \times C^\infty(M; \mathbb{S})$ *be such that* $r^{-1}|B_A| + |\psi| + r^{-1/2}|\nabla\psi| \leq c$. *The absolute value of the spectral flow for the family that is depicted in (5.18) is bounded by* $\kappa r^{3/2}$.

*Proof of Lemma 5.3*: To apply the strategy from Section 5b, take the range for $\Phi$ to be $\mathbb{R}$ and $\Phi$ to be the identity. Take $T = r^{1/2}$. Suppose that $s \in [0, 1]$ and that $(\lambda, s)$ is in a smooth stratum of $\mathcal{E}$. Let $\varsigma$ denote an eigenvector of $\mathcal{L}_s$ with eigenvalue $\lambda$ and with unit $L^2$ norm. The number $\lambda'$ given by (5.4) in this case is

$$\lambda' = 2^{1/2} r^{1/2} \int_M b_k (\eta^\dagger \tau^k \psi + \psi^\dagger \tau^k \eta) \ .$$

(5.19)

Granted this, it follows from the assumptions of the lemma that

$$|\lambda'| \leq c' r^{1/2}$$

(5.20)

where $c'$ is a constant that depends only on the constant c. As a consequence, the absolute value of the function $\wp(s)$ that is depicted in (5.6) is no greater than $c'$ times an upper bound for maximal number of linearly independent eigenvectors of $\mathcal{L}_s$ whose eigenvalue has absolute value less than $r^{1/2}$. This being the case, Proposition 5.2 implies that $|\wp(s)| \leq 2c'' r^{3/2}$, where $c''$ also just depends on the constant c. This bound for $|\wp(s)|$ implies the assertion made by Lemma 5.3.

### e) Spectral flow when $(A, \psi)$ is close to $(A_I, (1_\mathbb{C}, 0))$

This subsection constitutes a digression that first proves Lemma 3.3 and then establishes a somewhat stronger version of Lemma 3.3 that is used later. To start, fix a pair $(A, \psi) \in \text{Conn}(E) \oplus C^\infty(M; \mathbb{S})$ and some $r \geq 0$. Let $\mathfrak{L}$ denote the operator that is depicted in (3.1) with $\mathfrak{t} = \mathfrak{s} = 0$. The subsequent arguments in this subsection require the Bochner-Weitzenbock formula for $\mathfrak{L}^2$. To state this formula, fix an element $j = (b, \eta, \phi) \in C^\infty(M; iT^*M \oplus \mathbb{S} \oplus i\mathbb{R})$. If $D_A\psi = 0$, then the respective $iT^*M$, $\mathbb{S}$ and $i\mathbb{R}$ components of $\mathfrak{L}^2 j$ are

- $\nabla^\dagger \nabla b + \text{Ric}(b) + 2r|\psi|^2 b + 2^{-1/2} r^{1/2} (\eta^\dagger \nabla \psi - (\nabla \psi)^\dagger \eta)$
- $D_A^2 \eta - r[(\eta^\dagger \psi - \psi^\dagger \eta)\psi + (\psi^\dagger \tau^k \eta + \eta^\dagger \tau^k \psi)\tau^k \psi] - 2^{1/2} r^{1/2} b \cdot \nabla \psi$,
- $d^\dagger d\phi + 2r|\psi|^2 \phi$.

(5.21)

Here, $\text{Ric}(b)$ is obtained from b by viewing the Ricci curvature tensor of M as an endomorphism of $T^*M$ and using the latter on b. Meanwhile, $b \cdot \nabla \psi$ denotes the effect on b of the endomorphism from T*M to $\mathbb{S}$ that is defined using the metric and $\nabla \psi$. If $D_A \psi$ is not zero, then $\mathfrak{L}^2 \mathfrak{j}$ is the sum of what is written in (5.21) and a term that has the schematic form $r^{1/2} \Bbbk(D_A \psi, \mathfrak{j})$, where $\Bbbk( , )$ is fiber wise bilinear in its two argument and is such that $|\Bbbk(\zeta, \mathfrak{j})| \leq c\, |\zeta||\mathfrak{j}|$ with c a constant that is of r and $(A, \psi)$.

*Proof of Lemma 3.3*: Consider (5.21) in the case where $A = A_I$ and $\psi = (1_{\mathbb{C}}, 0)$. In this case, $\mathfrak{L}^2 \mathfrak{j}$ is

- $\nabla^\dagger \nabla b + \text{Ric}(b) + 2r\, b + r^{1/2}(\eta^\dagger \mathcal{R} - \mathcal{R}^\dagger \eta)$,
- $\mathbb{V}^\dagger \mathbb{V} \eta + \frac{1}{4} R \eta + 2r\, \eta - 2r^{1/2} b \cdot \mathcal{R}$,
- $d^\dagger d\phi + 2r\phi$.

(5.22)

Here, R denotes an endomorphism of $\mathbb{S}$ that depends only on the metric and the contact form. Meanwhile $\mathcal{R}$ denotes the section $\mathbb{V}(1_{\mathbb{C}}, 0)$ of $\mathbb{S} \otimes T^*M$. Note in paticular that $|R| \leq c_1$ and $|\mathcal{R}| \leq c_1$ where $c_1$ is independent of r. Contract both sides of (5.22) with $(b, \eta, \phi)$ and to see that

$$\int_M |\mathfrak{L}\mathfrak{j}|^2 \geq (r - c_2) \int_M |\mathfrak{j}|^2$$

(5.23)

Here, $c_2$ is a constant that is independent of r and $\mathfrak{j}$. The statement made by Lemma 3.3 follows from this last equation.

Consider now the case where $(A, \psi)$ is close to $(A_I, (1_{\mathbb{C}}, 0))$. To make this notion precise, first fix $r > 0$ and $\varepsilon > 0$. Let $A \in \text{Conn}(I_{\mathbb{C}})$ and $\psi = (\alpha, \beta) \in C^\infty(M; \mathbb{S}_I)$, and suppose that this pair is such that the following hold at each point in M:

- $1 - \varepsilon \leq |\alpha| \leq 1 + \varepsilon$ *and* $|\beta| \leq r^{-1/2}\varepsilon$,
- $|\nabla \alpha| \leq \varepsilon\, r^{1/2}$ *and* $|\nabla' \beta| \leq \varepsilon$,
- $|B_A| \leq \varepsilon\, r$.

(5.24)

**Lemma 5.4**: *There exist constants $\varepsilon_0 > 0$ and $r_0 \geq 1$ with the following significance: Suppose that $r > r_0$ and that $(A, \psi = (\alpha, \beta)) \in \text{Conn}(I_{\mathbb{C}}) \times C^\infty(M; \mathbb{S}_I)$ obeys (5.24) with $\varepsilon < \varepsilon_0$. Then the operator $\mathfrak{L}$ as given in (3.1) with $\mathfrak{t} = \mathfrak{s} = 0$ has no kernel. Moreover, if K has torsion first Chern class, then there is zero spectral flow between the latter operator and the $\mathfrak{t} = \mathfrak{s} = 0$ version of (3.1) that is defined by $(A_I, (1_{\mathbb{C}}, 0))$.*

***Proof of Lemma 5.4***: To see that $\mathcal{L}$ has no kernel, use (5.24) with the $D_A \psi \neq 0$ version of the Weitzenboch formula in (5.21) to see that

$$\int_M |\mathcal{L}j|^2 \geq ((1 - c_1(\varepsilon + r^{-1}) r \int_M |j|^2 .$$

(5.25)

Here $c_1$ is a constant that is independent of $(A, \psi)$. The latter inequality proves that the kernel of $\mathcal{L}$ is trivial when r is larger than some fixed $r_0$ and $\varepsilon c_1 < \frac{1}{4}$

To see that there is no spectral flow in the case where K has torsion first Chern class, note first that if $u \in C^\infty(M; S^1)$, then there is zero spectral flow between the respective versions of $\mathcal{L}$ that are defined by the two pairs $(A, \psi)$ and $(A - u^{-1}du, u\psi)$. This being the case, there exists a unique choice for u that makes $\alpha = |\alpha| 1_{\mathbb{C}}$. What with the preceding remarks, no generality is lost by assuming henceforth that $\alpha = |\alpha| 1_{\mathbb{C}}$.

Now write $A = A_I + \hat{a}$ with $\hat{a}$ a section of $iT^*M$. It then follows from the bound in (5.24) on $|\nabla \alpha|$ that

$$|\hat{a}| \leq 2\varepsilon \, r^{1/2}.$$

(5.26)

Granted this, introduce, for each $s \in [0, 1]$, the pair $(A^s, \psi^s)$, where $A^s = A_I + s\hat{a}$ and where $\psi^s = (\alpha^s, \beta^s)$ with $\alpha^s = (1 - s(1 - |\alpha|)) 1_{\mathbb{C}}$ and $\beta^s = s\beta$. Then $(A^s, \psi^s)$ obey the conditions in (5.24) with $2\varepsilon$ replacing $\varepsilon$. Hence the $(A^s, \psi^s)$ and $\mathfrak{t} = \mathfrak{s} = 0$ version of $\mathcal{L}$ obeys (5.25) with $2\varepsilon$ replacing $\varepsilon$. As a consequence, all of these versions of $\mathcal{L}$ have trivial kernel, and so there is zero spectral flow between the $(A, \psi)$ version of $\mathcal{L}$ and the $(A_I, (1_{\mathbb{C}}, 0))$ version of $\mathcal{L}$.

**f) Spectral flow for the Dirac operator**

This subsection considers the Dirac operator on $C^\infty(M; \mathbb{S})$ as defined by connections on $\det(\mathbb{S})$ and the spectral flow for a path of such operators. Note that the Dirac operator here is viewed as a $\mathbb{C}$-linear operator, and so eigenspaces are viewed as vector spaces over $\mathbb{C}$.

To put things into a slightly more general framework, make no assumption about a contact 1-form on M or the first Chern class of $\det(\mathbb{S})$. Assume only that M is a compact, oriented Riemannian manifold with a chosen $Spin_{\mathbb{C}}$ structure. Let $\mathbb{S}$ denote the corresponding $\mathbb{C}^2$ bundle. Let A denote a given connection on $\det(\mathbb{S})$, and suppose that constants $c \geq 1$ and $r \geq 1$ have been given and that the following conditions hold:

- $|B_A| \leq c\, r$ .

- $|\nabla B_A| \le c\, r^{3/2}$.

(5.27)

Let $A_0$ denote a fixed, fiducial connection on $\det(\mathbb{S})$.

**Proposition 5.5**: *Given c and the connection $A_0$ on $\det(\mathbb{S})$, there is a constant $\kappa$ with the following significance: Suppose that $r \ge 1$ and that A is a connection on $\det(\mathbb{S})$ that obeys the conditions in (5.27). Write $A = A_0 + \hat{a}_A$. Then the spectral flow along a path of Dirac operators that starts at that defined as in (2.1) by $A_0$ and ends at that defined as in (2.1) by A differs from $-\frac{1}{32\pi^2} \int_M \hat{a}_A \wedge d\hat{a}_A - \frac{1}{16\pi^2} \int_M \hat{a}_A \wedge *B_{A_0}$ by at most $\kappa\, r^{15/8} (\ln r)^{3/2}$.*

Note that the factor of $\frac{1}{32\pi^2}$ that appears here leads to the factor of $\frac{1}{4\pi^2}$ that appears in Propositions 5.1 and 3.4. A factor of 4 appears because the connection in Proposition 5.5 is defined on $\det(\mathbb{S})$ while those in Propositions 5.1 and 3.4 are defined on E. The extra factor of 2 that appears is due to the fact that the spectral flow in the cSWF context deals with operators that are $\mathbb{R}$-linear rather than $\mathbb{C}$-linear.

*Proof of Proposition 5.5*: The first point to make is that there exists a map $u: M \to S^1$ such that $A - u^{-1}du$ can be written as $A_0 + \hat{a}$ where $d*\hat{a} = 0$ and

$$|\hat{a}| \le c'r \quad \text{and} \quad |\nabla \hat{a}| \le c'r^{3/2}.$$

(5.28)

Here, $c'$ depends only on the constant c and $A_0$. This is proved as follows: First, write $A = A_0 + \hat{a}_A$. Next, fix a smooth map, $u_1: M \to S^1$ with the property that integral of the real valued 1-form $i(\hat{a}_A - u_1^{-1}du_1)$ around each of fixed set of basis elements for $H_1(M; \mathbb{Z})$ lies in the interval $[0, 2)$. This guarantees that the $L^2$-orthogonal projection of $\hat{a}_A - u_1^{-1}du_1$ onto the space of harmonic 1-forms has norm bound that depends only the metric. Granted that such is the case, then Hodge theory finds a unique, smooth and homotopically trivial map, $u_2: S^1 \to M$ such that $\hat{a} = \hat{a}_A - u_1^{-1}du_1 - u_2^{-1}du_2$ is coclosed. Note that $\hat{a}$ and $\hat{a}_A - u_1^{-1}du_1$ have the same orthogonal projection to the space of harmonic 1-forms. The bounds in (5.28) follow by exploiting standard estimates for the Green's kernel for the operator $*d$ on the vector space of co-closed 1-forms. (In fact, some care with the estimates finds $|\nabla \hat{a}| \le c'\, r\ln(r+1)$.)

The change in the spectral flow between the respective Dirac operators defined by connections A and $A - u^{-1}du$ is the same as the change in the respective A and $A - u^{-1}du$ versions of $-\frac{1}{32\pi^2} \int_M \hat{a}_{(\cdot)} \wedge d\hat{a}_{(\cdot)} - \frac{1}{16\pi^2} \int_M \hat{a}_{(\cdot)} \wedge *B_{A_0}$. Thus, it is sufficient to consider the case where $A = A_0 + \hat{a}$ with $\hat{a}$ obeying the bounds in (5.28).

For each $s \in [0, 1]$, set $A^s = A_0 + s\hat{a}$. The spectral flow will be estimated for the family $\{\mathcal{L}_s\}_{s \in [0,1]}$ where $\mathcal{L}_s$ is the Dirac operator that is defined as in (2.1) by the connection $A^s$. Thus, $\mathcal{L}_s = \mathcal{L}_0 + s\frac{1}{2}\text{cl}(\hat{a})$. Note that the factor of $\frac{1}{2}$ appears here because the connection $A^s$ now denotes a connection on $\det(\mathbb{S})$.

This is a family of self-adjoint, unbounded operators on $L^2(M; \mathbb{S})$ whose $s = 0$ member is the Dirac operator defined by (2.1) using $A_0$ and whose $s = 1$ member is that defined using $A$. The strategy that is described in Section 5b is used to estimate the spectral flow for this family. To apply this strategy, fix $t \in (0, r^{-1})$; a specific choice is made shortly. With t chosen, the range space for Section 5b's diffeomorphism $\Phi$ is the open interval $(-(\frac{\pi}{4t})^{1/2}, (\frac{\pi}{4t})^{1/2})$; and $\Phi$ itself is given by

$$\Phi(\lambda) = \int_0^\lambda e^{-\rho^2 t} d\rho \ .$$

(5.29)

Fix $R \geq 1$ and set $T = \Phi(Rt^{-1/2})$. A specific choice for R is also made below. Note for reference later that

$$|t^{1/2}T - (\tfrac{\pi}{4})^{1/2}| \leq \tfrac{1}{2R} e^{-R^2} \ .$$

(5.30)

The function depicted in (5.6) for this set up is

$$\wp(s) = \tfrac{1}{4T} \Sigma_\varsigma (\int_M \varsigma^\dagger \text{cl}(\hat{a}) \varsigma) \, e^{-\lambda_\varsigma^2 t}$$

(5.31)

where the sum in question is indexed by an orthonormal basis of eigenvectors of $\mathcal{L}_s$ whose eigenvalue has absolute value no greater than $Rt^{-1/2}$. The strategy for estimating $\wp$ exploits the fact that sum on the right hand side of (5.31) looks much like the trace on $L^2(M; \mathbb{S})$ of the composition of the multiplication operator $\tfrac{1}{4T}\text{cl}(\hat{a})$ with the heat kernel for $\mathcal{L}_s^2$, this the operator $E_t$ on $L^2(M; \mathbb{S})$ that is given by the expression in (5.10) with the sum indexed by an orthonormal basis of eigenvectors of $\mathcal{L}_s^2$. To make something of this resemblence, introduce $\Pi \subset L^2(M; \mathbb{S})$ to denote the span of the eigenvectors of $\mathcal{L}_s$ whose eigenvalue has absolute value no greater than $Rt^{-1/2}$. With $\Pi$ understood, note that

$$\Sigma_{\varsigma \notin \Pi} |\int_M \varsigma^\dagger \text{cl}(\hat{a}) \varsigma | \, e^{-\lambda_\varsigma^2 t} \leq c'r \Sigma_{\varsigma \notin \Pi} \, e^{-\lambda_\varsigma^2 t} \ ,$$

(5.32)

as can be seen with the help of (5.28). Here, the sum is over an orthonormal basis of eigenvectors of $\mathcal{L}_s$ whose eigenvalue has absolute value greater than $Rt^{-1/2}$. Let $\mathfrak{n}(\cdot)$ denote the function that is defined by Proposition 5.2 for $\mathcal{L} = \mathcal{L}_s$. It then follows that the sum on the right hand side of (5.32) is no larger than

$$rc' \Sigma_{m=1,2,\ldots} \mathfrak{n}(Rmt^{-1/2}) \, e^{-m^2 R^2} \leq c_1 c' r t^{-3/2} e^{-R^2/2} \ .$$

(5.33)

Here, $c_1$ is a constant that is independent of R, t, r and A. It follows from (5.31) and (5.33) that

$$\wp(s) = \tfrac{1}{4T} \, \mathrm{Tr}(\mathrm{cl}(\hat{a}) \, E_t) + \mathfrak{r},$$

(5.34)

where $E_t$ again denotes the time t heat kernel for $\mathcal{L}_s$, where $\mathrm{Tr}(\cdot)$ denotes the trace of the indicated operator on $L^2(M; \mathbb{S})$, and where

$$|\mathfrak{r}| \leq 2 \, c' \, c_2 \, t^{-1} \, r \, e^{-R^2/2}.$$

(5.35)

As with $c_1$, the constant $c_2$ is also independent of r, t, R, and A.

The task now is to provide an estimate with controlled errors for $\mathrm{Tr}(\mathrm{cl}(\hat{a}) \, E_t)$. This is done by using a small t approximation for $E_t$. The following lemma provides a useable estimate.

**Lemma 5.6**: *Let* $p \in M$. *Then*

$$E_t(p,p) = \tfrac{1}{2} \, (\tfrac{1}{4\pi})^{3/2} \, (\tfrac{1}{t})^{1/2} (\mathrm{cl}(B_{A^s})|_p + \mathfrak{w})$$

*where* $|\mathfrak{w}| \leq c_0 c \, r \, (rt)^{1/2}$. *Here, $c_0$ is a constant that depends only on the metric and the connection $A_0$; and c is the constant in (5.28).*

*Proof of Lemma 5.6*: Fix attention on a point, $p \in M$, and fix a Gaussian coordinate chart centered at p. This is a diffeomorphism $\varphi$, from the ball $U \subset \mathbb{R}^3$ of some radius $\rho > 0$ centered at the origin to M with $\varphi(0) = p$, and with the property that the Euclidean distance in U from the origin is the same as that defined by the pull-back of the metric from M. In particular, if m denotes the latter and if it is viewed as a symmetric, $3 \times 3$ matrix valued function on U, then

$$|m - \mathbb{I}| \leq c_* |x|^2 \quad \mathrm{and} \quad |dm| \leq c_* |x|,$$

(5.36)

where $\mathbb{I}$ here denotes the identity $3 \times 3$ matrix. Here, and in what follows, $c_* \geq 1$ denotes a constant that depends only on the Riemannian metric. It's precise value is allowed to change between successive appearances. The radius $\rho$ is determined by the metric and can be assumed to be independent of the point chosen in M. The Euclidean coordinates on B are denoted by $(x^1, x^2, x^3)$. To simplify notation, use $\varphi$ to identify B with $\varphi(B)$.

Use parallel transport by the connection $A^s$ to trivialize the bundle E over U, and use this trivialization with the coordinate chart's trivialization of the frame bundle of M over U to trivialize $\mathbb{S}$ over U.  With respect to this trivialization of E, the connection $A^s$ pulls back as an i-valued 1-form which appears when written with respect to the basis of coordinate differentials $\{dx^1, dx^2, dx^3\}$ as $\nu = \nu_j dx^j$ where here, and in what follows, repeated indices from the set $\{1, 2, 3\}$ are implicitly summed.  Note in particular that

$$\nu_j|_0 = 0 \quad and \quad \nu_j x^j = 0 .$$
(5.37)

With the trivialization of $\mathbb{S}$ implicit, the restriction to $U \times U$ of the integral kernel for $E_t$ is a function, $(x, y) \to E_t(x, y)$, on $U \times U$ with values in $End(\mathbb{C}^2)$.  As indicated by (5.15), this function obeys

$$|E_t(x, y)| \leq c_* ((\tfrac{1}{t})^{3/2} + 1) e^{-|x-y|^2/4t} e^{(c_* + cr)t},$$
(5.38)

where $c_*$ depends only on the Riemannian metric.  Moreover, if $|x| \leq \tfrac{1}{2}$ and if $q \in M-B$, then the value of the heat kernel at time t on either $(x, q) \in U \times M$ or $(q, x) \in M \times U$ is bounded in absolute value by $c_* e^{-\rho^2/2t} e^{(c_* + cr)t}$.

Let $\mathfrak{h}$ denote the $End(\mathbb{C}^2)$ valued function on $\mathbb{R} \times U$ given by $\mathfrak{h}(t, x) = E_t(x, 0)$.  This function obeys an equation of the form

$$\tfrac{\partial}{\partial t} \mathfrak{h} = \partial_j \partial_j \mathfrak{h} + \tfrac{1}{2} cl(B) \mathfrak{h} + V \mathfrak{h},$$
(5.39)

where $*B = d\nu$ is the curvature 2-form for the connection $A^s$ and

$$V\mathfrak{h} = (\delta_{ij} + \gamma_{ij})(2\nu_i \partial_j \mathfrak{h} + \partial_i \nu_j \mathfrak{h} + \nu_i(\nu_j + \gamma_j)\mathfrak{h} + 2\Gamma_i \partial_j \mathfrak{h}) + \gamma_{ij}\partial_i\partial_j \mathfrak{h} + + \Gamma_0 \mathfrak{h} .$$
(5.40)

Here, $\{\gamma_{ij}\}_{i,j=1,2,3}$, $\{\gamma_j, \Gamma_j\}_{j=1,2,3}$ and $\Gamma_0$ are $End(\mathbb{C}^2)$ valued functions on U that are determined by the Riemannian metric.  Note in particular that $\gamma_j = \Gamma_j|_0 = 0$ since both are linear combination of the metric's Christoffel symbols.  In addition $|\gamma_{ij}| \leq c'|x|^2$ with $c'$ depending only on the metric.  As $t \to 0$, the function $\mathfrak{h}(t, x)$ converges as an $End(\mathbb{C}^2)$ valued measure to $\mathbb{I} \delta_0$, where $\mathbb{I}$ now denotes the identity endomorphism of $\mathbb{C}^2$ and $\delta_0$ denotes the Dirac measure at $0 \in U$.

As the author learned from an unpublished paper by Tom Parker [Pa] (see also [BGM], [BGV]), there is a nice method of obtaining a controlled estimate for $\mathfrak{h}$ at small t.  To set the stage, introduce the function on $(0, \infty) \times (U \times U)$ that sends $(t, (x, y))$ to

$$K_t(x, y) = (\tfrac{1}{4\pi t})^{3/2} e^{-|x-y|^2/4t} .$$
(5.41)

Let $\chi\colon [0, \infty) \to [0, 1]$ denote a smooth, non-increasing function that equals 1 on $[0, \frac{1}{4}]$ and vanishes on $[\frac{1}{2}, \infty)$. Set $\chi_\rho$ to denote the function with compact support on U whose value at a given point x is $\chi(|x|/\rho)$. Let $h = \chi_\rho \mathfrak{h}$. This End($\mathbb{C}^2$) valued function obeys the equation

$$\tfrac{\partial}{\partial t} h = \partial_j \partial_j h + \tfrac{1}{2} \operatorname{cl}(B) h + V h - 2(\delta_{ij} + \gamma_{ij})\partial_i \chi_\rho \, \partial_j \mathfrak{h} - (\delta_{ij} + \gamma_{ij})\partial_i \partial_j \chi_\rho \, \mathfrak{h} - 2\gamma_{ij}(\Gamma_i + v_i)\partial_j \chi_\rho \mathfrak{h}. \tag{5.42}$$

Note in this regard that $\Gamma_j \partial_j \chi_\rho = v_j \partial_j \chi_\rho = 0$ because $\chi_\rho$ depends solely on the radial distance $|x|$ and $v_j x^j$ and $\Gamma_j x^j$ are both zero.

One then writes

$$h(t,x) = K_t(x, 0)\mathbb{I} + \tfrac{1}{2} \int_0^t \!\! \int_U K_{t-\sigma}(x,y)\operatorname{cl}(B)_y h(\sigma, y) d^3 y \, d\sigma + \mathfrak{R}(t,x), \tag{5.43}$$

where

$$\mathfrak{R}(t, x) = \int_0^t \!\! \int_U K_{t-\sigma}(x, y)(\overleftarrow{V}^\dagger \chi_\rho + 2\overleftarrow{\partial}_j \partial_j \chi_\rho + \partial_j \partial_j \chi_\rho)_y \mathfrak{h}(\sigma, y) d^3 y \, d\sigma. \tag{5.44}$$

Here, the arrows over $V^\dagger$ and $\partial_j$ indicates that one derivative acts on the term to their left. Bounds on $\mathfrak{R}$ are the next order of business. To this end, note that the terms with derivatives of $\chi_\rho$ are supported where $|y| \geq \tfrac{1}{4}\rho$, and thus where

$$|\mathfrak{h}(\sigma)| \leq c_* ((\tfrac{1}{\sigma})^{3/2} + 1) e^{-\rho^2/64\sigma} e^{(c_* + crt)} \tag{5.45}$$

Indeed, this follows from (5.15). As a consequence, the terms in (5.44) that involve derivatives of $\chi_\rho$ have norms that are no greater than

$$c_0 c \, e^{-\rho^2/c_* t} e^{(c_* + cr)t}. \tag{5.46}$$

where c is the constant from (5.27) and $c_0$ here, and in what follows, depends only on the Riemannian metric and the curvature of the connection $A_0$. Note that the different appearances of $c_0$ have distinct values. The remaining term can be bounded using (5.15) to bound $\mathfrak{h}(\sigma, y)$ and (5.37) to obtain bounds on v. As for the latter, the equation $dv = *B$ and (5.37) can be used to write v in terms of B and thus see that

$$|v| \leq c_0 c \, r |x| \quad \text{and} \quad |\partial_j v_j| \leq c_0 c \, r^{3/2} |x| \quad \text{and} \quad |\partial_i v_j| \leq c_0 c \, (r + r^{3/2}|x|), \tag{5.47}$$

where c is the constant in (5.27). Granted these last bounds, it then follows that the term in (5.44) with $\overleftarrow{V}^\dagger \chi_\rho$ has norm no greater than

$$c_0 c (t + r|x|t^{1/2} + (tr)^{3/2} + r^2 t^2)(\tfrac{1}{t})^{3/2} e^{-|x|^2/16t},$$

(5.48)

Turn next to the term in (5.43) with cl(B). What with (5.27) and (5.15), the norm of this term can be seen to be no greater then

$$c_0 c r t \, (\tfrac{1}{t})^{3/2} e^{-|x|^2/4t}.$$

(5.49)

An estimate for this term is needed as well as an upper bound on its norm. To obtain such an estimate, write

$$+ \tfrac{1}{2} \int_0^t \int_U K_{t-\sigma}(x,y) cl(B)_y K_\sigma(y,0) \, d^3y d\sigma + \mathfrak{P}_t(x),$$

(5.50)

where $\mathfrak{P}_t(x)$ is obtained from the left most term in (5.50) by replacing $K_\sigma(y, 0)$ with the term $h(\sigma, y) - K_\sigma(y, 0)\mathbb{I}$. Of interest here is the value of the left most term at $x = 0$. Since $|B|_y - B|_0| \le c_0 c \, r^{3/2} |y|$, this left most term at $x = 0$ has the form

$$(\tfrac{1}{2} cl(B)_0 t + \mathfrak{v})(\tfrac{1}{4\pi t})^{3/2},$$

(5.51)

where $|\mathfrak{v}| \le c_0 c (rt)^{3/2}$. Meanwhile, the norm of the term $\mathfrak{P}_t(x)$ at $x = 0$ can be bounded using (5.46), (5.48) and (5.49). In particular,

$$|\mathfrak{P}_t(0)| \le c_0 \, c \, (r^2 t^2) (\tfrac{1}{t})^{3/2}.$$

(5.52)

The assertion of Lemma 5.6 follows directly from the estimate in (5.51) with the bounds derived for the norms of $\mathfrak{R}_t$ and $\mathfrak{P}_t$.

It follows directly from Lemma 5.6 with (5.34) and (5.35) that

$$\wp(s) = -\tfrac{1}{16\pi^2} *(\hat{a} \wedge *B_{A^s}) + \mathfrak{z},$$

(5.53)

where

$$|\mathfrak{z}| \le c_0 \, c \, r^2 (r^2 + (rt)^{1/2} + e^{-R^2/2} (tr)^{-1}).$$

(5.54)

Now take $t = r^{-5/4}$ and $R = 2(\ln r)^{1/2}$. According to Section 5b, the spectral flow in question differs from $\int_0^1 \wp(s) ds$ by no more than $\mathfrak{n} = \mathfrak{n}(Rt^{-1/2})$. Given that Proposition 5.2 finds $n(Rt^{-1/2}) \le \kappa R^3 t^{-3/2}$, so $n$ in this case is bounded by $8\kappa r^2 ((\ln r)^{3/2} r^{-1/8})$. Meanwhile, the right hand side of (5.54) for this choice of t and R is no greater than $c_0 c\, r^2\, r^{-1/8}$. These bounds and (5.53) lead directly to the assertion made by Proposition 5.5.

### g) Spectral flow when $A = A_E - \frac{i}{2} r a$

This last subsection adds something to the statement of Proposition 5.5 for the case where the $\text{Spin}_{\mathbb{C}}$ structure is such that $\mathbb{S} = E \oplus K^{-1}E$ where the splitting is defined using the contact 1-form a. There is no need to assume in this subsection that E has torsion first Chern class. The following lemma is the focus of this subsection.

**Lemma 5.7**: *Given c, there is a constant $\kappa \ge 8$ with the following significance: Suppose that $E \to M$ is a complex line bundle and $A_E$ is a hermitian connection on E whose curvature has norm bounded everywhere by c. Let $\mathbb{S} = E \oplus K^{-1}E$, and define the family of Dirac operators $\{D_r \colon C^{\infty}(M; \mathbb{S}) \to C^{\infty}(M; \mathbb{S})\}_{r \in [0,\infty)}$ using the connection $A_E - \frac{i}{2} r a$. Suppose that $r > \kappa$ and that $\varsigma$ is an eigenfunction of $D_r$ with eigenvalue $\lambda$ with absolute value less than $(\frac{1}{6} r)^{1/2}$. Then*

$$-\tfrac{i}{2} \int_M \varsigma^{\dagger} \mathrm{cl}(a) \varsigma \ge \tfrac{1}{2}(1 - 8 r^{-1}).$$

*Thus, all eigenvalues that cross zero as r increases from $\kappa$ cross from below to above.*

***Proof of lemma 5.7:*** The Weitzenboch formula for $D_r^2$ asserts that

$$D_r^2 \eta = \nabla^{\dagger} \nabla \eta + i r \,\mathrm{cl}(a)\eta + \mathrm{cl}(B_{A_E}) \eta + \tfrac{1}{4} \mathcal{R}\, \eta,$$

(5.55)

where $\mathcal{R}$ denotes an endomorphism whose norm has an r-independent bound. Now suppose that $\varsigma$ is an eigenvector of $D_r$ with eigenvalue $\lambda$, and write $\varsigma = (\varsigma_0, \varsigma_1)$ with respect to the splitting $\mathbb{S}_E = E \oplus K^{-1}E$. Take the $L^2$-inner product of the expression in (5.55) first with $(\varsigma_0, 0)$ and then with $(0, \varsigma_1)$ and integrate by parts to obtain

- $\lambda^2 \|\varsigma_0\|_2^2 = \|\nabla \varsigma_0\|_2^2 - r \|\varsigma_0\|_2^2 + \langle \varsigma_0, R_0 \varsigma_0 \rangle_2 + \langle \varsigma_0, \theta \varsigma_1 \rangle_2 + \langle \varsigma_0, \Gamma \nabla' \varsigma_1 \rangle_2$,
- $\lambda^2 \|\varsigma_1\|_2^2 = \|\nabla' \varsigma_1\|_2^2 + r \|\varsigma_1\|_2^2 + \langle \varsigma_1, R_1 \varsigma_1 \rangle_2 + \langle \theta \varsigma_1, \varsigma_0 \rangle_2 + \langle \Gamma \nabla' \varsigma_1, \varsigma_0 \rangle_2$.

(5.56)

Here, $\langle\, ,\,\rangle_2$ denotes the $L^2$ inner product, and $R_0, R_1, \theta,$ and $\Gamma$ are homomorphisms that are determined solely by the Riemannian metric and the curvature 2-form of $A_E$. In particular, it follows from the second line in (5.56) that there is a constant, c, that depends solely on the metric and the curvature form of $A_E$, and is such that

$$\|\nabla'\varsigma_1\|_2^2 + (\tfrac{1}{2}r - \lambda^2)\|\varsigma_1\|_2^2 \leq \kappa \|\varsigma_0\|_2^2 .$$
(5.57)

In particular that if $\lambda^2 < \tfrac{1}{6}r$ and $r > 12c$, then this last equation implies that

$$\|\varsigma_1\|_2^2 \leq 4r^{-1}\|\varsigma_0\|_2^2 \quad and \quad \|\varsigma_0\|_2^2 \geq 1 - 4r^{-1} .$$
(5.58)

To finish the story, note that

$$-\tfrac{i}{2}\int_M \varsigma^\dagger cl(a)\varsigma = \tfrac{1}{2}(\|\varsigma_0\|_2^2 - \|\varsigma_1\|_2^2) .$$
(5.59)

and so if $|\lambda|^2 < \tfrac{1}{6}r$ and $r > 12c$, then the expression on the right hand side of (5.59) is no less than $\tfrac{1}{2}(1 - 8r^{-1})$

The final assertion of the lemma follows from the fact that the number $\lambda'$ that appears in (5.6) for the family $s \to \mathcal{L}_s$ with $\mathcal{L}_s = D_{r=s}$ is given by the expression on the right hand side of (5.59).

## 6) The behavior of solutions to the Seiberg-Witten equations

The purpose of this section is to tie up some loose ends with regards to the assertions made in Section 2 about the behavior of solutions to certain versions of (2.4). In particular, proofs are given here of Theorem 2.1 and Lemmas 2.2-2.4.

### a) Proof of Lemma 2.2

Because $D_A\psi = 0$, so $D_A^2\psi = 0$. The Bochner-Weitzenboch formula for $D_A^2$ finds that

$$D_A^2\psi = \nabla^\dagger\nabla\psi - r cl(\psi^\dagger\tau\psi - ia)\psi + \tfrac{1}{4}R\psi - icl(\mathfrak{u}))\psi = 0 ,$$
(6.1)

where $\mathfrak{u} = *d\mu + \varpi$ denotes the perturbation term in (2.4) and R now denotes the scalar curvature of the metric on M. Contract this equation with $\psi$ to see that

$$\tfrac{1}{2}d^\dagger d|\psi|^2 + |\nabla\psi|^2 + r|\psi|^2(|\psi|^2 - 1 - c_*r^{-1}) \leq 0$$
(6.2)

where $c_*$ depends only on the infimum of the scalar curvature and the maximum of $|\mathfrak{u}|$. Note that the latter has a bound that depends only on Lemma 2.2's constant c. Granted (6.2), the maximum principle requires that

$$|\psi|^2 \leq 1 + c_* r^{-1} .$$
(6.3)

This last equation gives the assertion made by the first bullet of Lemma 2.2

To continue, contract (6.1) first with $(\alpha, 0)$ and then with $(0, \beta)$ to see that

- $\frac{1}{2} d^\dagger d |\alpha|^2 + |\nabla \alpha|^2 - r(1 - |\alpha|^2 - |\beta|^2)|\alpha|^2 + \mathfrak{r}_0(\alpha, \beta) + \mathfrak{r}_1(\alpha, \nabla'\beta) + \mathfrak{r}_2|\alpha|^2 = 0$.
- $\frac{1}{2} d^\dagger d |\beta|^2 + |\nabla'\beta|^2 + r(1 + |\alpha|^2 + |\beta|^2)|\beta|^2 + \mathfrak{r}'_0|\beta|^2 + \mathfrak{r}'_1(\beta, \nabla\alpha) + \mathfrak{r}'_2(\alpha,\beta) = 0$.

(6.4)

Here, $\mathfrak{r}_0, \mathfrak{r}_1, \mathfrak{r}_2$ and their primed counterparts depend solely on the Riemannian metric. Introduce $w = (1 - |\alpha|^2)$. The top equation in (6.4) implies the following equation for w:

$$\tfrac{1}{2} d^\dagger dw + r\, w - |\nabla \alpha|^2 - r(w^2 + |\beta|^2 |\alpha|^2) - \mathfrak{r}_0(\alpha, \beta) - \mathfrak{r}_1(\alpha, \nabla'\beta) - \mathfrak{r}_2|\alpha|^2 = 0.$$
(6.5)

Equation (6.3), the bottom equation in (6.4) and (6.5) have the following consequence: There are constants $c_1$ and $c_2$ that depend solely on the Riemannian metric and the constant c, and are such that

$$d^\dagger d(|\beta|^2 - c_1 r^{-1} w - c_2 r^2) + r(|\beta|^2 - c_1 r^{-1} w - c_2 r^2) \leq 0 .$$
(6.6)

An application of the maximum principle to this last equation gives the second bullet of Lemma 2.2.

**b) Proof of Lemma 2.3**

No generality is lost by assuming $r \geq 1$. Fix a point $p \in M$ and fix a Gaussian coordinate chart centered at p. Let $\rho$ denote the radius of the ball in $\mathbb{R}^3$ on which the coordinate map, $\varphi$, is defined, and use $\varphi$ to identify this ball with its $\varphi$-image in M. Let $(x^1, x^2, x^3)$ denote the coordinates for this ball. Let $y^k = r^{1/2} x^k$ for each k, and view $(\alpha, \beta)$ as functions of $y = (y^1, y^2, y^3)$. Likewise, view the connection A using these coordinates. Use the coordinates $y = (y^1, y^2, y^3)$ for $\mathbb{R}^3$ and let $U \subset \mathbb{R}^3$ denote the ball where $|y| \leq 8$. The equations in (2.4) on this ball, when written using the y-coordinates read

- $B_{Ak} = (\psi^\dagger \tau^k \psi - i a_k) + r^{-1} i \mathfrak{u}_k$,
- $\tau^k \nabla^{(y)}_k \psi = 0$,

(6.7)

where $B_{Ak}$, $a_k$ and $\mathfrak{u}_k$ are the respective components of their namesake 1-forms when the latter are written as linear combinations of $dy^1$, $dy^2$ and $dy^3$. Meanwhile, $\tau^k = cl(dy^k)$ and $\nabla^{(y)}_k$ is the covariant derivative with respect to $y^k$. Here again, $\mathfrak{u} = *d\mu + \varpi$.

Granted that $|\psi|$ has an r-independent upper bound, standard elliptic regularity arguments (very much simpler versions of the sort found in Chapter 6 of [Mo]) can be applied to the equations in (6.7) on the ball where $|y| \leq 4$. These find, for each q, a constant, $c_q$ that depends only on the Riemannian metric and the $C^{q+2}$ norm of u, and is such that $|(\nabla^{(y)})^q \psi| \leq c_q$. When $\psi$ is viewed as a function of the Gaussian coordinates x, this last bound says $|\nabla^q \psi| \leq c_q r^{q/2}$. The assertion made by the first bullet of Lemma 2.3 follows directly from the latter bound.

To obtain the assertion made by the second bullet of the lemma, again view $\alpha$ and $\beta$ as functions of y. Project the equation in (6.1) onto the E´ summand in $\mathbb{S}$ to obtain the following equation:

$$\nabla^{(y)\dagger}\nabla^{(y)}\beta + (1 + |\alpha|^2 + |\beta|^2)\beta + r^{-1}\mathfrak{r}'_0 \beta + r^{-1/2}\mathfrak{r}'_1 \nabla^{(y)}\alpha + r^{-1}\mathfrak{r}'_2 \alpha = 0 \,. \tag{6.8}$$

By virtue of the uniform bound for $\psi$ and its y-covariant derivatives, there is a trivialization for the bundles E´ and E over the ball where $|y| \leq 3$ so that the connection, A, appears as an i-valued 1-form, $\nu$, that vanishes at $y = 0$, obeys $y^j \nu_j = 0$, and has uniform $C^q$ bounds. Granted this, and the fact that $|\beta| \leq c_* r^{-1}$, it is a relatively straightforward task using the Green's function $\frac{1}{4\pi |y-(\cdot)|}$ to bound the q'th order derivatives of $\beta$ at $y = 0$ by $r^{-1/2}c_q$ with $c_q$ depending only on the metric and the $C^{q+2}$ norm of u. The latter bounds imply the assertion in the second bullet of Lemma 2.3.

**c) Proof of Lemma 2.4**

The argument starts by recapitulating the derivation of (5.28); thus fix a basis of generators of $H_1(M; \mathbb{Z})$. With the connection $A_E$ given, write $A = A_E + \hat{a}_A$. As in the derivation of (5.28), fix a smooth map, $u_1: M \to S^1$ with the property that integral of the real valued 1-form $i(\hat{a}_A - u_1^{-1}du_1)$ around each of the chosen basis elements for $H_1(M; \mathbb{Z})$ lies in the interval [0, 2). As a consequence, the $L^2$-orthogonal projection of $\hat{a}_A - u_1^{-1}du_1$ onto the space of harmonic 1-forms has norm bound that depends only the metric. Now use Hodge theory to find a unique, smooth and homotopically trivial map, $u_2: S^1 \to M$ such that $\hat{a} = \hat{a}_A - u_1^{-1}du_1 - u_2^{-1}du_2$ is coclosed. The $L^2$-orthogonal projection of $\hat{a}$ to the space of harmonic 1-forms is the same as that of $\hat{a}_A - u_1^{-1}du_1$.

With $\hat{a}$ understood, standard properties of the Green's function for the operator $*d$ acting on co-closed 1-forms can be invoked to see that at any given $x \in M$, one has

$$|\hat{a}|(x) \leq c_* \left( \int_M \frac{1}{\text{dist}(x,\cdot)^2} |B_A| + 1 \right) ,\tag{6.9}$$

where c depends only on the Riemannian metric. With (6.9) understood, fix $\rho > 0$ and break the integral that appears in (6.9) into the part where the distance to x is greater than ρ, and that where the distance is less than or equal to ρ. The integral over the former is no greater than

$$c_* \rho^{-2} r \int_M | ia - \psi^\dagger \tau \psi | + c_1 ,$$

(6.10)

where $c_1$ is determined solely by the Riemannian metric and the given upper bound for the $C^0$ norm of $d\mu$ and the $L^2$ norm of $\varpi$. According to Lemma 2.2, the expression in (6.10) is bounded by

$$c_* \rho^{-2} (r \int_M |1 - |\alpha|^2| + 1) + c_1 ,$$

(6.11)

where this incarnation of $c_*$ differs from that in (6.10), but it none-the-less depends only on the metric. Likewise, this is a new incarnation of $c_1$, but its value is determined by the metric and the given upper bound for the $C^1$ norm of $d\mu$ and the $L^2$ norm of $\varpi$. Finally, Lemma 2.2 implies that the expression in (6.11) is no greater than

$$c_* \rho^{-2} E(A) + c_1 ,$$

(6.12)

where $c_*$ and $c_1$ are different then their namesakes in (6.11), by have the stated dependencies on the metric, μ and $\varpi$.

Now consider the contribution to the integral in (6.9) from the portion of the integration domain where the distance to x is no greater than ρ. As can be seen from (2.6) and Lemma 2.2, this part is bounded by

$$c_* r \rho + c_1$$

(6.13)

where $c_*$ depends only on the metric and $c_1$ on the metric and the given upper bound for the $C^1$ norm of μ and the $L^2$ norm of $\varpi$. Given (6.12) and (6.13), the claim made in Lemma 2.4 follows by taking $\rho = r^{-1/3} E^{1/3}$.

**d) Proof of Theorem 2.1**

The proof of Theorem 2.1 is broken into seven steps. The first five provide some preliminary results that are then used in the final steps to establish the desired conclusion.

*Step 1*: Fix a point $p \in M$. This step introduces the notion of an *adapted coordinate chart map* centered at p. Fix $\delta > 0$ and set $I = [-\delta, \delta] \subset \mathbb{R}$. Let $C \subset \mathbb{C}$ denote the disk of radius δ. An adapted, coordinate chart map centered at p is a smooth embedding, $\varphi: C \times I \to M$ that sends the origin, 0, to p and has certain additional properties. To state them, introduce z for the coordinate on I and write the complex coordinate on C as x + iy with x and y real. Use φ to identify C × I with φ(C × I). Then

- da = B dx ∧ dy *and the Reeb vector field is* $\partial_z$. *Here,* B *is a positive, z-independent function with value 1 at the origin.*
- *The metric pulls back as* $dx^2 + dy^2 + dz^2 + \mathfrak{h}$ *where* $\mathfrak{h}$ *obeys*
  a) $\mathfrak{h}(\partial_z, \partial_z) = 0$,
  b) *The restriction of* $\mathfrak{h}|_{z=0}$ *to the span of* $\{\partial_x, \partial_y\}$ *is* $m(dx^2 + dy^2)$ *where* m *vanishes at the origin, and has absolute value bounded by* $c_\varphi (x^2 + y^2)$.
  c) *The $C^4$ norm of* $\mathfrak{h}$ *is bounded by* ε.

(6.14)

Here, $c_\varphi$ is a constant. Such a coordinate chart map is constructed as follows: Use the exponential map at p to embed a disk centered at p whose tangent plane at p spans the kernel of the 1-form a. Fix coordinates (x, y) on this disk so that the metric is conformal to the Euclidean metric and differs from the latter by $\mathcal{O}(x^2 + y^2)$. If the radius of this disk is sufficiently small, then it will be everywhere transversal to the Reeb vector field. This understood, then there is a unique extension of these coordinates to coordinates (x, y, z) where $\partial_z$ is the Reeb vector field. These coordinates satisfy the conditions in (6.14).

By taking δ small, there is a map centered at each point in M. In particular, there exists $c_* > 0$ and $\delta > 0$ with the following significance: For each p ∈ M, there is map that obeys (6.14) with constant $c_\varphi < c_*$. Such a map is deemed to be an *adapted coordinate chart map* centered at p.

*Step 2*: This step introduces the vortex equations on ℂ. The latter consist of equations for a pair (𝒜, σ) where 𝒜 is an i-valued 1-form on ℂ and where σ is a complex valued function on ℂ. These equations read:

$$*d\mathcal{A} = -i(1 - |\sigma|^2) \quad and \quad \bar{\partial}_\mathcal{A} \sigma = 0 \quad and \quad |\sigma| \le 1.$$

(6.15)

Here, ∗ denotes the Euclidean Hodge star operator on ℂ = ℝ² and where $\bar{\partial}_\mathcal{A}$ is the d-bar operator for the trivial bundle ℂ × ℂ → ℂ that is defined using 𝒜 as the connection 1-form. Note that these equations are gauge invariant in the following sense: If (𝒜, σ) is a solution and u: ℂ → $S^1$ is a smooth map, then so is (𝒜 - $u^{-1}du$, uσ). Two configurations that differ in this way are said to be *gauge equivalent*. Unless stated explicitly, the discussion that follows won't distinguish between gauge equivalent solutions. Here are some basic facts about solutions to these equations (see Section 4a in [T1], or Section 2b in [T5]. See also [JT])

- *If $|\sigma| = 1$ at any point, then $|\sigma|$ is identically 1 and it is $\mathcal{A}$-covariantly constant. In this case, $(\mathcal{A}, \sigma)$ is gauge equivalent to $(0, 1)$.*
- *There exists a constant $c_1$ such that then $|\nabla_\mathcal{A} \sigma| \leq c_1$. Moreover, for each positive integer q, there exists a constant $c_q$ such that $|(\nabla_\mathcal{A})^q \sigma| \leq c_q$. Note that these constants do not depend on the particular solution $(\mathcal{A}, \sigma)$.*
- *The function $|\sigma|$ has no non-zero, local minima.*
- *The zeros of $\sigma$ are isolated, and each zero has positive local degree.*
- *If $\int_\mathbb{C} (1 - |\sigma|^2)$ is finite, then this integral is equal to $2\pi k$ with k a non-negative integer. In this case, $\sigma$ has precisely k zeros counting multiplicity.*
- *There is a constant $c \in (0, 1)$ with the following significance: Let $d: \mathbb{C} \to [0, \infty)$ denote the function that gives the distance to the set where $|\sigma| \leq \frac{1}{2}$. If $d > c^{-1}$, then*
   a) $1 - |\sigma| \leq e^{-cd}$.
   b) $|\nabla_\mathcal{A} \sigma| \leq c^{-1} e^{-cd}$.
   *Moreover, this constant c does not depend on the particular solution $(\mathcal{A}, \sigma)$.*

(6.16)

A solution $(\mathcal{A}, \sigma)$ to the vortex equations will be viewed at times as having domain of definition $\mathbb{C} \times \mathbb{R}$. In this case, there is no dependence on the $\mathbb{R}$-factor.

*Step 3*: This step explains the relevance of the vortex equation to the version of (2.4) under consideration. To this end, fix $c > 0$, fix $r \geq 1$ and fix an adapted coordinate chart map, $\varphi: C \times I \to M$. Let $I_r = [-r^{1/2}\delta, r^{1/2}\delta]$ and let $C_r \subset \mathbb{C}$ denote the disk of radius $r^{1/2}\delta$. Define $\varphi_r: C_r \times I_r \to M$ so that $\varphi_r(x, y, z) = \varphi(r^{-1/2}x, r^{-1/2}y, r^{-1/2}z)$. Now, suppose that $(A, \psi = (\alpha, \beta)) \in \text{Conn}(E) \times C^\infty(M; \mathbb{S})$. Pull back $(A, \psi)$ by $\varphi_r$ and write this pull-back as $(A_{\varphi,r}, (\alpha_{\varphi,r}, \beta_{\varphi,r}))$.

**Lemma 6.1**: *Fix $c > 0$, $R \geq 1$ and $\varepsilon > 0$; and there exists $r_*$ such that the following is true: Suppose that $r \geq r_*$ and that $(A, \psi = (\alpha, \beta))$ is such that*

$$B_A = r(\psi^\dagger \tau \psi - i a) + i\mathfrak{u} \quad \text{and} \quad D_A \psi = 0$$

(6.17)

*where $\mathfrak{u}$ is a co-closed 1-form on M with $C^3$-norm less than c. Suppose that $\varphi: C \times I \to M$ is an adapted, coordinate chart map. There exists a trivialization, $\mathfrak{u}_{\varphi,r}$, of $\varphi_r^*E$, and a solution $(\mathcal{A}_\varphi, \sigma_\varphi)$ to the vortex equations, here viewed on $\mathbb{C} \times \mathbb{R}$, such that when written with respect to this trivialization,*
- $|\alpha_{\varphi,r} - \sigma_\varphi| < \varepsilon$
- $|\nabla_{A_{\varphi,r}} \alpha_{\varphi,r} - \nabla_{\mathcal{A}_\varphi} \sigma_\varphi| < \varepsilon$

*at all points $(x+iy, z) \in \mathbb{C} \times \mathbb{R}$ with $x^2 + y^2 + z^2 \leq R^2$.*

*Proof of Lemma 6.1*: There are trivializations of $\varphi_r^*E$ and $\varphi_r^*K^{-1}$ such that the triple $(A_{\varphi,r}, (\alpha_{\varphi,r}, \beta_{\varphi,r}))$ on its domain of definition in $\mathbb{R}^3 = \mathbb{C} \times \mathbb{R}$ obeys

- $|\alpha_{\varphi,r}| < 1 + r^{-1}\kappa$ *and* $|\beta_{\varphi,r}| \le r^{-1/2}\kappa$.
- $B_{A_{\varphi,r}} = -i(1 - |\alpha_{\varphi,r}|^2)\,dz + r^{-1}q_0$,
- $|\nabla_{A_{\varphi,r}} \alpha_{\varphi,r}| \le \kappa$.
- $(\nabla_{A_{\varphi,r}})_z \alpha_{\varphi,r} = r^{-1} q_+$,
- $\bar{\partial}_{A_{\varphi,r}} \alpha_{\varphi,r} = r^{-1} q_-$,

(6.18)

where $\kappa$ and the three versions of q are bounded independent of r, $\varphi$, as are their derivatives. Indeed, these bounds follow from Lemmas 2.2 and 2.3 by rescaling.

Suppose that no such $r_*$ exists for some given $\varepsilon$ and R. One could then find sequences consisting of adapted coordinate chart maps, values of r tending to infinity and corresponding solutions to (6.17) where the conclusions of the lemma fail on each element in the sequence for $\varepsilon$ and R. Even so, by virtue of (6.18), the resulting sequence of triples $(A_{\varphi,r}, (\alpha_{\varphi,r}, \beta_{\varphi,r}))$ has a subsequence that converges on compact domains in $\mathbb{C} \times \mathbb{R}$ to some $(\mathcal{A}, (\sigma, 0))$ where the pair $(\mathcal{A}, \sigma)$ solves the vortex equation. But such convergence could happen only if the conclusions of the lemma held for each member of this subsequence for the given $\varepsilon$ and R

*Step 4*: This step starts out with:

**Lemma 6.2**: *Fix $c \ge 0$ and there exists $\kappa$ with the following significance: Fix a co-closed 1-form $\mathfrak{u}$ on M with $C^3$-norm bounded by c. With $r \ge 0$ fixed, let $(A, \psi = (\alpha, \beta))$ denote a solution to (6.17). Fix an adapted coordinate chart map $\varphi: \mathbb{C} \times I \to M$. Then $|\nabla_{Az}\alpha| \le \kappa$.*

*Proof of Lemma 6.2*: The Dirac equation sets $\nabla_{Az}\alpha$ equal to linear combinations of the $\nabla'$ covariant derivatives of $\beta$ and products of $\beta$ with metric dependent terms. This understood, then the assertion follows from Lemmas 2.2 and 2.3.

With Lemma 6.2 in hand, fix a smooth function, $\chi: [0, \infty) \to [0, 1]$ with compact support that equals 1 on $[0, \tfrac{1}{4})$ and 0 on $[\tfrac{1}{2}, \infty)$. With $r \ge 1$ given, set $\chi_r: \mathbb{C} \to [0, 1]$ to denote the function $\chi(r^{1/2}|x+iy|)$. Fix $c \ge 0$, and let $\mathfrak{u}$ denote a co-closed 1-form on M whose $C^3$-norm is bounded by c. Fix $r \ge \delta^{-1}$ and suppose that $(A, (\alpha, \beta))$ is a solution to (6.17). Let $\varphi: \mathbb{C} \times I \to M$ denote an adapted, coordinate chart map that sends $\varphi(0)$ to a point where $|\alpha| \le \tfrac{3}{4}$. Introduce the function on the interval I that sends z to

$$L(z) = r \int_{C \times \{z\}} \chi_r (1 - |\alpha|^2)^2$$

(6.19)

It follows from Lemma 6.1 and (6.16) that there exists $\kappa_* \in (0, 1)$ which is independent of r and c, and there exits $r_* > 0$ which depends only on c such that if $r \geq r_*$, then

$$L(0) \geq \kappa_*.$$

(6.20)

Meanwhile, Lemma 6.2 implies that $r_*$ and $\kappa_*$ can be chosen so that

$$|\partial_z L| \leq \kappa_*^{-1} (E(A) + 1).$$

(6.21)

It follows from (6.21) that

$$L(z) \geq \tfrac{1}{2} \kappa_* \quad \text{provided that} \quad |z| \leq \tfrac{1}{2} \kappa_*^2 (E(A) + 1)^{-1}.$$

(6.22)

To present a key consequence of this last assertion, fix $E \geq E(A)$. Introduce

$$R_E = \min(\tfrac{1}{2} \kappa^2 (E + 1)^{-1}, \tfrac{1}{2} \delta, \tfrac{1}{64} c_*^{-1})$$

(6.23)

where $c_*$ is a chosen constant that dominates (6.14)'s constant $c_\varphi$ if $\varphi$ is an adapted, coordinate chart map. Let $\Delta \subset C$ denote the disk with center at 0 and radius $r^{1/2}$. Note that if $R \leq R_E$, then $\Delta \times [-R, R] \subset C \times I$. Moreover, the Riemannian metric on $\Delta \times [-R, R]$ from its embedding via $\varphi$ in M differs from the product metric that comes by embedding $\Delta$ via $\varphi$ as $\Delta \times \{0\} \subset C \times I$ by no more than $c'\delta$ with $c'$ a $\delta$ and r independent constant.

Now suppose that $R \leq R_E$ and that $\varphi: C \times I \to M$ is an adapted coordinate chart map with $|\alpha| \leq \tfrac{3}{4}$ at $\varphi(0)$. Then (6.21) implies

$$\tfrac{1}{4} \kappa_* R \leq r \int_{\Delta \times [-R, R]} |1 - |\alpha|^2| \ .$$

(6.24)

*Step 5*: This step establishes various consequences of (6.24). Here is the first:

**Lemma 6.3**: *Given $c \geq 0$ and $E \geq 0$, there is a constant $\kappa \geq 1$ with the following significance: Fix a co-closed 1-form $u$ on M whose $C^3$ norm is bounded by c. With $r \geq \kappa$ fixed, let $(A, \psi = (\alpha, \beta))$ denote a solution to (6.17) such that $E(A) \leq E$. Fix an adapted coordinate chart map $\varphi: C \times I \to M$. Then there are no more than $\kappa$ disjoint disks of radius $2r^{1/2}$ in $C \times \{0\}$ with distance $\tfrac{1}{2} \delta$ or less from the origin and such that $|\alpha| \leq \tfrac{3}{4}$ at the center point.*

*Proof of Lemma 6.3*: This follows from (6.24) by taking $R = R_E$. Indeed, Lemma 2.2 can be used to find a constant $c_1$ that depends only on the constant $c$ and is such that

$$r \int_M |1 - |\alpha|^2| \leq E(A) + c_1 .$$

(6.25)

Meanwhile, by virtue of (6.24), the integral on the left hand side is no less than $\frac{1}{4} \kappa_* R_E$ times the number of disks that obey the lemma's stated conditions.

Lemma 6.3 has a the following corollary:

**Lemma 6.4**: *Given $c \geq 0$, $E \geq 0$ and $\varepsilon > 0$, there is a constant $\kappa \geq 1$ with the following significance: Fix a co-closed 1-form $\mathfrak{u}$ on M whose $C^3$ norm is bounded by $c$. With $r \geq \kappa$ fixed, let $(A, \psi = (\alpha, \beta))$ denote a solution to (6.17) such that $E(A) \leq E$. Fix an adapted coordinate chart map $\varphi \colon C \times I \to M$. Let $\Lambda \subset C \times \{0\}$ denote the set of points where $\alpha$ vanishes. There are at most $\kappa$ points in $\Lambda$. Moreover, if $p \in C \times \{0\}$, if $|p| \leq \frac{1}{2} \delta$, and if $\mathrm{dist}(p, \Lambda) > \kappa \, r^{-1/2}$, then $|\alpha|_p > 1 - \varepsilon - \kappa r^{-1}$.*

*Proof of Lemma 6.4*: Given Lemma 6.3, this follows from the third, fourth and final bullets in (6.16).

The next lemma provides a result that is closely related to that given in the previous lemma.

**Lemma 6.5**: *Fix $c \geq 0$ and $E \geq 0$. There is a constant, $\kappa > 1$ with the following significance: Fix a co-closed 1-form $\mathfrak{u}$ on M whose $C^3$ norm is bounded by $c$. With $r \geq \kappa$ fixed, let $(A, \psi = (\alpha, \beta))$ denote a solution to (6.17) such that $E(A) \leq E$. Let $\varphi \colon C \times I \to M$ denote an adapted, coordinate chart map whose center is a zero of $\alpha$. Let $z \in [-\kappa^{-1}, \kappa^{-1}]$. Then $\alpha^{-1}(0)$ intersects $C \times \{z\}$ at a point with distance less than $\kappa r^{-1/2}$ in $C \times \{z\}$ from the point $(0, z)$.*

*Proof of Lemma 6.5*: This follows from (6.22) using the third, fourth and sixth bullets in (6.16).

Here is a final consequence of (6.24):

**Lemma 6.6:** *Given $c \geq 0$ and $E \geq 0$, there is a constant $\kappa \geq 1$ with the following significance: Fix a co-closed 1-form $\mathfrak{u}$ on M whose $C^3$ norm is bounded by $c$. With $r \geq \kappa$ fixed, let $(A, \psi = (\alpha, \beta))$ denote a solution to (6.17) such that $E(A) \leq E$. Suppose that $R \in (4r^{-1/2}, R_E)$. Let $U \subset M$ denote a ball of radius $R$ with center where $\alpha = 0$. Then* $r \int_U |1 - |\alpha|^2| \geq \kappa^{-1} R$. *As a consequence, any set of disjoint balls of radius R whose centers lie where $\alpha$ is zero has at most $\kappa R^{-1}(E + 1)$ elements.*

*Proof of Lemma 6.6*: The first claim just restates (6.24). To prove the second, let $\Theta$ denote a maximal set of disjoint, radius R balls whose centers lie where $\alpha$ is zero. Let n denote the number of balls in this set. It follows from (6.24) and Lemma 2.2 that each ball in $\Theta$ contributes at least $\frac{1}{4} \kappa_* R - c_* R^3$ to the integral that computes $E(A)$. Here, $c_*$ depends only on the $C^3$ norm $\mu$ and the $L^2$ norm of $\varpi$. Thus, the union of the balls from $\Theta$ contributes at least $\frac{1}{4} n \kappa_* R - c_* \text{vol}(M)$ to the integral for $E(A)$. As a consequence, n can be no greater than $4\kappa_*^{-1} R^{-1}(E + c_* \text{vol}(M))$. This last bound gives the second assertion of the lemma.

*Step 6*: Now consider a sequence $(A_n, \psi_n = (\alpha_n, \beta_n))$ as given in the statement of Theorem 2.1. Fix E so that $E(A_n) \leq E$ for all n, and fix c so as to be greater than the $C^3$ norms for $*d\mu_n + \varpi_n$. No generality is lost by assuming that each $r_n$ is large enough so that the conclusions of Lemmas 6.1-6.6 hold.

Let $R_E$ be as in (6.23). For each $m \in \{1, 2, \ldots\}$, let $\rho_m = (\frac{1}{32})^m R_E$. For each m and for all n such that $r_n^{-1} < \frac{1}{64} \rho_m$, choose a maximal set of disjoint balls of radius $\rho_m$ with centers on $\alpha_n^{-1}(0)$. This set is non-empty for n large due to the fact that $\sup_M |1 - |\psi_n||$ is assumed to have an n-independent, positive lower bound. Indeed, granted this bound, it follows from the second, third and final points in (6.16) that there are points in M where $\alpha_n$ is zero. Denote this maximal set of balls by $\Theta_{n,m}$. For each ball $U \in \Theta_{n,m}$, let $p_U$ denote its center, and let $U'$ denote the ball whose radius is $4\rho_m$ and whose center is $p_U$. Note that the collection $\Theta_{n,m}' = \{U': U \in \Theta_{n,m}\}$ has the property that its members cover $\alpha_n^{-1}(0)$. Moreover, each ball from $\Theta_{n,m+1}'$ is contained in some ball from $\Theta_{n,m}'$.

To continue, for each m, let $Q_m$ denote the upper bound given by Lemma 6.6 for the case $R = \rho_m$. Thus, each $\Theta_{n,m}$ has at most $Q_m$ elements. Label the points in the set $\{p_U : U \in \Theta_{n,m}\}$ by consecutive integers starting from 1, and then add as many extra copies of the first point as needed so as to define a point, $\theta_{m,n} \in \times_{Q_m} M$.

Choose a diagonal subsequence of $\{(A_n, \psi_n)\}$ so that for each m, the corresponding subsequence $\{\theta_{m,n}\}$ converges in $\times_{Q_m} M$. For each such m, let $\theta_m$ denote the limit. Let $\Theta_m$ denote the set of radius $4\rho_m$ balls in M whose centers give the entries of $\theta_m$. Then each ball in $\Theta_{m+1}$ is contained in a ball from $\Theta_m$. This understood, use $Z_m$ to denote the union of the balls that comprise $\Theta_m$. As $Z_{m+1} \subset Z_m$, it makes sense to define

$$Z = \cap_{m=1,2,\ldots} Z_m .$$
(6.26)

As is argued in the next step, Z is the desired union of closed integral curves of the Reeb vector field.

*Step 7*:  The story on Z starts with

**Lemma 6.7**: *The set Z is a non-empty union of closed integral curves of the Reeb vector field* v.

*Proof of Lemma 6.7*:  The fact that Z is non-empty follows by compactness.  Fix an adapted coordinate chart map φ: C × I → M that sends the origin to a point in Z.  It follows from Lemma 6.4 that the intersection of Z with C × {0} consists of at most κ points, where an upper bound for κ is determined by the constants c and E.  It then follows from Lemma 6.5 that the intersection of Z with a neighborhood in C × I is a union of at most κ properly embedded, integral curves of the vector field v.  This bound on the number implies that Z is a union of a finite set of closed integral curves of v.

To complete the argument for Theorem 2.1, it is necessary to explain how to assign non-zero integer weights to the closed integral curves that comprise Z so that the resulting formal, weighted sum of loops in M gives the Poincaré dual in $H_1(M; \mathbb{Z})$ to the first Chern class of the line bundle E.  To this end, note that if α is a section of E with transversal zero locus, then $\alpha^{-1}(0)$ is Poincaré dual to the first Chern class of E.

To make use of this last observation, suppose that γ ⊂ Z is a component.  Select an adapted coordinate chart map, φ, that sends the origin to a point on γ.  Let C´ ⊂ C denote a closed subdisk centered at the origin such that C´ × {0} intersects Z only at the origin.  Let $\{(A_n, \psi_n = (\alpha_n, \beta_n))\}$ denote the diagonal subsequence that was chosen in the previous step to define Z.  Fix a trivialization of E over C´ × {0} so as to view $\alpha_n$ as a map from C´ × {0} to $\mathbb{C}$.  It then follows from the third point in (6.16) that for each n sufficiently large, $\alpha_n$ has positive winding number around ∂C´.  Note that this winding number does not depend on the chosen trivialization.  Let $k_{\gamma,n}$ denote this winding number.  The fifth point in (6.16) provides an index n-independent upper bound to $k_{\gamma,n}$.

Choose a subsequence of $\{(A_n, \psi_n)\}$ so the corresponding subsequence of $\{k_{\gamma,n}\}$ converges for each component γ of Z.  For each such component, let $k_\gamma$ denote the limit.  This is a positive integer, and is the weight that is assigned to the component γ.  With this assignment understood, it follows that $\sum k_\gamma [\gamma]$ is the Poincaré dual to the first Chern class of E.

## 7. Perturbations

The purpose of this final section is to tie up the loose ends from Section 3 by proving Lemma 3.1, Proposition 3.11 and 3.12, and by justifying the assumptions that are made in Properties 1-5 in Section 3d.

### a) Proof of Lemma 3.1

It proves useful to fix a fiducial connection, $A_E$, on E so as to identify Conn(E) with $C^\infty(M; iT^*M)$. Take the connection chosen just prior to (2.8). Let $\mathcal{H}_3$ and $\mathcal{H}_2$ denote the respective Hilbert spaces of Sobolev class $L^2_3$ and $L^2_2$ sections of $iT^*M \times \mathbb{S}$.

Given $r \geq 0$, introduce the 'universal' moduli space, $\mathcal{N}$; this the space of triples $((A, \psi), \mu)$ where $\mu \in \Omega$ and $(A = A_E + b, \psi)$ have the following properties: First, $(b, \psi) \in \mathcal{H}_3$. Second, $(A, \psi)$ solves the r and $\mu$ version of (2.5). Let $\mathcal{H}_{3\text{irr}} \subset \mathcal{H}_3$ denote the subset of pairs $(A, \psi)$ with $\psi$ not identically zero. Likewise, let $\mathcal{N}_{\text{irr}} \subset \mathcal{N}$ denote the subset of $((A, \psi), \mu)$ where $\psi$ is not identically zero. The set $\mathcal{N}_{\text{irr}}$ is the zero set of a certain section of a smooth vector bundle, $\mathcal{V} \to \mathcal{H}_{3\text{irr}} \times \Omega$. In this regard, the fiber of $\mathcal{V}$ over any given point $((A, \psi), \mu)$ is the subspace in $\mathcal{H}_2$ of pairs $(q, \varsigma)$ that obey the equation

$$*d*q - 2^{-1/2} r^{1/2} (\varsigma^\dagger \psi - \psi^\dagger \varsigma) = 0.$$

(7.1)

The section of $\mathcal{V}$ that defines $\mathcal{N}_{\text{irr}}$ sends a given element $((A, \psi), \mu)$ to the section whose $iT^*M$ and $\mathbb{S}$ components are

- $B_A - r(\psi^\dagger \tau^k \psi - ia) + i *d\mu$,
- $2r^{1/2} D_A \psi$.

(7.2)

The section of $\mathcal{V}$ just defined is denoted in what follows by $\mathfrak{s}$. This is a smooth section of $\mathcal{V}$.

Because $\mathfrak{s}$ is a smooth, the subspace $\mathcal{N}_{\text{irr}} \subset \mathcal{H}_{3\text{irr}} \times \Omega$ has the structure of a smooth Hilbert manifold near any $((A, \psi), \mu) \in \mathcal{N}_{\text{irr}}$ where the differential of $\mathfrak{s}$ is surjective. As is explained next, the differential is surjective on the whole of $\mathcal{N}_{\text{irr}}$. To this end, note first that the restriction of the differential of $\mathfrak{s}$ to the tangent vectors of the form $((b, \eta), 0)$ has respective $iT^*M$ and $\mathbb{S}$ components that are, up to a factor of -1, the $\phi = 0$ and $\mathfrak{t} = \mathfrak{s} = 0$ versions of the top two equations in (3.1). This implies that the cokernel of the restriction of $d\mathfrak{s}$ to $\mathcal{H}_3$ is finite dimensional. Let $(q, \varsigma)$ denote an element of this cokernel. This pair obeys the coupled equations

- $*dq - 2^{-1/2} r^{1/2} (\psi^\dagger \tau \varsigma + \varsigma^\dagger \tau \psi) = 0$,
- $D_A \varsigma + 2^{1/2} r^{1/2} \text{cl}(q) \psi = 0$.
- $*d*q - 2^{-1/2} r^{1/2} (\varsigma^\dagger \psi - \psi^\dagger \varsigma) = 0.$

(7.3)

If $(q, \varsigma)$ is not in the image of the differential of $\mathfrak{s}$ as applied to vectors of the form $(0, \mu)$ with $\mu \in \Omega$, then h must be $L^2$-orthogonal to all co-exact 1-forms on M. Indeed, this follows from the fact that $\Omega$ is dense in $C^\infty(M; iT^*M)$. Thus, $q = df + v$ with f a smooth, i-valued function on M and $v$ an i-valued harmonic 1-form.

Granted this form for h, then the middle equation in (7.3) finds $\varsigma = -r^{1/2}f\psi + \lambda$ where $\lambda$ obeys the equation $D_A\lambda = -r^{1/2}cl(\nu)\psi$. Meanwhile, the top equation in (7.3) asserts that $r^{1/2}(\psi^\dagger\tau\lambda + \lambda^\dagger\tau\psi) = 0$. This last equation requires that $\lambda = m\psi$ with m an i-valued function that is defined where $\psi \neq 0$. This and the fact that $D_A\lambda = -r^{1/2}cl(\nu)\psi$ requires that $\nu = dm$ where $\psi \neq 0$. This then implies that $\nu = 0$ and m is constant. Here is why: The unique continuation principle requires that $\psi$ can neither vanish on an open set, nor vanish so that its zero locus disconnects some ball in M. As a consequence, any loop in M can be homotoped a small amount so as to lie where $\psi \neq 0$. This implies that $\nu$ has zero pairing with $H_1(M; \mathbb{R})$, and so $\nu = 0$ and also m is constant. But with m constant and $\nu = 0$, then the third bullet in (7.3) demands that $(h, \varsigma) = 0$.

Let $\mathcal{N}_*$ denote the quotient of $\mathcal{N}_{irr}$ by the action of $C^\infty(M; S^1)$. This is a smooth Banach manifold. Moreover, the projection $\pi: \mathcal{N}_* \to \Omega$ is a Fredholm map of index zero. This understood, the Smale-Sard theorem [Sm] finds a residual subset of points in $\Omega$ with small norm that are regular values for $\pi$. Any $\mu$ from this residual set has only non-degenerate solutions to its version of (2.5) that are irreducible.

To continue with the proof, recall that the space of reducible solutions to (2.5) consists of pairs $(A, \psi = 0)$ where $A = A_* - \frac{1}{2}ir a + \mu$ where $A_*$ is a flat connection on E. The corresponding version of $\mathfrak{L}$ for a triple $(b, \eta, \phi)$ has components given by (5.8). Thus, if M has positive first Betti number, there are no non-degenerate irreducible solutions to (2.5). If the first Betti number is zero, then an argument much like the one just given proves that there is residual set in $\Omega$ whose version of $D_A$ has trivial kernel. This set is open and dense because a trivial kernel is preserved by small deformations. Granted this, it follows in the case where M has zero first Betti number that there is an open dense set of $\mu$ from $\Omega$ for which all solutions to (2.5) are non-degenerate.

In general, the assertion of the non-degeneracy of all solutions to (2.4) as defined by r and $\mathfrak{g}$ from an open dense set in $\mathcal{P}$ is a consequence of Theorem 12.1.2 and Lemmas 12.5.2 and 12.6.1 in [KM]. The subset is open since, as noted previously, the non-degeneracy condition is stable under perturbations.

### b) Proof of Proposition 3.11

The proof is given in five steps.

*Step 1*: This step finds a residual set of $\mu$ and a locally finite set $\{\rho_j\}$ such that Item 1 of the proposition holds when $r \notin \{\rho_j\}$. Let $\Omega^a \subset \Omega$ denote the vector subspace of forms that are $L^2$ orthogonal to a. Let $\mathcal{W}$ denote the space of tuples $((A = A_E + b, \psi), r, \mu)$ with the following properties: First, $r \in (0, \infty)$ and $\mu \in \Omega^a$. Second $(b, \psi) \in \mathcal{H}_3$ and $\psi$ is not identically zero. Finally, the pair $(A, \psi)$ solves the r and $\mu$ version of (2.5). This space is the zero locus of a section, $\mathfrak{s}$, of the vector bundle $\mathcal{V} \to \mathcal{H}_3 \times (0, \infty) \times \Omega^a$ whose fiber at $((A, \psi), r, \mu)$ consists of the subspace of pairs $(q, \varsigma) \in \mathcal{H}_2$ that satisfy (7.1). The section $\mathfrak{s}$ is given by (7.2). The space $\mathcal{W}$ is a $C^\infty$ Banach manifold if the differential of $\mathfrak{s}$ is surjective along $\mathcal{W} = \mathfrak{s}^{-1}(0)$. The argument given in the preceding subsection shows that such is the case. Let $\mathcal{W}_*$ denote the quotient of $\mathcal{W}$ by the space of maps from M to $S^1$. This is also a smooth Banach manifold. (See, e.g. Chapter 9.3 in [MK].)

Let $\pi: \mathcal{W}_* \to \Omega^a$ denote map that is induced by the projection. The map $\pi$ is a Fredholm map, now of index 1. Its fiber over any given $\mu \in \Omega^a$ consists of the gauge equivalence classes of triples $((A, \psi), r)$ such that $(A, \psi)$ obey the r and $\mu$ version of (2.5). The Sard-Smale theorem finds a residual subset of $\Omega^a$ that consists of regular values of $\pi$. Suppose that $\mu$ is in this set, and introduce $\mathcal{W}_\mu$ to denote $\pi^{-1}(\mu) \subset \mathcal{W}_*$. This is a smooth, 1-dimensional manifold. Let $\pi_r: \mathcal{W}_\mu \to (0, \infty)$ denote the function that assigns r to the gauge equivalence classs of $((A, \psi), r)$. The map $\pi_r$ has an open, dense set of regular values. A given r is a regular value of $\pi_r$ if and only if operator $\mathfrak{L}$ as defined by $((A, \psi), r)$ has trivial kernel. This follows from the fact that $\mu$ is a regular value of $\pi$ and r is a regular value of $\pi_r$. To elaborate, note that a tangent vector to $\mathcal{W}$ at a given $((A, \psi), r, \mu)$ has the form $\mathfrak{v} = ((b, \eta), s, \upsilon)$, where $(b, \eta) \in C^\infty(M; iT^*M \oplus \mathbb{S})$, where $s \in \mathbb{R}$, and where $\upsilon \in \Omega^a$. These are such that

- $*db - r(\psi^\dagger \tau \eta + \eta^\dagger \tau \psi) - i(\upsilon - sa) = 0$
- $D_A \eta + cl(b)\psi = 0$

(7.4)

The form $\mu$ is a regular value of $\pi$ and r is a regular value of $\pi_r$ if and only if all possible choices for $\upsilon$ and s appear in (7.4). Such is the case if and only the kernel of $\mathfrak{L}$ is trivial.

It also follows from (7.4) that $\mu$ is a regular value of $\pi$ if and only if the kernel of $\mathfrak{L}$ at each point in $\mathcal{W}_\mu$ has dimension 1 or less, and at points where its dimension is one, the $L^2$-orthogonal projection of (ia, 0) to the kernel spans the kernel.

To continue with the proof of the proposition, suppose that $\mu$ is a regular value of $\pi$. Then $\pi_r$ is a function on $\mathcal{W}_\mu$ and so Sard's theorem implies that it has an open and dense set of regular values. Note in this regard that the level sets of $\pi_r$ are compact. What follows explains why the critical values of $\pi_r$ form a locally finite set. For this purpose, suppose that $\mathfrak{c} = (A, \psi)$ and $(\mathfrak{c}, r) \in \mathcal{W}_\mu$ is a critical point of $\pi$. Let $\mathfrak{b}$ span the kernel of $\mathfrak{L}$ at $\mathfrak{c}$. Note that $\mathfrak{b}$ has the form $(q_0, \varsigma_0, 0)$ where $(q = q_0, \varsigma = \varsigma_0)$ obeys (7.1). Let $\mathcal{V}_0$ denote the set of solutions to (7.1) that are $L^2$-orthogonal to $(q_0, \varsigma_0)$. As explained momentarily, now standard perturbation theory (as pioneered by Kuranishi) with the slice theorems from Chapter 9 of [KM] finds a neighborhood of zero in $\mathbb{R}$; a neighborhood, $\Delta$, of the graph in $\mathbb{R}^2$ of the function f; and real analytic real analytic map, $\Phi: \Delta \to \mathcal{V}_0$; all with the following properties:

- F *and its first derivatives vanish at the origin.*
- $\Phi$ *vanishes at the origin.*

- *For each $z = (x, y) \in \Delta$, let $\mathfrak{c}(z) = (A + xb_0, \psi + x\eta_0) + \Phi(z)$. Then the map from $\Delta$ into $\mathrm{Conn}(E) \times C^\infty(M; \mathbb{S}) \times (0, \infty)$ that sends $z \to (\mathfrak{c}(z), r + y)$ maps the graph of f diffeomorphically onto a neighborhood of $(\mathfrak{c}, r)$ in $\mathcal{W}_\mu$.*

(7.5)

Note in this regard that f and $\Phi$ are real analytic because the non-linearities in (2.4) are given locally by real analytic (quadratic) functions of the components of A and $\psi$.

What is given in (7.5) endows $\mathcal{W}_\mu$ with a real analytic structure near $(\mathfrak{c}, r)$ that identifies $\pi_r$ with a real analytic function. As such, the set of critical points of $\pi_r$ is a real analytic set. Since $\pi_r$ is a proper map, this implies that the regular values of $\pi_r$ form a locally finite set. Granted this last conclusion, then Item 1) of the proposition holds for a given regular value, $\mu$, of $\pi$ if $\{\rho_j\}$ includes the set of critical values of the function $\pi_r$ on $\mathcal{W}_\mu$.

The fact that (7.5) holds can be seen as follows: Given a point $(x, y)$ near 0 in $\mathbb{R}^2$, a solution to the $r + y$ and $\mu$ version of (2.5) near to $(A, \psi)$ in $\mathrm{Conn}(E) \times C^\infty(\mathbb{S})$ is gauge equivalent to one that has the form $(A + xb_0, \psi + x\eta_0) + \Phi$ where $\Phi \in \mathcal{V}_0$. As such, $\Phi$ obeys a non-linear equation that has the schematic form

$$\mathcal{L}\Phi = (-iya, 0) + x^2 \mathfrak{R}_0 + 2x \mathfrak{R}_1(\Phi) + \mathfrak{R}_2(\Phi, \Phi)$$

(7.6)

where both sides are to be viewed as elements in $\mathcal{V}_0$, and where $\mathcal{L}$ is obtained from the first two lines in (3.1) by setting $\mathfrak{t}$, $\mathfrak{s}$ and $\phi$ equal to 0. Meanwhile, $\mathfrak{R}_0$, $\mathfrak{R}_1$ and $\mathfrak{R}_2$ are elements in $\mathcal{V}$ that are independent of $\Phi$, linear in the components of $\Phi$ and quadratic in the components of $\Phi$. Let $\Pi_0$ denote the projection of (7.6) onto the span of $(q_0, \varsigma_0)$. If $\Phi$ has small $L^2_2$ norm, then it must be a fixed point of the mapping from a small radius ball in $\mathcal{V}_0$ to $\mathcal{V}_0$ that sends a given $\Phi'$ to

$$T(\Phi') = \mathcal{L}^{-1}(1 - \Pi_0)((-iya, 0) + x^2 \mathfrak{R}_0 + 2x \mathfrak{R}_1(\Phi') + \mathfrak{R}_2(\Phi', \Phi')).$$

(7.7)

This is a contraction mapping on a small radius ball in $\mathcal{V}_0$ if x and y are small. As a consequence, there is a unique solution for any such pair $(x, y)$. As the mapping depends in a real-analytic fashion on x, y and the components of $\Phi'$, so the fixed point will vary with x and y in a real analytic fashion. Writing $(x, y) = z$, let $z \to \Phi(z)$ denote the resulting map from a neighborhood of 0 in $\mathbb{R}^2$ to a neighborhood of 0 in $\mathcal{V}_0$.

Having solved most of (7.6), there remains yet the projection of (7.6) to the span of $(q_0, \varsigma_0)$. As this vector is in the kernel of $\mathcal{L}$, the vanishing of the projection of (7.6) onto $(q_0, \varsigma_0)$ asserts that

$$\Pi_0((-iya, 0) + x^2 \mathfrak{R}_0 + 2x \mathfrak{R}_1(\Phi(z)) + \mathfrak{R}_2(\Phi(z), \Phi(z))) = 0.$$

(7.8)

Now, as remarked previously, (ia, 0) has non-zero inner product with $(q_0, \zeta_0)$. Thus, this equation can be rewritten to read:

$$y + h(x, y) = 0 \tag{7.9}$$

where h is a real-analytic function of x and y that vanishes with its first derivatives at the origin. This being the case, the contraction mapping theorem can be used to find a function $y = f(x)$ with f a real analytic function defined near zero in $\mathbb{R}$ such that $(y, x)$ obeys (7.9) near 0 in $\mathbb{R}^2$ if and only if $y = f(x)$.

*Step 2*: This step finds a residual subset of $\mu \in \Omega^a$ which are regular values of $\pi$, and such that the assertions of Items 1) and 2) of the proposition hold. To start introduce $\mathcal{W} \otimes \mathcal{W} \subset \mathcal{W} \times \mathcal{W}$ to denote the subset of pairs of the form $((\mathfrak{c}, r, \mu), (\mathfrak{c}', r, \mu))$ such that $\mathfrak{c}$ is not gauge equivalent to $\mathfrak{c}'$ and such that both the $\mathfrak{c}$ and $\mathfrak{c}'$ versions of $\mathcal{L}$ have trivial kernel. Note in particular that $\mathfrak{c}$ and $\mathfrak{c}'$ obey the same version of (2.5).

This $\mathcal{W} \otimes \mathcal{W}$ is a smooth submanifold of $\mathcal{W} \times \mathcal{W}$. To see why, note first that the set $\mathcal{W}_0 \subset \mathcal{W}$ of elements $(\mathfrak{c}, r, \mu)$ for which the kernel of $\mathcal{L}$ is trivial is an open (dense) set. Thus, $\mathcal{W}_0 \times \mathcal{W}_0$ is open in $\mathcal{W} \times \mathcal{W}$. This understood, then $\mathcal{W}_0 \otimes \mathcal{W}_0$ is the inverse image in $\times_2 ((0, \infty) \times \Omega^a)$ of the diagonal via the projections $((\pi_r, \pi), (\pi_r, \pi))$. Hence $\mathcal{W}_0 \otimes \mathcal{W}_0$ is a manifold if this map is transversal to the diagonal. Such is the case for $((\mathfrak{c}, \cdot), (\mathfrak{c}', \cdot))$ when both the $\mathfrak{c}$ and $\mathfrak{c}'$ versions of $\mathcal{L}$ have trivial cokernel.

Define the function

$$\mathfrak{h}: \mathcal{W} \otimes \mathcal{W} \to \mathbb{R} \tag{7.10}$$

so as to send $((\mathfrak{c}, \cdot), (\mathfrak{c}', \cdot))$ to $\mathfrak{h} = \mathfrak{a}(\mathfrak{c}) - \mathfrak{a}(\mathfrak{c}')$. What follows explains why 0 is not a critical value of $\mathfrak{h}$. To this end, consider that the derivative of $\mathfrak{a}$ at $(\mathfrak{c}, r, \mu) \in \mathcal{W}$ in the direction of the tangent vector $\mathfrak{v} = ((b, \eta), s, \upsilon)$ is

$$\partial_\mathfrak{v} \mathfrak{a} = -\tfrac{1}{2} sE + \mathfrak{e}_\upsilon . \tag{7.11}$$

If the kernel of $\mathcal{L}$ is trivial, then any pair $(s, \upsilon)$ can appear in (7.11). As a consequence, the differential of $\mathfrak{h}$ at $((\mathfrak{c}, r, \mu), (\mathfrak{c}', r, \mu))$ is zero if and only if $B_A = B_{A'}$. This last condition requires that $\psi^\dagger \tau \psi = \psi'^\dagger \tau \psi'$. As a consequence, $\psi' = u\psi$ with $|u| = 1$ where $\psi \neq 0$. Meanwhile, Hodge theory finds that $A' = A - i\nu$ where $\nu$ is a closed 1-form. Because $D_A \psi = 0$ and $D_{A'} \psi' = 0$, these last two conclusions demand that $i\nu = u^{-1}du$ at points where $\psi \neq 0$. As noted previously, $\psi$ can not vanish on an open set, nor can its zero locus disconnect any ball in M. Each class in $H_1(M; \mathbb{Z})$ has a generating loop that avoids the zero locus of $\psi$. It then follows that $\nu$ has integral periods around each such generator. This means that $(A, \psi)$ and $(A', \psi')$ are gauge equivalent, which is forbidden.

By virtue of what was just said, 0 is a regular value of $\mathfrak{h}$ on $\mathcal{W}_0 \otimes \mathcal{W}_0$ and so $\mathfrak{h}^{-1}(0) \subset \mathcal{W}_0 \otimes \mathcal{W}_0$ is a smooth, codimension 1 submanifold. Let $\mathcal{W}_* \otimes \mathcal{W}_*$ denote the quotient of $\mathcal{W}_0 \otimes \mathcal{W}_0$ by the action of $C^\infty(M; S^1) \times C^\infty(M; S^1)$. This is a smooth Banach manifold, and the projection from this manifold to $\Omega^a$ is Fredholm with index 0. As such, it has a residual set of regular values. If $\mu$ is a regular value, then the fiber in $\mathfrak{h}^{-1}(0)$ over $\mu$ is a zero dimensional manifold, thus a locally finite set of points.

Since the intersection of two residual sets is residual, there is a residual set of points in $\Omega^a$ that are simultaneously regular values for the projection on $\mathcal{W}_*$ and the projection on $\mathfrak{h}^{-1}(0) \subset \mathcal{W}_* \otimes \mathcal{W}_*$. If $\mu$ is a regular value for both projections, then the conclusions of Items 1) and 2) hold for some locally finite set $\{\rho_j\} \subset (r_k, \infty)$.

**c) The proof of Proposition 3.12**

The first point to make is one made before by Lemma 3.6: If $r \in (\rho_i, \rho_{i+1})$ and $\mathfrak{q} \in \mathcal{P}$ is chosen so that $(\mu, \mathfrak{q})$ is strongly $(r, k)$ admissible, then the pair $(\mu, \mathfrak{q})$ will be strongly $(r', k)$ admissible for all $r'$ in some some neighborhood of $r$.

The next point to make is that there is a smooth function, $\varepsilon_0 \colon (\rho_i, \rho_{i+1}) \to (0, 1)$ with limit 0 as $r \to \rho_i$ and as $r \to \rho_{i+1}$ such that if $r \in (\rho_i, \rho_{i+1})$ and if $\mathfrak{q} \in \mathcal{P}$ has norm less than $\varepsilon_0(r)$, then the following is true: First, $\mathfrak{q}$ lies in the radius $\varepsilon = \varepsilon(r)$ ball that is described in Proposition 3.5. Thus, the solutions to the $r$ and $\mathfrak{g} = \mathfrak{e}_\mu + \mathfrak{q}$ version of (2.4) with degrees $k$ or greater are non-degenerate, and their gauge equivalence classes are in 1-1 correspondence via the map $\mathfrak{c}(\cdot)$ from Proposition 3.5 with those of the $r$ and $\mu$ version of (2.5). Second, this correspondence between gauge equivalence classes is such that the ordering imposed on the gauge equivalence classes of solutions to the $r$ and $\mathfrak{g} = \mathfrak{e}_\mu + \mathfrak{q}$ version of (2.4) by this same $r$ and $\mathfrak{g}$ version of (2.9) is the same as that imposed on the set of gauge equivalence classes of solutions to the $r$ and $\mu$ version of (2.5) by the $r$ and $\mathfrak{g} = \mathfrak{e}_\mu$ version of (2.9). Indeed, if $\varepsilon_0(r)$ is small, then each solution to the $r$ and $\mathfrak{g} = \mathfrak{e}_\mu + \mathfrak{q}$ version of (2.4) will be very close to the gauge orbit of its corresponding solution to the $r$ and $\mu$ version of (2.5), in particular, much closer to the latter then it is to any other such gauge orbit.

**Lemma 7.1**: *Let $\varepsilon_1(\cdot) \colon (\rho_i, \rho_{i+1}) \to (0, 1)$ denote a continuous function with limit 0 as $r \to \rho_i$ and as $r \to \rho_{i+1}$ such that $\varepsilon_1(\cdot) < \varepsilon_0(\cdot)$ at all $r \in (\rho_i, \rho_{i+1})$. Then there is a contiguous set $\mathbb{J}(i) \in \mathbb{Z}$, an increasing sequence $\{t_n\}_{n \in \mathbb{J}(i)} \subset (\rho_i, \rho_{i+1})$, and a sequence $\{\mathfrak{q}_n\}_{n \in \mathbb{J}(i)} \subset \mathcal{P}$. These are such that the following is true for each $m \in \mathbb{J}(i)$, let $r_m = \frac{1}{2}(t_m + t_{m+1})$.*

- *$(\mu, \mathfrak{q}_m)$ is strongly $(r, k)$-admissible for all $r \in [t_m, t_{m+1}]$.*
- *$\|\mathfrak{q}_m\|_\mathcal{P} < \varepsilon_1(r)$ for all $r \in [t_m, t_{m+1}]$.*

***Proof of lemma 7.1***: Since the condition of being strongly (r, k) admissable is an open condition, the existence of this data follows from Lemma 3.6 and the fact that the open interval $(\rho_i, \rho_{i+1})$ is locally compact.

A particular version of the function $\varepsilon_1(\cdot)$ is needed when it is time to prove that the cSWF homology changes in the required manner as r crosses a given $t_m \in \{t_n\}_{n \in \mathbb{J}(i)}$.

The rest of the proof of Proposition 3.12 has two parts.

*Part 1*: This part of the proof explains how to compare the cSWF complexes and their homology as r crosses any given $t_m \in \{t_n\}_{n \in \mathbb{J}(i)}$. Suppose that $\varepsilon > 0$ has been chosen, that $\varepsilon_1(r) < \varepsilon$ for $r \in [t_{m-1}, t_m]$, and that both $\mathfrak{q}_{m-1}$ and $\mathfrak{q}_m$ lie in the radius $\varepsilon$ ball about the origin in $\mathcal{P}$. Let $s \to \mathfrak{q}(s)$ denote a path in this ball, parametrized by $s \to \mathbb{R}$ such that $\mathfrak{q}(s) = \mathfrak{q}_{m-1}$ where $s < -1$, such that $\mathfrak{q}(s) = \mathfrak{q}_m$ where $s > 1$, and such that $|\frac{d}{ds}\mathfrak{q}| < \varepsilon$ for all s. At each s, $\mathfrak{g}(r, s) = \mathfrak{e}_\mu + \mathfrak{q}(s)$ defines perturbation terms $(\mathfrak{T}_s, \mathfrak{S}_s)$ for use in (2.11). This s-dependent perturbation gives the following generalization of (2.11):

- $\frac{\partial}{\partial s} A = -B_A + r(\psi^\dagger \tau^k \psi - ia) + \mathfrak{T}_s(A, \psi)$,
- $\frac{\partial}{\partial s} \psi = -D_A \psi + \mathfrak{S}_s(A, \psi)$.

(7.12)

Of interest here are solutions to (7.13) where $\lim_{s \to -\infty} (A, \psi)$ is a solution to the version of (2.4) that is defined by the given r and $\mathfrak{g} = \mathfrak{e}_\mu + \mathfrak{q}_{m-1}$, and where $\lim_{s \to \infty} (A, \psi)$ is a solution to the version of (2.4) that is defined by the given r and $\mathfrak{g} = \mathfrak{e}_\mu + \mathfrak{q}_m$. In particular, given solutions $\mathfrak{c}_-$ and $\mathfrak{c}_+$ to the respective $r = t_m$ and $\mathfrak{g} = \mathfrak{e}_\mu + \mathfrak{q}_{m-1}$ and $\mathfrak{g} = \mathfrak{e}_\mu + \mathfrak{q}_m$ versions of (2.4), let $\mathcal{M}_{\mathfrak{q}(\cdot)}(\mathfrak{c}_-, \mathfrak{c}_+)$ denote the of solutions to (7.13) with $s \to -\infty$ limit equal to $\mathfrak{c}_-$ and with $s \to \infty$ limit equal to $u \mathfrak{c}_+$ with u a smooth map from M to $S^1$.

According to Proposition 24.4.7 in [KM], there are paths $s \to \mathfrak{q}(s)$ as just described where $\mathfrak{q}(s)$ is in the radius $\varepsilon$ ball in $\mathcal{P}$ for all s, and such that the following is true: If $\mathfrak{c}_-$ has degree $d_- \geq k$ and $\mathfrak{c}_+$ has degree $d_+ \geq k$, then $\mathcal{M}_{\mathfrak{q}(\cdot)}(\mathfrak{c}_-, \mathfrak{c}_+)$ has the structure of a smooth, manifold of dimension $d_- - d_+$. Fix such a path. Of interest in what follows is the case where $\mathfrak{c}_-$ and $\mathfrak{c}_+$ have the same dimension. In this case, it follows from Theorem 24.6.2 in [KM] that $\mathfrak{q}(\cdot)$ can be found with the added feature that $\mathcal{M}_{\mathfrak{q}(\cdot)}(\mathfrak{c}_-, \mathfrak{c}_+)$ is compact.

As explained in Chapter 25.2 of [KM], each element in the each $(\mathfrak{c}_-, \mathfrak{c}_+)$ version of $\mathcal{M}_{\mathfrak{q}(\cdot)}(\mathfrak{c}_-, \mathfrak{c}_+)$ can be given assigned a sign, either +1 or -1. For a given such pair $(\mathfrak{c}_-, \mathfrak{c}_+)$, let $\sigma(\mathfrak{c}_-, \mathfrak{c}_+) \in \mathbb{Z}$ denote the sum of these signs, with the understanding that $\sigma(\mathfrak{c}_-, \mathfrak{c}_+) = 0$ if $\mathcal{M}_{\mathfrak{g}(\cdot)}(\mathfrak{c}_-, \mathfrak{c}_+) = \emptyset$.

To explain the significance of this number, let $\mathfrak{c}_v$ denote a generator of the canonical basis in degree k or greater. Use $\mathfrak{c}_{v-}$ to denote the corresponding gauge equivalence class of solutions to the r and $\mathfrak{g} = \mathfrak{e}_\mu + \mathfrak{q}_{m-1}$ version of (2.4); and use $\mathfrak{c}_{v+}$ denote the corresponding gauge equivalence class of solutions to the r and $\mathfrak{g} = \mathfrak{e}_\mu + \mathfrak{q}_m$ version of (2.4). Let $\mathbb{V}$ denote the vector space generated over $\mathbb{Z}$ by the canonical basis elements in degrees k and greater. Define a linear map $\mathbb{T}: \mathbb{V} \to \mathbb{V}$ by the rule

$$\mathbb{T}\mathfrak{c}_v = \sum_{v'} \sigma(\mathfrak{c}_{v-}, \mathfrak{c}_{v'+}) \, \mathfrak{c}_{v'},$$

(7.13)

where the sum is restricted to the generators that have the same degree as $v$. Now, let $\delta_{m-1}$ and $\delta_m$ denote the respective differentials of the cSWF complex in degrees k and greater as defined by using r and $\mathfrak{g} = \mathfrak{e}_\mu + \mathfrak{q}_{m-1}$, and by r and $\mathfrak{g} = \mathfrak{e}_\mu + \mathfrak{q}_m$. Chapter 25.3 of [KM] proves that $\mathbb{T}$ intertwines these differentials, thus $\mathbb{T}\delta_{m-1} = \delta_m \mathbb{T}$; and that it induces an isomorphism between the respective $\delta_{m-1}$ and $\delta_m$ homology groups.

With the preceding understood, the task now is to prove that the function $\varepsilon_1(\cdot)$ can be chosen so as to guarantee that $\mathbb{T}$ is upper triangular with 1's on the diagonal. The Lemma 7.2 below implies that there exists such a function.

*Part 2*: Fix $\varepsilon \in (0, \frac{1}{2}\varepsilon_0(r))$, and suppose that both $\mathfrak{q}_-$ and $\mathfrak{q}_+$ lie in the ball of radius $\varepsilon$ about the origin in $\mathcal{P}$, and both chosen so that $(\mu, \mathfrak{q}_\pm)$ is strongly $(r = t_m, k)$ admissable. Let $s \to \mathfrak{q}(s)$ denote a smooth map from $\mathbb{R}$ into $\mathcal{P}$ with $\mathfrak{q}(s) = \mathfrak{q}_-$ for $s \leq -1$, with $\mathfrak{q}(s) = \mathfrak{q}_+$ for $s \geq 1$ and such that $\|\frac{d}{ds}\mathfrak{q}\|_\mathcal{P} < \varepsilon$ for all s. As in Part 1, use $\mathfrak{q}(s)$ to define the moduli spaces $\mathcal{M}_{\mathfrak{q}(\cdot)}(\mathfrak{c}_-, \mathfrak{c}_+)$ where $\mathfrak{c}_\pm$ are respective solutions to the $r = t_m$ and $\mathfrak{g} = \mathfrak{e}_\mu + \mathfrak{q}_-$ version of (2.4), and to the $r = t_m$ and $\mathfrak{g} = \mathfrak{e}_\mu + \mathfrak{q}_+$ version.

**Lemma 7.2**: *There exists $\varepsilon > 0$ such that if $\mathfrak{q}_-$, $\mathfrak{q}_+$ and $\mathfrak{q}(\cdot)$ lie in the ball of radius $\varepsilon$ about the origin in $\mathcal{P}$, then each version of the $\mathfrak{q}(s)$ version of $\mathcal{M}_{\mathfrak{q}(\cdot)}(\mathfrak{c}, \mathfrak{c})$ has precisely one non-degenerate element. Moreover, $\mathcal{M}_{\mathfrak{q}(\cdot)}(\mathfrak{c}_-, \mathfrak{c}_+) = \emptyset$ if $\mathfrak{a}(\mathfrak{c}_+) > \mathfrak{a}(\mathfrak{c}_-)$ where $\mathfrak{a}$ here is the $r = t_m$ and $\mathfrak{g} = \mathfrak{e}_\mu$ version of (2.9).*

*Proof of Lemma 7.2*: Define (2.11) using any given r and perturbation term $\mathfrak{g}$. Let $\mathfrak{a}$ denote the corresponding version of (2.9). The equations in (2.11) imply that $\frac{d}{ds}\mathfrak{a} = -\|\nabla \mathfrak{a}\|_2^2$. As a consequence, the equations require that $\mathfrak{a}$ decrease as s increases unless the solution, $s \to \mathfrak{c}(s)$ is constant.

Note that this last point implies that when $\mathcal{M}(\mathfrak{c}, \mathfrak{c})$ is defined by the solutions to any r and $\mathfrak{g}$ version of (2.11), then it has just one element, the constant map $s \to \mathfrak{c}$. In addition, if $\mathfrak{c}$ is an unobstructed, irreducible solution to the r and $\mathfrak{g}$ version of (2.4), then this constant instanton is an unobstructed solution to (2.11)

Keeping the preceding points in mind, suppose that no such $\varepsilon$ exists. One would then have a sequence, $\{(\varepsilon_p, \mathfrak{q}_{p-}, \mathfrak{q}_{p+}, \mathfrak{q}_p(\cdot))\}_{p=1,2,\ldots}$ as described above such that $\lim_{p\to\infty} \varepsilon_p = 0$ and such that one or more of the following occurs: Either the $\mathfrak{q}_p(\cdot)$ version of $\mathcal{M}_{\mathfrak{q}(\cdot)}(\mathfrak{c}_0, \mathfrak{c}_0)$ has two or more elements for some fixed canonical basis element $\mathfrak{c}_0$. Or, there exists a pair $\mathfrak{c}_-, \mathfrak{c}_+$ such that $\mathfrak{a}(\mathfrak{c}_+) > \mathfrak{a}(\mathfrak{c}_-)$ and such that $\mathcal{M}_{\mathfrak{q}(\cdot)}(\mathfrak{c}_-, \mathfrak{c}_+)$ Let $\mathfrak{c} = \mathfrak{c} = \mathfrak{c}_0$ in the first instance, and let $c = c_-$ and $c' = c_+$ in the second. One could then use arguments from Chapters 16 and 17 in [KM] to obtain a subsequence of elements indexed by p, each from the corresponding $q_p(\cdot)$ version of $\mathcal{M}_{g(\cdot)}(c, c')$ that converged to what Kronheimer and Mrowka call a 'broken trajectory'. This consists, in part, of a set of solutions $(\eth_1, \ldots, \eth_n)$ of solutions to the r and $g = e_\mu$ version of (2.11) such that each is an instanton, such that the $s \to -\infty$ limit of $\eth_1$ is $c$, the $s \to +\infty$ limit of $\eth_n$ is $c'$, and such that for each $j = 2, \ldots, n$, the $s \to -\infty$ limit of $\eth_j$ is the $s \to +\infty$ limit of $\eth_{j-1}$. Here $n > 1$ and at least one of the $\eth_j$ can not be $\mathbb{R}$-invariant. Note that in the case $c = c'$, the sequence can't converge to the constant instanton $s \to c$ as the latter is unobstructed as a solution to (2.11). In any case, since $a(c) - a(c') \leq 0$, the sum of the changes in $a$ as s runs from $-\infty$ to $+\infty$ for the various $\eth_j$ must equal 0. At least one of these drops must be non-trivial since at least on $\eth_j$ is not constant. Thus, at least one of these drops must be positive, and this is not possible.

### d) Properties 1-5 from Section 3d

Fix $\rho_i \in \{\rho_j\}$ and then $r_-$ and $r_+$ as described at the beginning of Section 3d. The purpose of what follows is to explain how to obtain a path $r \to q(r)$ with the five properties that are listed in Section 3d. The discussion has six parts after the stage setting that follows.

Fix a smooth function, $\chi: [0, \infty) \to [0, 1]$ that equals zero on $(0, r_-]$ and $[r_+, \infty)$, and equals one on a neighborhood of $\rho_i$. Let $\mathcal{B} \subset \mathcal{P}$ denote the ball about the origin of radius 1. Fix $\upsilon \in (0, 1)$. Given $\mathfrak{p} \in \mathcal{B}$, use $q_\mathfrak{p}(\cdot)$ to denote the map from $[r_-, r_+]$ to $\mathcal{B}$ that sends r to $\upsilon \chi(r) \mathfrak{p}$. The map $r \to q(r)$ will have the form $q_\mathfrak{p}(r)$ for a particular choice of $\mathfrak{p} \in \mathcal{B}$ and $\upsilon > 0$ very small.

To see how to choose $\mathfrak{p}$, it is convenient to introduce $\mathcal{S}_*$ to denote the space of gauge equivalence classes of tuples $(r, \mathfrak{p}, \mathfrak{c}) \in (r_-, r_+) \times \mathcal{B} \times (\text{Conn}(E) \times C^\infty(M; \mathbb{S}))$ such that $\mathfrak{c}$ obeys the r and $g = e_\mu + q_\mathfrak{p}(r)$ version of (2.4) and has degree k or greater. To keep the notation under control, a given $(r, \mathfrak{p}, \mathfrak{c}) \in (r_-, r_+) \times \mathcal{B} \times (\text{Conn}(E) \times C^\infty(M; \mathbb{S}))$ will not be distinguised in what follows from its gauge equivalence class. Let $\pi: \mathcal{S}_* \to \mathcal{B}$ denote the projection, and let $\mathcal{S}_{*\mathfrak{p}}$ denote the fiber of $\pi$ over $\mathfrak{p} \in \mathcal{B}$. Take $\mu$ so as to satisfy the conditions of Proposition 3.11. In particular, take $\mu$ from the residual subsets in $\mathcal{B}$ that are described in Parts 1 and 2 of Section 7b. By virtue of what is proved in Section 7b, the fiber $\mathcal{S}_{*0}$ over $\mathfrak{p} = 0$ is a smooth, 1-dimensional manifold, a manifold that is embedded in the quotient of $(r_-, r_+) \times (\text{Conn}(E) \times C^\infty(M; \mathbb{S}))$ by the action of $C^\infty(M; S^1)$. Granted that this is the case, there exists $\varepsilon_1 > 0$ such that when $\upsilon < \varepsilon_1$, then $\pi$ is a submersion over $\mathcal{B}$. This understood, assume that $\upsilon < \varepsilon_1$. In this case, $\mathcal{S}_*$ is fibered by $\pi$ over $\mathcal{B}$.

Under certain circumstances, it is permissable to use $\mathfrak{p} = 0$ and so take Section 3d's map $q(\cdot)$ to be the zero map. The circumstances are that one and only one of the following holds:

- *All solutions to the $r = \rho_i$ and $\mu$ version of (2.5) with degree k or greater are non-degenerate, and there is precisely one pair of distinct, gauge equivalence classes of solutions to the $r = \rho_i$ and $\mu$ version of (2.5) that are not distinquished by the values of the r and $\mathfrak{g} = \mathfrak{e}_\mu$ version of (2.9).*
- *There is precisely one gauge equivalence class of solution to the $r = \rho_i$ and $\mu$ version of (2.5) with degree k or greater that is not non-degenerate. Let $\mathfrak{c}$ denote a representative of this one class where the corresponding version of $\mathfrak{L}$ has a non-trivial kernel. The function h that appears in (7.9) is such that $\frac{\partial^2 h}{\partial x^2} \neq 0$ at the origin in $\mathbb{R}^2$. Also, the gauge equivalence classes of solutions to the $r = \rho_i$ and $\mu$ version of (2.5) are distinguished by the values of the $r = \rho_i$ and $\mathfrak{g} = \mathfrak{e}_\mu$ version of (2.9).*

(7.14)

If (7.14) is satisfied, then Property 4 of Section 3d is satisfied by taking $\mathfrak{q}(\cdot) = 0$. As is explained below, arguments much like those from Section 7c can be used to establish the Properties 2, 3 and 5. Property 1 is satisfied automatically given the choice of $\mu$. To put (7.14) into a larger context, introduce the projection $\pi_r: \mathcal{S}_{*0} \to (r_-, r_+)$. This $\pi_r$ is a function on $\mathcal{S}_{*0}$. Its critical points are the triples $(r, 0, \mathfrak{c})$ where the $\mathfrak{c}$, r and $\mu$ version of the operator $\mathfrak{L}$ has a non-trivial kernel. By construction, these critical points occur only at $r = \rho_i$. Such a critical point is non-degenerate (in the sense of Morse theory, not in the sense that $\mathfrak{c}$ is a non-degenerate solution to (2.5)) if and only if the corresponding function h from (7.9) obeys $\frac{\partial^2 h}{\partial x^2} \neq 0$.

With the preceding understood, let $\mathfrak{p} \in \mathcal{B}$ and let $\pi_r: \mathcal{S}_{*\mathfrak{p}} \to (r_-, r_+)$ again denote the projection. This is a function on $\mathcal{S}_{*\mathfrak{p}}$ and its critical points consist of the triples $(r, \mathfrak{p}, \mathfrak{c})$ where the r, $\mathfrak{g} = \mathfrak{e}_\mu + \upsilon \chi \mathfrak{p}$ and $\mathfrak{c}$ version of $\mathfrak{L}$ has a non-trivial kernel. Note that if $\upsilon$ is sufficiently small, then these occur where $\chi(r) = 1$. Such a small value for $\upsilon$ is assumed in what follows.

If all the critical points are non-degenerate in the sense of Morse theory, and if they have distinct critical values, then there are but a finite set of critical values for $\pi_r$. More over, if y is a critical value of $\pi_r$ on $\mathcal{S}_{*\mathfrak{p}}$, then the second bullet in (3.5) describes the situation at y, except that there may be more than one gauge equivalence class of solution with the same value of $\mathfrak{a}_{\mathfrak{g}(y)}$. If y is not a critical value, then all solutions to the $r = y$ and $\mathfrak{g}(y) = \mathfrak{e}_\mu + \upsilon \chi(y) \mathfrak{p}$ version of (2.4) are non-degenerate in the sense used in the previous sections.

Given the preceding, the first step towards finding $\mathfrak{p}$ so that $\mathfrak{q}_\mathfrak{p}(\cdot)$ is described by Properties 1-5 in Section 3d finds $\mathfrak{p}$ such that all critical points of $\pi_r$ on $\mathcal{S}_{*\mathfrak{p}}$ are non-degenerate in the Morse theory sense. Part 1 below contains this step. Having found such $\mathfrak{p}$, note that $\pi_r$ on $\mathcal{S}_{*\mathfrak{p}'}$ will have non-degenerate critical values if $\mathfrak{p}'$ is sufficiently close to $\mathfrak{p}$. This understood, Part 2 below finds some $\mathfrak{p}'$ near $\mathfrak{p}$ where $\pi_r$'s critical points have distinct critical values. Part 3 finds $\mathfrak{p}''$ near $\mathfrak{p}'$ where the values of the corresponding of $\mathfrak{a}$ in (2.9) distinguish all of the gauge equivalence classes of solutions to the r and $\mathfrak{g} = \mathfrak{e}_\mu + \upsilon\chi(r)\mathfrak{p}''$ at the critical points of $\pi_r$. Part 4 perturbs $\mathfrak{p}'''$ so that the both the first and second bullets in (3.5) are satisfied. The remaining parts address Properties 1, 2, 3 and 5 in Section 3d.

*Part 1*: The space $\mathcal{S}_*$ is fibered by $\pi$ over $\mathcal{B}$. As a consequence, it has a 'vertical' tangent bundle, this the kernel of the differential of $\pi$. The latter is a trivial, real line bundle over $\mathcal{S}_*$. Fix a nowhere zero section, $\mathfrak{v}$, of this bundle. Thus, $\mathfrak{v}$ restricts to each $\mathcal{S}_{*\mathfrak{q}}$ as a non-zero tangent vector to $\mathcal{S}_{*\mathfrak{q}}$.

Let $\pi_r: \mathcal{S}_* \to (r_-, r_+)$ denote the map induced by the projection to $(r_-, r_+)$, and then introduce $f': \mathcal{S}_* \to \mathbb{R}$ to denote the directional derivative of $\pi_r$ by the vector field $\mathfrak{v}$. Thus, $f' = 0$ at some $(r, \mathfrak{p}, \mathfrak{c})$ if and only if $\pi_r$ has a critical point $(r, \mathfrak{p}, \mathfrak{c})$ when viewed as a function on $\mathcal{S}_{*\mathfrak{p}}$. Let $f'': \mathcal{S}_+ \to \mathbb{R}$ denote $\mathfrak{v}(f)$. Thus, $f'' = 0$ at a zero of $f'$ if and only if the corresponding critical point of $\pi_r$ on the relevant fiber of $\pi$ is degenerate in the sense of Morse theory.

**Lemma 7.3**: *There is a neighborhood, $\mathcal{B}' \subset \mathcal{B}$, of the origin such that the zero locus of the function of $f'$ in $\pi^{-1}(\mathcal{B}')$ is a smooth codimension 1 submanifold, and the zero locus of $(f', f''): \pi^{-1}(\mathcal{B}') \to \mathbb{R}^2$ is a smooth, codimension 2 submanifold.*

To see where this lemma leads, let $\mathcal{Z}_1 \subset \pi^{-1}(\mathcal{B}')$ denote the zero locus of $f'$ and let $\mathcal{Z}_2$ denote the zero locus of $(f', f'')$. The projection $\pi: \mathcal{Z}_1 \to \mathcal{B}'$ is a Fredholm map of index 0, so there is a residual set of regular values. If $\mathfrak{p}$ is such a regular value, then $f'$ has at most a finite set of zeros on $\mathcal{S}_{*\mathfrak{p}}$. Likewise, the projection $\pi: \mathcal{Z}_2 \to \mathcal{B}'$ is a Fredholm map of index -1. Thus, it too has residual set of regular values. If $\mathfrak{p}$ is a regular value for both of these projections, then $\pi_r$ on $\mathcal{S}_{*\mathfrak{p}}$ has at most a finite set of critical points, and all such points are non-degenerate critical points in the Morse theoretic sense.

*Proof of Lemma 7.3*: It is enough to prove that the respective differentials of $f'$ and of $(f', f'')$ are surjective maps to $\mathbb{R}$ and $\mathbb{R}^2$ at all points of $\mathcal{S}_{*0}$ where the relevant map is zero. To carry out this task, remember that $f'$ is zero on $\mathcal{S}_{*0}$ only at $r = \rho_i$ and only at a gauge equivalence class of some solution $\mathfrak{c}$ to the $r = \rho_i$ and $\mu$ version of (2.5) where the corresponding operator $\mathfrak{L}$ is zero. Return now to the notation used subsequent to (7.5) in Step 1 of Section 7b. Write $\mathfrak{c} = (A, \psi)$. Then a neighborhood of $(\rho_i, 0, \mathfrak{c})$ in $\mathcal{S}_*$ is parametrized by pairs $(x, \mathfrak{p})$ where x is near zero in $\mathbb{R}$ and where $\mathfrak{p}$ is near zero in $\mathcal{P}$. This parametrization has the form $(x, \mathfrak{p}) \to (\rho_i + y(x, \mathfrak{p}), \mathfrak{p}, \mathfrak{c}(x, \mathfrak{p}))$ where the notation is as follows: First, $\mathfrak{c}(x, \mathfrak{p})$ is used here to denote $(A + xb_0, \psi + x\eta_0) + \Phi(x, \mathfrak{p})$ where $\Phi \in \mathcal{V}_0$ is a smooth function of x and $\mathfrak{p}$. The latter obeys (7.6) with the term $\upsilon(\nabla\mathfrak{p})_{\mathfrak{c}(x,\mathfrak{p})}$ added on the right hand side. Here, $(\nabla\mathfrak{p})_{(\cdot)}$ is defined so that $\langle \mathfrak{b}, \nabla\mathfrak{p}\rangle_{L^2} = (\frac{d}{dt}\mathfrak{p}(\cdot + t\mathfrak{b}))_{t=0}$. With $\Phi(x, \mathfrak{p})$ understood, the function y = y(x, $\mathfrak{p}$) obeys (7.8) with the addition on the right side of the term $\Pi_0(\upsilon\nabla\mathfrak{p}_{\mathfrak{c}(x,\mathfrak{p})})$. This is to say that y is the solution to an equation that has the form

$$y + h(x, y, \mathfrak{p}) = 0 ,$$
(7.15)

where h(x, y, 0) is the function that is depicted in (7.9). Granted (7.15), the map $\pi_r$ near $(\rho_i, 0, \mathfrak{c})$ sends (x, $\mathfrak{p}$) to y = y(x, $\mathfrak{p}$). The map $f'$ can be taken to be $(x, \mathfrak{p}) \to \frac{\partial y}{\partial x}|_{(x,\mathfrak{p})}$, and the map $f''$ can be taken to be $(x, \mathfrak{p}) \to \frac{\partial^2 y}{\partial x^2}|_{(x,\mathfrak{p})}$. For fixed $\mathfrak{p} \in \mathcal{B}'$, let $(x, t) \to y_\mathfrak{p}(x, t)$ denote the function on a neighborhood of the origin in $\mathbb{R}^2$ that given by y(x, t$\mathfrak{p}$).

With the preceding understood, the lemma follows by proving the following:

- *There exists $\mathfrak{p}$ such that $\frac{\partial^2 y_\mathfrak{p}}{\partial x \partial t}|_{(0,0)} \neq 0$.*
- *There exists $\mathfrak{p}$ such that $\frac{\partial^2 y_\mathfrak{p}}{\partial x \partial t}|_{(0,0)} = 0$ and $\frac{\partial^3 y_\mathfrak{p}}{\partial x^2 \partial t}|_{(0,0)} \neq 0$.*

(7.16)

To satisfy the first bullet, it is sufficient to find, given $\varepsilon > 0$, a perturbation $\mathfrak{p}$ with the following properties: For any $\mathfrak{b} \in C^\infty(M; iT^*M \oplus \mathbb{S})$ and $\lambda \in \mathbb{R}$ near zero,

$$\mathfrak{p}(\mathfrak{c} + \lambda\mathfrak{b}) = \lambda^2 \Pi_0 \mathfrak{b} + \mathcal{O}(\varepsilon\lambda^2 + \lambda^3).$$
(7.17)

To argue for the second bullet, it is sufficient to find, given $\varepsilon > 0$, a perturbation $\mathfrak{p}$ such that

$$\mathfrak{p}(\mathfrak{c} + \lambda\mathfrak{b}) = \lambda^3 \Pi_0 \mathfrak{b} + \mathcal{O}(\varepsilon\lambda^3 + \lambda^4) .$$
(7.18)

Note that with these choices, the section $\Phi(x, t\mathfrak{p})$ that solves the version of (7.6) with $\upsilon\nabla\mathfrak{p}$ is respectively $\mathcal{O}(\varepsilon t x)$ and $\mathcal{O}(\varepsilon t x^2)$ for x and t near zero. The fact that such $\mathfrak{p}$ exist in $\mathcal{P}$ follows from the denseness conditions that are stated in Definition 11.6.3 in [KM].

*Part 2*: Choose $\mathfrak{p}_1 \in \mathcal{P}$ very close to zero such that there are but a finite number of critical points of $\pi_r$ on $\mathcal{S}_{*\mathfrak{p}}$, and so that all are non-degenerate. This part explains why there are points $\mathfrak{p}' \in \mathcal{P}$ in any given neighborhood of $\mathfrak{p}_1$ such that the critical points of $\pi_r$ on $\mathcal{S}_{*\mathfrak{p}'}$ are finite in number, non-degenerate, and have distinct critical values. To this end, let $\{(\mathfrak{r}_\sigma, \mathfrak{p}_1, \mathfrak{c}_\sigma)\}_{\sigma=1,2...N}$ label the critical points of $\pi_r$ on $\mathcal{S}_{*\mathfrak{p}_1}$. Let (r, $\mathfrak{p}_1$, $\mathfrak{c}$) denote one of these points. A neighborhood of this point in $\mathcal{S}_*$ is parametrized by a map from a neighborhood of (0, 0) in $\mathbb{R} \times \mathcal{P}$ just like the map introduced in Part 1. To elaborate, this map sends $(x, \mathfrak{q}) \to (\mathfrak{r} + y(x, \mathfrak{p}), \mathfrak{p}_1 + \mathfrak{p}, \mathfrak{c}(x, \mathfrak{p}))$ where the notation is as follows: First, $\mathfrak{c}(x, \mathfrak{p}) = \mathfrak{c} + x\mathfrak{b}_0 + \Phi(x, \mathfrak{p}))$ where $\mathfrak{b}_0$ spans the kernel of the r, $\mathfrak{g} = \mathfrak{e}_\mu + \upsilon\mathfrak{p}_1$, and $\mathfrak{c}$ version of $\mathfrak{L}$, and where $\Phi(x, \mathfrak{p})$ solves the version of (7.6) that has $\upsilon(\nabla\mathfrak{p}_1 + \nabla\mathfrak{p})_{\mathfrak{c}(x,\mathfrak{p})}$ added to its right hand side. Second, y = y(x, $\mathfrak{p}$) now solves (7.7) with the term $\upsilon \Pi_0 (\nabla \mathfrak{p}_1 + \nabla \mathfrak{p})_{\mathfrak{c}(x,\mathfrak{p})}$ added to the right hand side. This equation is equivalent to a version of (7.15) where h vanishes at the origin in $\mathbb{R}^2 \times \mathcal{P}$ as well as its first derivatives in x and y. Meanwhile, its second derivative in x is non-zero at the origin in $\mathbb{R}^2 \times \mathcal{P}$.

Granted all of this, fix one of the indices $\sigma \in \{1, \ldots, N\}$ that label the critical points of $\pi_r$ on $\mathcal{S}_{*\mathfrak{p}_1}$. Suppose that $\lambda > 0$, $\varepsilon > 0$ and that $\mathfrak{p} \in \mathcal{P}$ is such that

$$\mathfrak{p}(\mathfrak{c}_\sigma + \lambda \mathfrak{b}) = \lambda \, \Pi_{0\sigma} \mathfrak{b} + \mathcal{O}(\varepsilon \lambda) \quad \text{and that} \quad \mathfrak{p}(\mathfrak{c}_{\sigma'} + \lambda \mathfrak{b}) = \mathcal{O}(\varepsilon \lambda) \text{ when } \sigma' \neq \sigma.$$
(7.19)

Suppose that $\mathfrak{p}$ is as just depicted, and that $\varepsilon$ is very small. For all sufficiently small $t > 0$, the critical points of $\pi_r$ on $\mathcal{S}_{*\mathfrak{p}_1 + t\mathfrak{p}}$ are in 1-1 correspondence with those of $\pi_r$ on $\mathcal{S}_{*\mathfrak{p}_1}$, and vice versa. The difference between the critical values of the members of each such pair is $\mathcal{O}(t\varepsilon)$ except for the pair with $(r_\sigma, \mathfrak{p}_1, \mathfrak{c}_\sigma)$. The difference here will be $\mathcal{O}(t)$ only. Thus, granted a version of (7.25) for each critical point of $\pi_r$ on $\mathcal{S}_{*\mathfrak{p}_1}$, it then follows that there exists perturbations in any given neighborhood of $\mathfrak{p}_1$ with the property that $\pi_r$ on $\mathcal{S}_{*\mathfrak{p}_2}$ has a finite set of critical points, all non-degenerate, and no two with the same value of $\pi_r$.

As before, the denseness conditions that are stated in Definition 11.6.3 of [KM] guarantee that the required perturbations do indeed exist.

*Part 3*: Now choose $\mathfrak{p} \in \mathcal{P}$ with very small norm so that there are finitely many critical points of $\pi_r$ on $\mathcal{S}_{*\mathfrak{p}}$, all are non-degenerate, and such that the values of $\pi_r$ distinguish these points. This part of the subsection explains how to choose $\mathfrak{p}'$ in any given neighborhood of $\mathfrak{p}$ so that the following is true:

*At a critical value of $\pi_r$, the values of the r and $\mathfrak{g}(r) = \mathfrak{e}_\mu + \upsilon \chi(r) \mathfrak{p}'$ version of (2.9) distinguish the points in $\pi_r^{-1}(r) \subset \mathcal{S}_{*\mathfrak{p}'}$.*
(7.20)

Note that whether or not (7.20) is obeyed, it is still the case that for any $\mathfrak{p}'$ sufficiently close to $\mathfrak{p}$, there are a finite number of critical points of $\pi_r$ on $\mathcal{S}_{*\mathfrak{p}'}$, all are non-degenerate, and they are distinguished by the values of $\pi_r$.

To achieve (7.20), let $r_*$ denote a critical value of $\pi_r$ on $\mathcal{S}_{*\mathfrak{p}}$. Let $(r_*, \mathfrak{p}, \mathfrak{c})$ denote the corresponding critical point. The denseness conditions that are stated in Definition 11.6.3 of [KM] guarantee the existence of an element q in any given neighborhood of 0 in $\mathcal{P}$ such that q takes distinct values on $\pi_r^{-1}(r) \subset \mathcal{S}_{*\mathfrak{p}}$. Fix such an element q. For $\lambda$ sufficiently small, the critical points of $\pi_r$ on $\mathcal{S}_{*\mathfrak{p}+\lambda\mathfrak{q}}$ are in 1-1 correspondence with the critical points of $\pi_r$ on $\mathcal{S}_{*\mathfrak{p}}$. This correspondence is such as to pair critical points that are very much closer to each other then to any other critical points. It then follows that if $\lambda$ is sufficiently small, then the values of the r and $\mathfrak{g}(r) = \mathfrak{e}_\mu + \upsilon(\mathfrak{p} + \lambda\mathfrak{q})$ version of (2.9) distinguish the elements in $\pi_r^{-1}(r) \subset \mathcal{S}_{*\mathfrak{p}+\lambda\mathfrak{q}}$ for values of r near the value of $\pi_r$ on the critical point in $\mathcal{S}_{*\mathfrak{p}+\lambda\mathfrak{q}}$ that is paired with $(r, \mathfrak{p}, \mathfrak{c})$.

*Part 4*: Suppose now that $\mathfrak{p} \in \mathcal{P}$ has very small norm, and is such that $\pi_r$ on $\mathcal{S}_{*\mathfrak{p}}$ has a finite set of critical points, all non-degenerate, all with distinct critical values and such that (7.20) holds. This part explains how to find $\mathfrak{p}'$ in any given neighborhood of $\mathfrak{p}$ so that the following is true:

*The r and $\mathfrak{g}(r) = \mathfrak{e}_\mu + \upsilon \chi(r)\mathfrak{p}'$ version of (2.9) distinguishes the points in $\pi_r^{-1}(r) \subset \mathcal{S}_{*\mathfrak{p}'}$ for all but finitely many values of r, and at the latter, at most one pair of points is not distinguished by this same version of (2.9).*

(7.21)

For this purpose, introduce $\mathcal{S}_* \otimes \mathcal{S}_* \subset \mathcal{S}_* \times \mathcal{S}_*$ to denote the subset that consists of pairs $((r, \mathfrak{p}, \mathfrak{c}_1), (r, \mathfrak{p}, \mathfrak{c}_2))$ with $\mathfrak{c}_1 \neq \mathfrak{c}_2$. Also, introduce $\mathcal{S}_* \otimes \mathcal{S}_* \otimes \mathcal{S}_* \subset \mathcal{S}_* \times \mathcal{S}_* \times \mathcal{S}_*$ to denote the subset of triples $((-, \mathfrak{c}_1), (-, \mathfrak{c}_2), (-, \mathfrak{c}_3))$ where no two of the three are the same. The first space is a manifold at points where both $\mathfrak{c}_1$ are $\mathfrak{c}_2$ non-degenerate solutions to the r and $\mathfrak{g}(r, \mathfrak{p}) = \mathfrak{e}_\mu + \upsilon \chi(r)\mathfrak{p}$ version of (2.4), and the second is a manifold at points where all three are non-degenerate solution to this version of (2.4). By virtue of the choice of $\mathfrak{p}$, it is only necessary to consider the parts of these spaces where such is the case.

Let $\mathfrak{a}_{\mathfrak{g}(r,\mathfrak{p})}$ denote the version of (2.9) that is defined using r and $\mathfrak{g}(r, \mathfrak{p})$. Now consider the functions $\mathfrak{h}: \mathcal{S}_* \otimes \mathcal{S}_* \to \mathbb{R}$ that assigns $\mathfrak{a}_{\mathfrak{g}(r,\mathfrak{p})}(\mathfrak{c}_1) - \mathfrak{a}_{\mathfrak{g}(r,\mathfrak{p})}(\mathfrak{c}_2)$ to given $((r, \mathfrak{p}, \mathfrak{c}_1), (r, \mathfrak{p}, \mathfrak{c}_2))$. Likewise, define $\mathfrak{h}_2: \mathcal{S}_* \otimes \mathcal{S}_* \otimes \mathcal{S}_* \to \mathbb{R}^2$ by declaring that its first component be $\mathfrak{h}((-, \mathfrak{c}_1), (-, \mathfrak{c}_2))$ and that its second be $\mathfrak{h}((-, \mathfrak{c}_1), (-, \mathfrak{c}_3))$. The first point to make is that both $\mathfrak{h}$ and $\mathfrak{h}_2$ are submersions at all points in $\mathfrak{h}^{-1}(0)$ and $\mathfrak{h}_2^{-1}(0)$, respectively. This again follows from the from the denseness conditions that are stated in Definition 11.6.3 in [KM]. The point is that one can find some $\mathfrak{q}$ in any given neighborhood of 0 in $\mathcal{P}$ that distinguishes any three elements in $\text{Conn}(E) \times C^\infty(M; \mathbb{S})$. To continue, $\mathfrak{h}^{-1}(0)$ is a smooth, codimension 1 submanifold of $\mathcal{S}_* \otimes \mathcal{S}_*$. The projection, $\pi: \mathfrak{h}^{-1}(0) \to \mathcal{B}$ is a Fredholm map of index zero, and so it has a residual set of regular values. If $\mathfrak{p}'$ is such a regular value, then there are at most a finite set of points where $\mathfrak{h}^{-1}(0)$ intersects $\mathcal{S}_{*\mathfrak{p}'} \otimes \mathcal{S}_{*\mathfrak{p}'}$. Meanwhile, $\mathfrak{h}_2^{-1}(0)$ is a smooth codimension 2 submanifold of $\mathcal{S}_* \otimes \mathcal{S}_* \otimes \mathcal{S}_*$ and the restriction of $\pi$ to $\mathfrak{h}^{-1}(0)$ is a Fredholm map of index -1. It too has a residual set of regular values. If $\mathfrak{p}'$ is in both of these residual sets and close to $\mathfrak{p}$, then both (7.20) and (7.21) are satisfied.

*Part 5*: This part of the subsection addresses the claims in Part 1 of Property 5. Here, the story is really no different than what has been done so far. To elaborate, let $\{\mathfrak{c}_v\}_{v=1,2,...}$ label the $I_+$-canonical basis for the cSWF complex at $y = y_+$. Extend this basis to $y = y_-$ using the maps in Item 2 of the first bullet of (3.5). If $\mathfrak{c}$ and $\mathfrak{c}'$ have the same degree, then the argument for Property 3 works in this case if it is understood that the basis used at $y_-$ is the extension via the maps in Item 2 of the $I_+$-canonical basis used at $y_+$. The point here is that there are no instantons from $\mathfrak{c}$ to $\mathfrak{c}'$ or vice versa for values of r near y, and so whether their ordering changes or not in their guise as $I_\pm$-canonical basis elements when r crosses y makes no difference to the differential of the cSWF complex as long as it is understood that the basis used is <u>not</u> changed as r crosses y. The assertion in Part 1 of Property 5 in this case follows directly from this observation.

The argument for Property 3 also works with no change if the degrees of $\mathfrak{c}$ and $\mathfrak{c}'$ do not differ by 1. This understood, consider the case where degree($\mathfrak{c}$) = degree($\mathfrak{c}'$) + 1. To start, note that there are no instantons from $\mathfrak{c}$ to $\mathfrak{c}'$ at r = y since $\mathfrak{a}_{\mathfrak{g}(y)}$ has the same value on these two generators. This implies that there are no instantons between $\mathfrak{c}(r)$ and $\mathfrak{c}'(r)$ for r near y as well. Indeed, were there such instantons for a sequence $\{r_n\}_{n=1,2,\ldots}$ converging to y, essentially the same argument used by [KM] to prove their Theorem 16.1.3 would find a broken trajectory limit of this sequence. In this case, the broken trajectory consists of a sequence $\{\mathfrak{d}_1, \ldots, \mathfrak{d}_n\}$ of instanton solutions to the r = y and $\mathfrak{g}(y)$ version of (2.11) such that the s → -∞ limit of $\mathfrak{d}_1$ is $\mathfrak{c}$, the s → ∞ limit of $\mathfrak{d}_n$ is $\mathfrak{c}'$, and such that for each j ∈ {2, …, n}, the s → -∞ limit of $\mathfrak{d}_j$ is the s → ∞ limit of $\mathfrak{d}_{j-1}$. Moreover, at least one $\mathfrak{d}_j$ in this sequence is not constant as s varies in $\mathbb{R}$. But this is impossible as the drop in the y and $\mathfrak{g}(y)$ version of (2.9) along any non-constant instanton is negative. Moreover, the sum of these drop is $\mathfrak{a}_{\mathfrak{g}(y)}(\mathfrak{c}') - \mathfrak{a}_{\mathfrak{g}(y)}(\mathfrak{c})$, and since this number is zero, there are no such broken trajectories.

Granted that there are no instantons between $\mathfrak{c}(r)$ and $\mathfrak{c}'(r)$ for r near y, it then follows that the differential in the cSWF complex is insensitive to the fact that $\mathfrak{a}_{\mathfrak{g}(r)}$ takes equal values on $\mathfrak{c}(r)$ and $\mathfrak{c}'(r)$ at r = y.

*Part 6*: This last part of the subsection considers the assertions made in Part 2 and Part 3 of Property 5. The task here is to prove Lemmas 3.15 and 3.17. To start this task, fix σ > 0 be fixed such that y + σ < $y_+$ and y - σ > $y_-$. Fix a smooth, increasing function on $\mathbb{R}$ with derivative bounded by 1 that equals y - σ where s < -1 and y + σ where s > 1. Denote this function by $\mathfrak{r}$. Fix $\mathfrak{p}' \in \mathcal{P}$ with very small norm. There is a residual set of choices for $\mathfrak{p}'$ such that at both r = y+σ and r = y-σ, all instanton solutions to the r and $\mathfrak{g} = \mathfrak{e}_\mu + \mathfrak{p} + \mathfrak{p}'$ version of (2.11) that limit as s → ±∞ to degree k or greater solutions to this r and $\mathfrak{g}$ of (2.4) have non-degenerate moduli spaces. Take $\mathfrak{p}'$ much closer to 0 then the version of $\mathfrak{p}_m$ supplied by Property 2 for the interval in $I_-$ that contains r = y - σ and also the version for the interval in $I_+$ that contains y + σ.

Consider the equations for a map s → $\mathfrak{d}(s)$ = (A(s), ψ(s)) given by

- $\frac{\partial}{\partial s} A = -B_A + \mathfrak{r}(s)(\psi^\dagger \tau^k \psi - i a) + \upsilon \chi(\mathfrak{r}(s)) \mathfrak{T}(A, \psi)$,
- $\frac{\partial}{\partial s} \psi = -D_A \psi + \upsilon \chi(\mathfrak{r}(s)) \mathfrak{S}(A, \psi)$.

(7.23)

Here, $\mathfrak{T}$ and $\mathfrak{S}$ are the respective components of $\nabla(\mathfrak{p} + \mathfrak{p}')$ in $C^\infty(M; iT^*M)$ and $C^\infty(M; \mathbb{S})$. Of particular interest are the instanton solutions, those that limit as $s \to -\infty$ to a solution of the $r = y - \sigma$ and $\mathfrak{g} = \mathfrak{e}_\mu + \upsilon\chi(y-\sigma)(\mathfrak{p} + \mathfrak{p}')$ version of (2.4), and limit as $s \to \infty$ to a solution of the $r = y + \sigma$ and $\mathfrak{g} = \mathfrak{e}_\mu + \upsilon\chi(y+\sigma)(\mathfrak{p} + \mathfrak{p}')$ version of (2.4). Let $\mathfrak{c}_-$ denote a solution to the former version of (2.4) and let $\mathfrak{c}_+$ denote a solution to the latter. Let $\mathcal{M}_{y,\sigma}(\mathfrak{c}_-, \mathfrak{c}_+)$ denote the space of solutions to (7.23) with $s \to -\infty$ limit $\mathfrak{c}_-$ and $s \to \infty$ limit $u\mathfrak{c}_+$ where u can be any smooth map from M to $S^1$. As with the case of (7.12), Proposition 24.4.7 of [KM] finds a residual set of choices for $\mathfrak{p}'$ from the ball of radius 1 about the origin in $\mathcal{P}$ such that the following is true: If $\mathfrak{c}_-$ has degree $d_- \geq k$ and $\mathfrak{c}_+$ has degree $d_+ \geq k$, then the moduli space $\mathcal{M}_{y,\sigma}(\mathfrak{c}_-, \mathfrak{c}_+)$ has the structure of a smooth, manifold whose dimension is $d_- - d_+$. Assume that $\mathfrak{p}'$ is now from this residual set. In the case when $d_- = d_+$, it follows from Theorem 24.6.2 in [KM] that $\mathcal{M}_{y,\sigma}(\mathfrak{c}_-, \mathfrak{c}_+)$ is compact.

Assume now that $\mathfrak{c}_-$ and $\mathfrak{c}_+$ have the same degree. Just as in Part 1 of the proof of Proposition 3.12, each element in $\mathcal{M}_{y,\sigma}(\mathfrak{c}_-, \mathfrak{c}_+)$ has an associated sign, either +1 or -1. This sign is explained in Chapter 25.2 of [KM]. Use $\sigma(\mathfrak{c}_-, \mathfrak{c}_+)$ to denote the sum of these signs with the understanding that $\sigma(\mathfrak{c}_-, \mathfrak{c}_+) = 0$ when $\mathcal{M}_{y,s}(\mathfrak{c}_-, \mathfrak{c}_+) = \emptyset$.

Let cSWF$_-$ denote the cSWF complex in degrees k and greater as defined using the $r = y - \sigma$ and $\mathfrak{g} = \mathfrak{e}_\mu + \upsilon\chi(y-\sigma)(\mathfrak{p} + \mathfrak{p}')$ versions of (2.4) and (2.11) to obtain the generators and differential. Likewise, define cSWF$_+$ using the $r = y + \sigma$ and $\mathfrak{g} = \mathfrak{e}_\mu + \upsilon\chi(y+\sigma)(\mathfrak{p} + \mathfrak{p}')$ versions of (2.4) and (2.11) to obtain the generators and differential. Note that by virtue of the fact that $\mathfrak{p}'$ has very small norm, the $I_-$-canonical basis can be used for the cSWF$_-$ complex and the $I_+$-canonical basis can be used for the cSWF$_+$ complex. Note also that the ordering of the generators in the canonical basis for cSWF$_-$ is the same as that given by the values of the $r = y - \sigma$ and $\mathfrak{g} = \mathfrak{e}_\mu + \upsilon\chi(y-\sigma)(\mathfrak{p} + \mathfrak{p}')$. The analogous statement holds for the ordering of the generators in the canonical basis for the cSWF$_+$ complex.

These integer weights $\{\sigma(\mathfrak{c}_-, \mathfrak{c}_+)\}$ are used, as in (7.13), to define a degree preserving homomorphism from cSWF$_+$ to cSWF$_-$. Chapter 25.3 in [KM] proves that $\mathbb{T}_\sigma$ intertwines the differential on the cSWF$_-$ complex with that on the cSWF$_+$ complex, and induces an isomorphism between the respective homology groups.

*Proof of Lemma 3.15*: In the notation used by Lemma 3.15, the cSWF vectors spaces in any given degree k or greater as defined at $y + \sigma$ is denoted by $\mathbb{V}_+$. The canonical basis of $\mathbb{V}_+$ in any given degree is denoted by $\{\mathfrak{c}_v\}$. With this notation,

- $\mathbb{T}_\sigma \mathfrak{c}_v = \sum_{v'} \sigma(\mathfrak{c}_v, \mathfrak{c}_{v'}) \mathfrak{c}_{v'}$ *in degrees not equal to d or d+1*.
- $\mathbb{T}_\sigma \mathfrak{c}_v = \sum_{v'} \sigma(\mathfrak{c}_v, \mathfrak{c}_{v'}) \mathfrak{c}_{v'}$ *and* $\mathbb{T}_\sigma \mathfrak{c} = \sum_{v'} \sigma(\mathfrak{c}, \mathfrak{c}_{v'}) \mathfrak{c}_{v'}$ *in degree d+1*.
- $\mathbb{T}_\sigma \mathfrak{c}_v = \sum_{v'} \sigma(\mathfrak{c}_v, \mathfrak{c}_{v'}) \mathfrak{c}_{v'}$ *and* $\mathbb{T}_\sigma \mathfrak{c}' = \sum_{v'} \sigma(\mathfrak{c}', \mathfrak{c}_{v'}) \mathfrak{c}_{v'}$ *in degree d*.

(7.24)

Let $\mathbb{A}\colon \mathbb{V}_+ \to \mathbb{V}_+$ denote the restriction of $\mathbb{T}_\sigma$ to the $\mathbb{V}_+$ summand in $\mathbb{Z}\mathfrak{c} \oplus \mathbb{Z}\mathfrak{c}' \oplus \mathbb{V}_+$. An argument that differs only cosmetically from that used in the proof of Proposition 3.12 proves that if $\sigma$ is sufficiently small, then $\mathbb{A}$ is an upper triangular matrix with 1's on the diagonal. These arguments also prove that $\mathcal{M}_{y,\sigma}(\mathfrak{c}, \mathfrak{c}_v)$ and $\mathcal{M}_{y,\sigma}(\mathfrak{c}', \mathfrak{c}_v)$ are empty unless $\mathfrak{a}_{g(y)}(\mathfrak{c}_v) < \mathfrak{a}_{g(y)}(\mathfrak{c}) = \mathfrak{a}_{g(y)}(\mathfrak{c}')$. Here, $\mathfrak{a}_{g(y)}$ is the $r = y$ and $g(y) = \mathfrak{e}_\mu + \upsilon\chi(y)\mathfrak{p}$ version of (2.9). Note in this regard that $\mathfrak{a}_{g(y)}(\mathfrak{c})$ is not equal to any $\mathfrak{a}_{g(y)}(\mathfrak{c}_v)$. What has just been said implies that the matrix $\mathbb{T}_\sigma$ satisfies the conditions stated for $\mathbb{T}$ by Lemma 3.15.

The matrix $\mathbb{T}$ is not necessarily equal to $\mathbb{T}_\sigma$. However, as explained next, $\mathbb{T}$ is obtained from $\mathbb{T}_\sigma$ by composing with an upper triangular matrix that has 1's on the diagonal. If this is the case, $\mathbb{T}$ also satisfies the conditions that are stated by Lemma 3.15.

To obtain $\mathbb{T}$ from $\mathbb{T}_\sigma$, let $m \in \mathbb{K}(I)$ be such that $y - \sigma \in [w_m, w_{m+1})$. The cSWF homology in degrees greater than k is defined using $\mathfrak{g} = \mathfrak{e}_\mu + \upsilon\chi(y-\sigma)(\mathfrak{p} + \mathfrak{p}_m)$ to define the generators and differential. On the other hand, the cSWF$_-$ generators and differential are defined using $r = y - \sigma$ and $\mathfrak{g} = \mathfrak{e}_\mu + \upsilon\chi(y-\sigma)(\mathfrak{p} + \mathfrak{p}')$. However, an argument just like that used to prove Proposition 3.12 finds an upper triangular matrix with 1's on the diagonal that maps the first version of the complex to the second, intertwines their differentials, and induces an isomorphism on homology. There is a completely analogous story to be told at $r = y + \sigma$. Composing these matrices with $\mathbb{T}_\sigma$ gives a new matrix, $\mathbb{T}_\sigma'$, that relates the cSWF complex in degrees k and greater at $r = y - \sigma$ to the cSWF complex at $r = y + \sigma$, and that satisfies the conditions stated by Lemma 3.15. The matrix $\mathbb{T}$ is obtained from $\mathbb{T}_\sigma'$ by composing with the upper triangular matrices that are given in Lemma 3.14 to move from $r = y - \sigma$ to $y_-$ and to move from $r = y + \sigma$ to $y_+$. Granted that the matrices from Lemma 3.14 are upper triangular with 1's on the diagonal, Lemma 3.15 follows from what is said in the preceding paragraph.

*Proof of Lemma 3.17*: The notation used here uses $\mathbb{V}_-$ to denote the cSWF vector spaces in any given degree k or greater as defined at $y - \sigma$. The canonical basis of $\mathbb{V}_-$ in any given degree is denoted by $\{\mathfrak{c}_v\}$. With this notation,

- $\mathbb{T}_\sigma \mathfrak{c}_v = \sum_{v'} \sigma(\mathfrak{c}_v, \mathfrak{c}_{v'}) \mathfrak{c}_{v'}$ *in degrees not equal to d+1 or d.*
- $\mathbb{T}_\sigma \mathfrak{c}_v = \sum_{v'} \sigma(\mathfrak{c}_v, \mathfrak{c}_{v'}) \mathfrak{c}_{v'} + \sigma(\mathfrak{c}_v, \mathfrak{c})\mathfrak{c}$ *in degree d+1.*
- $\mathbb{T}_\sigma \mathfrak{c}_v = \sum_{v'} \sigma(\mathfrak{c}_v, \mathfrak{c}_{v'}) \mathfrak{c}_{v'} + \sigma(\mathfrak{c}_v, \mathfrak{c}')\mathfrak{c}'$ *in degree d.*

(7.25)

Write $\mathbb{A}\colon \mathbb{V}_- \to \mathbb{V}_-$ for the composition of $\mathbb{T}_\sigma$ with the projection from $\mathbb{Z}\mathfrak{c} \oplus \mathbb{Z}\mathfrak{c}' \oplus \mathbb{V}_-$ to $\mathbb{V}_-$. As in the preceding proof, arguments that differ only cosmetically from those used to prove Proposition 3.12 prove that $\mathbb{A}$ is upper triangular with 1's on the diagonal if $\sigma$ is sufficiently small. These same arguments show that $\mathcal{M}_{y,\sigma}(\mathfrak{c}_v, \mathfrak{c})$ and $\mathcal{M}_{y,\sigma}(\mathfrak{c}_v, \mathfrak{c}')$ are empty when $\mathfrak{a}_{\mathfrak{g}(y)}(\mathfrak{c}_v) < \mathfrak{a}_{\mathfrak{g}(y)}(\mathfrak{c}) = \mathfrak{a}_{\mathfrak{g}(y)}(\mathfrak{c}')$. As a consequence, the matrix $\mathbb{T}_\sigma$ satisfies the conditions for $\mathbb{T}$ stated by Lemma 3.17.

As in the previous proof, the matrix $\mathbb{T}$ is obtained from $\mathbb{T}_\sigma$ by composing with upper triangular matrices. Thus, $\mathbb{T}$ also satisfies the conditions stated by Lemma 3.17.

## 8. The proof of Theorem 1

The last section puts all of the pieces together and so completes the proof of Theorem 1. To start, fix a complex line bundle $E \to M$ whose first Chern class differs by a torsion class from half the first Chern class of K. Fix a co-exact 1-form μ from the collection supplied by Proposition 3.11. Fix $k \ll 0$ and use μ to define the cSWF homology in degrees greater than k. Fix a non-zero cSWF homology class, θ, with degree $k' > k$ but with $k' < 0$. These are supplied by Proposition 3.8. Section 4 explains how to define θ for all $r > r_k$ save for a discrete set with no accumulation points. It follows from Propositions 4.6 and 5.1 that θ is not a divergence class. As a consequence, there exists an unbounded sequence, $\{(r_n, (A_n, \psi_n))\}_{n=1,2...}$ such that $(A_n, \psi_n)$ satisfies the r = $r_n$ and μ version of (2.5); and such that $(A_n, \psi_n)$ is non-degenerate and has degree k´.

Write $\psi_n = (\alpha_n, \beta_n)$ to correspond with the splitting in (2.2). If E is not the trivial bundle, $1_\mathbb{C}$, then $\alpha_n$ must vanish at some points in M, and so $\sup_M(1 - |\psi_n|) = 1$. As a consequence, all of the conditions in Theorem 2.1 are met, and Theorem 2.1 thus supplies the set of closed integral curves of the Reeb vector field for Theorem 1. Now suppose that $E = 1_\mathbb{C}$. As is explained momentarily, there is in this case a constant $\kappa > 0$ such that $\sup_M(1 - |\psi_n|) > \kappa$. Granted this, Theorem 2.1 again supplies a set of closed integral curves of the Reeb vector field for Theorem 1.

Suppose, for the sake of argument that no such κ exists. The following is then a consequence: Given $\epsilon > 0$, then for all n sufficiently large, the first bullet and the third bullet in (5.24) are satisfied by $(A_n, \psi_n)$. A repeat of the rescaling argument used in Section 6b to prove Lemma 2.3 can be used to establish the second bullet in (5.24). This is because the rescaled sequence of solutions will converge strongly in the ball where $|y| \le 4$ to the solution with $\alpha = 1$ and $\beta = 0$. If $\epsilon < \epsilon_0$ from Lemma 5.4, it then follows that the degree $k' = 0$. This is nonsense since k´ was chosen to be negative.